# Outer Billiards on Regular Polygons

G.H. Hughes

In [M1] (1973), J. Moser proposed that his Twist Theorem could be used to show that orbits of the outer billiards map on a sufficiently smooth closed curve were always bounded. He asked the same question for a convex polygon in [M2] (1978). In [VS] (1987) F.Vivaldi and A. Shaidenko showed that all orbits for a regular polygon must be bounded. More recently, R. Schwartz [S1] showed that a quadrilateral known as a Penrose Kite has unbounded orbits and he proposed that 'most' convex polygons support unbounded orbits.

Except for a few special cases, very little is known about the dynamics of the outer billiards map on regular polygons. In this paper we present a unified approach to the analysis of regular polygons – using the canonical 'resonances' which are shared by all regular N-gons. In the case of the regular pentagon and regular octagon these resonances exist on all scales and the fractal structure is well documented, but these are the only non-trivial cases that have been analyzed. We present a partial analysis of the regular heptagon, but the limiting structure is poorly understood and this does not bode well for the remaining regular polygons. The minimal polynomial for the vertices of a regular N-gon has degree $\varphi(N)/2$ where $\varphi$ is the Euler totient function, so N = 5, 7 and 11 are respectively quadratic, cubic and quintic. In the words of R. Schwartz, "A case such as N = 11 seems beyond the reach of current technology."

The study of regular N-gons can be reduced to the case where N is even and this implies that the dynamics of the outer billiards map can be reduced to piecewise rotations on a rhombus which defines the 'first generation'- and all subsequent generations. Thus the outer billiards map for regular polygons can be regarded as a piecewise affine map on a torus. One such map is the Digital Filter map (Df) described in Appendix F. The linear form of Df is conjugate to a rotation with 'twist' $\rho$ in the range (0, ¼]. Every rational $\rho$ = p/q in this range with p and q relatively prime integers, corresponds to a regular q-gon with Df web of 'step-size' p if q is even or a regular 2q-gon with step size 2p if q is odd. This implies that the rotational dynamics of a regular N-gon may allow for multiple interpretations – where traditional outer billiards is just the step-1 version.

Piecewise affine maps have been studied by a number of authors including L.Chua, R. Adler, A. Goetz, J. Lowenstein, and F. Vivaldi. Recent studies are based on what Vivaldi calls 'algebraic dynamics'- using algebraic number theory and symbolic computations to analyze the dynamics. Among the major issues are the Lebesgue measure and Hausdorff dimension of the 'residual set' $\overline{W}/W$ where W is the singularity set. This set contains all the non-trivial dynamics and typically has empty interior and (maximal) fractal dimension between 1 and 2. These methods have been applied successfully to the dynamics of quadratic cases, and recently there has been progress in the cubic case. In [LKV] (2005) and [L] (2007) Lowenstein et al. used algebraic techniques and Mathematica to obtain an accurate spectrum for the recurrence time dimensions for a map based on rotations by $\pi/7$. Applying these techniques to the cubic regular polygons N = 7 and N = 9 will be a challenge because the scaling is not well understood.

Section 1 - Introduction

Below are the first few iterations of the outer billiards map (Tangent map) $\tau$ for two convex polygons. For an initial point p, each iteration is a reflection (central symmetry) about one of two possible 'support' vertices, so the formula is $\tau(p) = 2c - p$ where c is the support vertex. We will usually assume a clockwise orientation for the original polygon. The inverse map $\tau^{-1}$ is $\tau$ applied to the polygon with opposite orientation, so for a regular polygon, $\tau$ is essentially its own inverse.

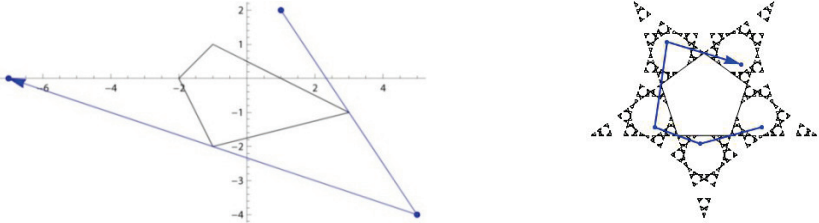

Iterating $\tau$ with regular polygons, generates canonical 'First Families' of related polygons. These families represent the major resonances of $\tau$, so they have periodic orbits with relatively small periods. Some of these resonances can be seen above for the regular pentagon – also known as N = 5 or $M_5$. The blue orbit shown here is period 10 but the center of the decagon 'tile' is period 5 because $\tau$ inverts the decagon on each iteration- leaving just the center fixed. This tile is called S[1] because the orbit advances one vertex on each iteration.

Below are the canonical resonances for N = 5, 7 and 14. Note that N = 7 and N = 14 have essentially the same families. This occurs because N = 7 generates a tile which is a perfect copy of N = 14 and N = 14 generates a tile which is a perfect copy of N = 7, so these two have 'conjugate' dynamics. This is called the Twice-Odd Lemma. See Section 2.

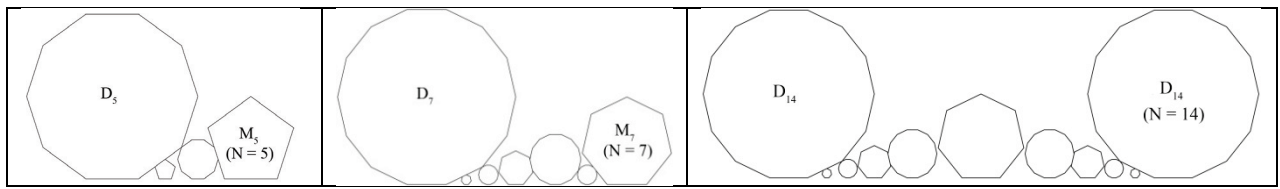

**Notation:** (i) When k is odd, the generating regular k-gon will be known as $M_k$ (or N = k or simply M, when k is understood) . In this case there will always be a matching regular 2k-gon with the same side length. This tile will be called $D_k$ or simply D, so M and D form the nucleus of the First Family when k is odd

(ii) When k is even, the generating regular k-gon will be known as $D_k$ (or N = k or simply D when k is understood) and there will be a matching tile which is identical to $D_k$. These two will form the nucleus of the First Family for k even.

Therefore every regular polygon generates a matching D tile and these tiles play a major role in the dynamics. They have maximal measure and serve as bounds for the dynamics.

For any generating polygon P, the resonant tiles on all scales can be generated by iterating the extended 'trailing' edges of P using the inverse map $\tau^{-1}$. This is called the (forward) web. The web can also be obtained by iterating the extended forward edges under $\tau$. This is called the inverse web.

**Definition of the web for a convex polygon P**

Given a convex n-gon P with vertices $\{c_1, c_2,..,c_n\}$ which we assume to be numbered clockwise, the corresponding edges $E = \{E_1, E_2, .., E_n\}$ have both a 'forward' and 'trailing' extension. (Just the trailing extensions are shown below in blue.) Let $E^f = \cup E^f_k$ be the union of the forward extensions and $E^t = \cup E^t_k$ be the union of the trailing extensions. Each extension is assumed to be an open ray.

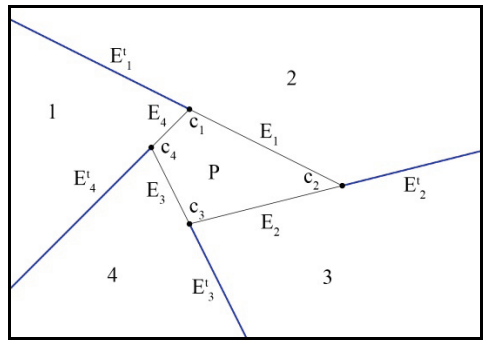

The map $\tau$ is not defined on P or on the trailing edges, so the level-0 web is defined to be $W_0 = E \cup E^t$. This is the exceptional (or singular) set of $\tau$. Since $W_0$ is connected, the complement of $W_0$ external to P consists of n disjoint open (convex) sets which are known as level-0 tiles. These tiles are also called 'atoms' because all the subsequent dynamics are determine by repeated action of $\tau$ on this level-0 partition of the space external to P. These primitive tiles define the domain of $\tau$ relative to each vertex so they can be labeled by the indices 1,2,..n. These labels will be the first elements of the 'corner sequence' of a point in these regions. As the web progresses, each level-k tile will have a corner sequence of length k+1.

The union of the level-0 tiles is the domain of $\tau$ which is abbreviated Dom($\tau$). Dom($\tau^2$) is Dom($\tau$) – $\tau^{-1}$ ($W_0$). The union of $W_0$ and $\tau^{-1}$ ($W_0$) is called the level-1 web, $W_1$. In general the level k (forward) web is defined to be:

$$(i) \quad W_k = \bigcup_{j=0}^{j=k} \tau^{-j}(W_0) \text{ where } W_0 = E \cup E^t$$

The inverse web is defined in a similar fashion using $\tau$ and the extended forward edges:

$$(ii) \quad W^i_k = \bigcup_{j=0}^{j=k} \tau^{j}(W^i_0) \text{ where } W^i_0 = E \cup E^f$$

At each iteration, $W_k$ and $W^i_k$ partition the plane into disjoint open convex regions which are the level-k tiles. Below is $W_1$ for the polygon P, showing the level-1 tiles and the corner sequences.

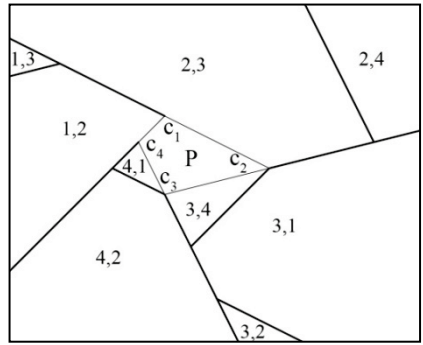

The analysis of these corner sequences is part of symbolic dynamics. For any convex piecewise isometry such as τ, these corner sequences have polynomial complexity [B] and the degree of that polynomial is a common measure of the complexity of the dynamics. The first differences of these corner sequences yields 'step-sequences' which in turn can be used to define winding numbers - the average rotation per iteration. See Appendices C and E.

$W_k$ is the (forward) exceptional set for $\tau^k$ and the limiting web W is $\lim_{k \to +\infty} W_k$. The limiting inverse web $W^i$ is defined in a similar fashion by applying τ to the forward edges. Clearly these two limiting webs are identical but they differ on each iteration and it is useful to make the distinction.

Since W is the union of a countable number of lines or line segments, it has zero (Lebsegue) measure so the complement $W^c$ has full measure and τ is defined 'almost' everywhere.

At every iteration of the web there are unbounded tiles but for regular polygons the limiting tiles have bounded measure and the tiles with this maximal measure are the D's. These D tiles also have the largest number of possible sides which is 2N for N odd and N for N even. No tile can have more than 2N sides because all iterations of an extended edge are parallel, so there are never more than 2N directions to choose from. For N even, there are just N directions.

For regular polygons most of the tiles are 'stable' at some finite iteration k – which means they are no longer partitioned by the web W. These stable tiles must be periodic (see below). There may be bounded tiles which are not stable at any iteration. These tiles will have zero limiting measure, so they are points or line segments. These limiting tiles can only arise only if the distance between iterations of the web goes to 0 and these distances are linear combinations of the vertices of the N-gon, so the linear space generated by $\cos(2\pi/N)$ must have a rank at least 2. For regular N-gons, this rank is $\varphi(N)/2$ where φ is the Euler totient function. (We will define this to be the algebraic complexity of the N-gon.) The regular polygons N = 3,4 and 6 have linear complexity so they are affinely equivalent to rational polygons. N = 5 is the first regular polygon with quadratic complexity.

The web for the regular pentagon is partitioned in a self-similar fashion which yields generations of D tiles (and M tiles) on all scales as shown on the left below. Therefore the limiting tiles are points whose orbits are not finite. Since the boundaries of tiles are always in W, these limit points are in the closure $\overline{W}$. The Hausdorff dimension of $\overline{W}/W$ for N = 5 is derived in section 4. Even though there are an uncountable number of points in this 'residual' set, the Lebesgue measure is zero - so the periodic obits have full measure. This does not rule out the possibility that the residual set $\overline{W}/W$ may have non-zero measure for other polygons. This issue is still unresolved and has a long history dating back to the issue of positive measure for chaotic orbits of Hamiltonian systems. The image shown here on the left is generated from a single nonperiodic orbit - which is dense in $\overline{W}/W$. The coding for this orbit is given in Appendix C.

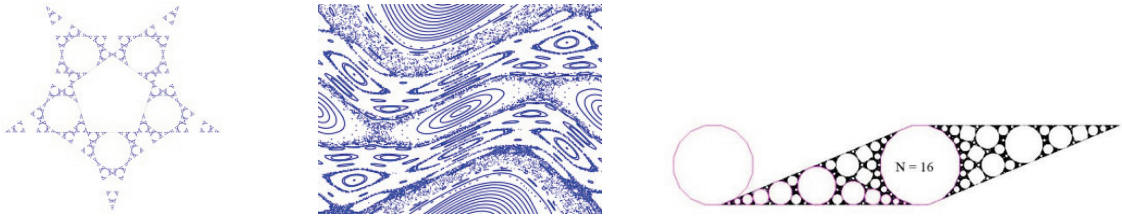

For continuous systems arising in Hamiltonian dynamics, the resonances are 'center manifolds' and these manifolds are typically separated by 'stochastic' webs as shown above using the Standard Map. The Tangent map for a convex polygon is a discontinuous model of Hamiltonian dynamics and therefore it can be used to study the breakdown of the KAM stability. There is a toral version of $\tau$, which plays a role similar to the Standard Map. This $\tau$-torus is illustrated on the right above with N = 16, where the First Family is shown in magenta. See Appendix F.

These webs can be implemented without explicit use of $\tau^{-1}$, because this map can be obtained from $\tau$ by reversing the orientation of the generating polygon. Assuming that the origin is inside the polygon, the orientation can be reversed by taking the reflection about the y-axis – which we call Tr.

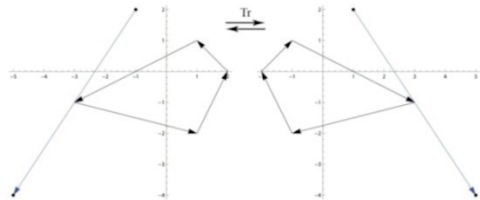

Tr is an affine transformation and as such it commutes with $\tau$ – which is a collection of affine transformations. Two polygons which are affinely related will have conjugate dynamics and they are said to be in the same affine family. Any discussion of the dynamics of a polygon will apply equally to all members of the same affine family. For example no regular triangle or hexagon has rational coordinates, but they are affinely equivalent to polygons with rational coordinates, so they are classified as 'rational'. The only true rational regular polygon is the square. As far as $\tau$ is concerned, N = 3, 4 and 6 are all lattice polygons and their web geometry is trivial.

Our arbitrary choice is to implement τ clockwise, but for a counterclockwise matrix M, it is just a matter of reflecting M to get Tr[M], then using the clockwise τ, and reflecting the result back as shown above. For regular polygons centered at the origin, M and Tr[M] are identical (except for orientation) so the procedure for generating the web is very simple: implement τ for a clockwise M and apply τ to the forward edges to obtain the inverse web $W^i_k$. Sometimes this inverse web is sufficient, but if a true web is needed, just apply Tr to the resulting inverse web to get $W_k$. This works because there is no loss of generality in assuming that the original matrix was counterclockwise. Since $W_k$ and $W^i_k$ are conjugate, any analysis can be done with either web.

**Example**: On the left below are the two level-5 webs in the vicinity of S[1] for N = 11. The inverse web $W^i_5$ is in blue and the forward web $W_5$ is in Red. Because these webs are just reflections of each other, using both webs is computationally very efficient. S[1] and S[2] get their names from the fact that the interior points have step sequences which are constant {1} and {2}. These tiles evolve at different rates but both will be stable 22-gons after 11 iterations of the combined web. On the right is the level 200 combined web.

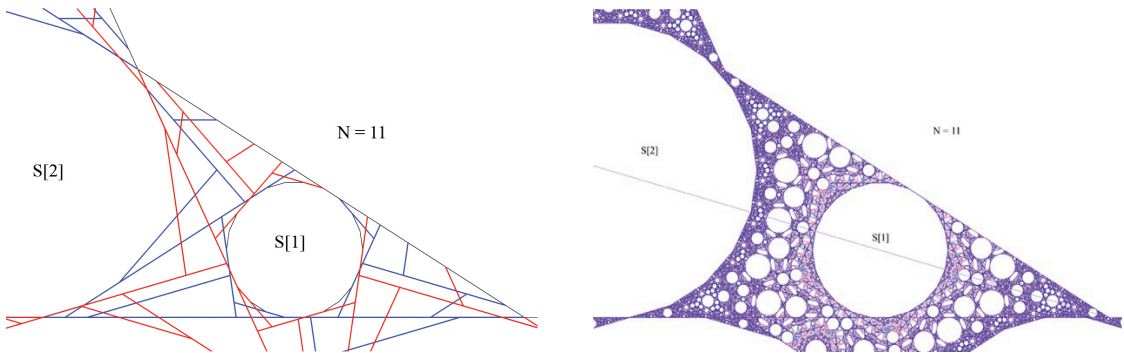

Another way to compare webs is to iterate the forward edges and see how long it takes for them to map to the trailing edges. This is done below for the regular pentagon, N = 5. The forward edges are in blue and the trailing edges are in magenta. As the blue points are mapped under τ to obtain the inverse web, they cannot be periodic because they have no inverse image, so they either map to trailing edges or are nonperiodic. (There needs to be a distinction here between these nonperiodic points and the 'boundary' points which map to the trailing edges because they were not in Dom($\tau^k$) for some k.) Shown here are levels 0 through 6 followed by level 20 of the inverse web $W^i$. The $L_1$ segments generate the 'star' region bounded by the ring of 5 D's.

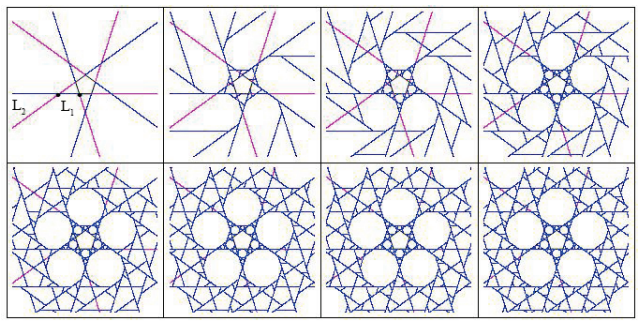

The $L_2$ rays generate the outer edges of the five D's and endless rings of congruent D's at horizontal spacing $2(|L_1|+ s)$ where s is the edge length of N = 5. The periods of these rings are odd multiples of 5.

Because N = 5 has fairly simple dynamics these two webs align very quickly, but there will always be magenta 'hot' spots because of the nonperiodic orbits. The rate at which these two webs align can be used as a measure of the complexity of the dynamics in a given region.

As indicated earlier, a tile is said to be 'stable' at level-k if the tile is no longer partitioned after level-k of the web. Stable tiles cannot have zero measure so when the generating polygon is regular, the orbit of a stable tile is bounded and not self-intersecting so it must be periodic.

If a tile has a periodic orbit, the points in the orbit of this tile will have the same step sequence (except for cyclic rotation), but this does not imply that all points will have the same prime period. On each iteration of $\tau$, the tiles are inverted, so periods will tend to be even, but it is possible for a center of symmetry to have an odd period. This can only occur if there are an odd number of tiles in the orbit. A tile does not have to be regular to experience this 'period doubling', but clearly any tile with an odd number of vertices is excluded so all the M-type tiles will have even period. The D's may or may not have period doubling. See Appendix D.

**Lemma (Propensity of even periods):** For any generating polygon, if a tile has even period k, then all points in the tile have period k. If the tile has odd period k, the all points in the tile have period 2k except for the center of symmetry which has period k.

Examples: (a) On the left below is part of the web for N = 7, showing the period 7 orbit of the center of D in magenta. Since this is an odd period, the remaining points in D will have period 14. This orbit passes through the centers of S[1] and S[2] and their orbits are also period 7 and entrained with D's orbit, so the magenta points include all three orbits in one star polygon. This is generic for regular polygons. (b) On the right is N = 7 again showing the period 14 orbit of DS[3] (D step-3). This is an off center point but the period of the tile is even so all points have period 14. The step sequence for this orbit is {3,2}, but it visits each D in a step-3 fashion, so relative to the D's the step sequence is {3}. In this sense the seven D's act as one.

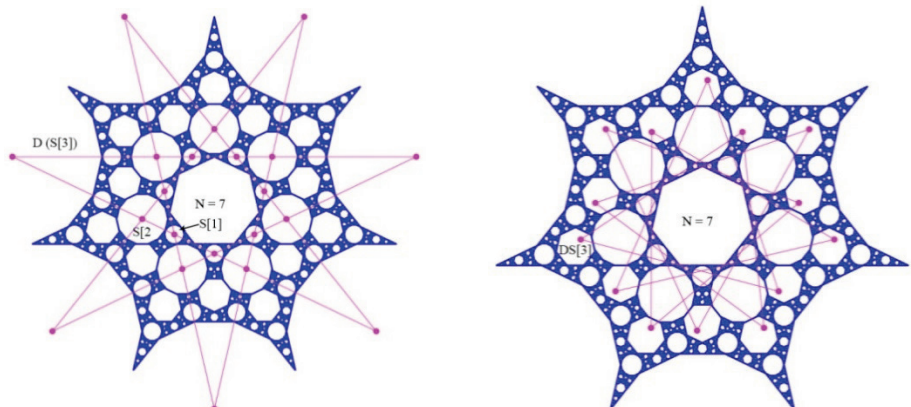

D is step-3 relative to N = 7 (M) but if D was the generating polygon, M would be step-5 relative to D so these two tiles generate each other. The web evolution is described below.

# Evolution of the web for regular polygons

The large-scale web for any regular polygon is dominated by rings of large resonances called D's. These 'necklaces' of D's guarantee that the region between rings is invariant and this guarantees that no orbits are unbounded [VS]. There is a natural conjugacy relating the dynamics of each inter-ring region. The region inside the first ring is called the central 'star' region and it serves as a template for the global dynamics. Below are typical webs for N odd and N even. In the N-even case the generating polygon is itself a D and the local geometry of all the D's is conjugate. Section 6 is devoted to the large-scale dynamics of regular polygons.

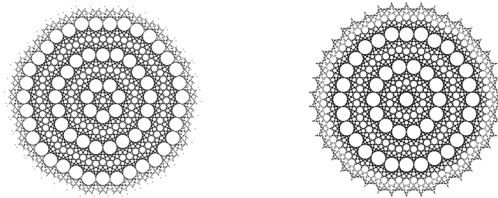

Below are the first 4 iterations of the inverse web for N = 7 in blue with the level 0 forward web shown in magenta for reference. The evolution is similar to N = 5. The central 'star' region will be formed from $L_1$ and $L_2$ and their symmetric counterparts. The $L_3$ rays will form the outer edges of the ring of seven D's and the remaining large scale web. The local evolution of the canonical D is described below.

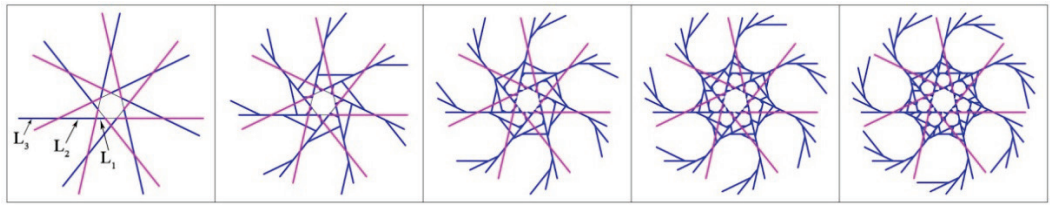

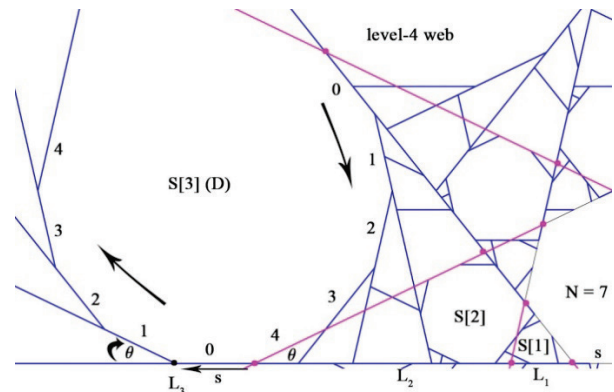

To track the web development of D's external edges, start with the end point of $L_3$ on the extended forward edge and recursively 'slide and rotate' as each iteration acts on the previous. Translations and rotation are the key elements of all piecewise isometries acting on polygons. This web algorithm is a two-dimensional version of the 'swap domain and range' cobwebbing used for functions of one variable. Each swap involves a shear in opposite directions along two edges so it is rigid rotation. The shear for τ is always the edge length of N = 7, which we call s. This guarantees that D inherits the same edge length as N = 7. The rotation angle depends on the region. The region for D is bounded by the two magenta trailing edges. The rotation angle is θ - which is half the exterior angle of N = 7, so D is a regular 14-gon. When N is even the shear is the same, but θ matches the exterior angle of N, so D is a clone of N.

The inner edges of D are generated using the symmetric point on the top edge of this region. Now the shear direction is reversed, so the $L_2$ segments generate the inner edges of D and also define the outer edges of S[2]. ($L_1$ defines S[1] and the inner edges of S[2].) The 4$^{th}$ iteration of $L_2$ is aligned with the magenta trailing edge of N = 7, and the evolution of D's inner edges will cease. (The edge numbered 4 is shown in magenta here, but it is also blue.)

This web evolution in the vicinity of D is generic for all regular polygons because the D's always have one edge which lies on a trailing edge of the generation polygon and one edge which lies on a forward edge. For N odd, it only takes N-3 iterations for the trailing edge to map to a forward edge so the inner edges of D no longer evolve. This implies that the inner star region is invariant after N -3 iterations for N odd (N/2-2 iterations for N even).

D's region spans 3 forward edges since it is step -3. This span is maximal for N = 7 so D has the largest possible measure. The development of S[1] and S[2] are similar, but complicated by the congestion of the 'inner star' region. In all cases the shear is the same, so the centers of these regions are displaced by s/2 outwards from the corresponding 'star' point -where the forward edges meet the trailing edges. These points are analogs of hyperbolic fixed points.

The rotation angles of these domain regions are the 'star' angles. When N is odd, they are of the form $\pi-k\varphi$ for consecutive integers k, where $\varphi =2\pi/N$ is the exterior angle of N. When N is even they are of the form $k\varphi$. For example with N = 7, $\theta$ is clearly $\pi-3\varphi$, and the next two are $\pi- 2\varphi$, followed by $\pi-\varphi$ which is the interior angle of N = 7. Therefore N = 7 can be generated by these algorithms. Eventually D will have a canonical tile in each of his six step regions. S[2] is the only tile which is shared by both D and N = 7. For N = 14, the shared tile would be a scaled copy of N = 7 and in this case D and N= 14 would have identical webs and conjugate dynamics.

On the left below is the first iteration of the web local to S[1] – showing the 'swap domain and range' algorithm applied to $L_1$. Iterates of $L_1$ will generate S[1] and the interior edges of S[2]. This 'inner star' region is invariant as shown on the right below. The black 'outer star' region is generated by $L_2$, with S[2] on the boundary. Since the D's are step-3 their rotation per iteration is 3/7 (on a scale of 0 to 1). This is their 'winding number'. This star region can be regarded as a 'winding number trap' because any escaping point would have to have a winding number which exceeds 3/7 and this is impossible since their orbits are linked with D's. In the same fashion, S[2] with winding number 2/7, is an upper bound for the inner star points. The general issue of invariance is not this simple, but there should be ⌊N/2⌋-1 regions for prime N- and typically an infinite number of non-trivial secondary invariances.

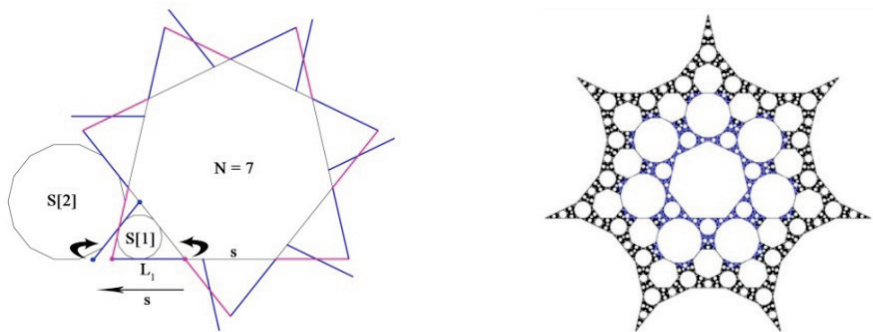

Section 2 - The geometry of the First Families

Since the First Families are the major resonances of the Tangent map, they form early in the web generation process. This means that it is possible to derive the geometry of these canonical tiles from a recursive application of the web generation algorithm as shown above. Extending this web geometry to secondary tiles is far from trivial and the standard renormalization techniques only apply in a few select cases. However the recursive nature of the web for regular polygons means there is potential for self-similarity. In many cases there appear to be canonical tiles at all scales and the First Family can provide a template for future generations of tiles.

The geometry of the First Family is a function of the 'star points' where the extended trailing edges of the N-gon meet one of the forward edges. Typically these 'hyperbolic' star points support their own local dynamics and infinite chains of families. Secondary tiles have their own star points which can play similar roles.  Section 3 will investigate the dynamics of a few primary and secondary star points for N = 7.  (Below we will assume that the generating polygon is centered at the origin.)

**Star Point Theorem**: *All regular N-gons have ⌊N/2⌋ (Floor[N/2]) 'star' points for N odd and N/2 -1 star points  for N even. These star points are intersections of extended trailing edges of the generating N-gon with a single extended forward edge- which can always be assumed to be a horizontal edge as shown below for N = 11 and N = 16.  By convention the star points are numbered starting with star[1] which is a vertex of the generating N-gon.  The last star point is also known as GenStar[N] because in some contexts there are infinite chains of generations converging to this point.*

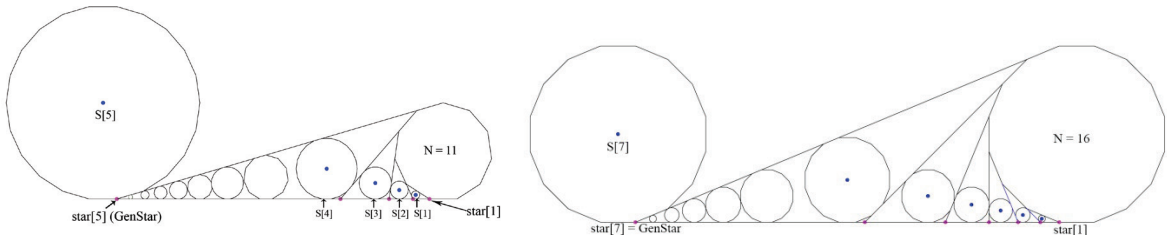

*The star points define transitions in the dynamics of the Tangent map τ . The region between star[k] and star[k+1] is called the 'k-step region' because it contains a unique orbit with constant k-step (advancing k vertices on each iteration). This orbit is always period N –but the prime period may be a divisor of N. The points with this k-step orbit represent a canonical 'resonance' of τ. They form a convex polygon (tile) called S[k]. (The last S[k] is also known as D.)  The center of region S[k] is displaced horizontally outwards by s/2 from star[k], where s is the side length of the generating N-gon. These centers lie on a line of symmetry extending from star[1] to cD (center of D) which is well-defined because the vertical coordinate of cD matches either the height or the center of the generating N-gon – depending on whether N is odd or even. Therefore the star points define the centers of the S[k].*

When N is prime the S[k] will be regular 2N-gons with prime period N. When N is composite, some of the S[k] may be non-regular with periods which are factors of N. These non-regular S[k] are always based on a regular 'template'. In all cases, the last S[k] resonance is a regular 2N-gon for N odd or N-gon for N even. This S[[N/2]] or S[N/2-1] tile is also known as D – but D-type tiles can exist at all scales. Note that the orbits of these canonical family members define regular star polygons with Schläfli symbols {N/k} for k = 1 to ⌊N/2⌋ for N odd and k = 1 to N/2-1 for N even. As indicated earlier for N = 7, the centers are themselves star points. They are the star points of the star polygon defined by the orbit of the center of D.

**Definition of the First Family**

For all regular polygons, the D tiles have their own local webs which form as part of the global webs. These local webs are identical to the webs which would form if D was the generating polygon so the D's can be used as equivalent generators for all regular polygons.

The definition of the First Families will include the S[k] defined above as well as the corresponding DS[k] of D – but the S[k] and DS[k] are closely related. The basic structure described below depends only on whether N is odd or even, but the detailed structure in the next section requires four cases: N prime, N odd and composite, N twice odd and N twice even.

When N is odd, the D's will be regular 2N-gons with the same side length as M and there will be ⌊N/2⌋ star points, so there will be ⌊N/2⌋ S[k]'s. The matching D will define another N-1 DS[k]'s. D and M share S[⌊N/2⌋-1], so the total count for the two families is ⌊N/2⌋ + N-2. We will define the First Family to be the combined S[k] and DS[k], but excluding the last DS[N-1] region for D which is an identical D on the right of M –sometimes called DRight. This means the First Family for N odd will have ⌊N/2⌋ + N-3 members.

The 13 members of the First Family for N = 11 are shown below. Note that DS[8] is the same as S[4] since D is just two copies of M. By this same reasoning, DS[10] is the same size as S[5] (D). This duplicate DRight is omitted here by convention.

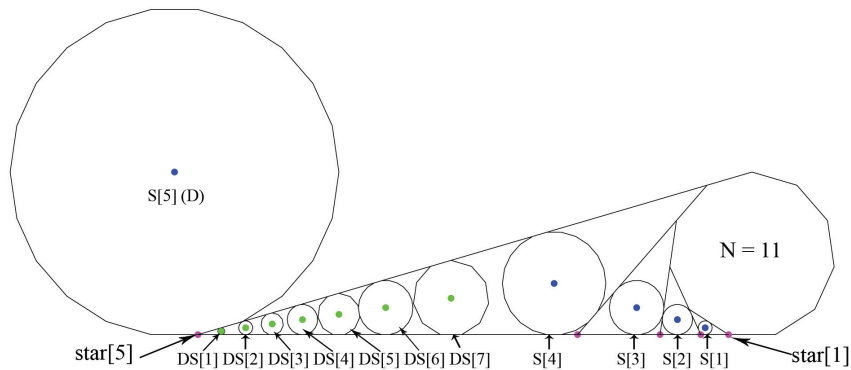

The DS[k] for N = 11 can be found from the Star Point Theorem using N = 22, but they can also be found using N = 11 because D's web geometry can be derived from M's. Each S[k] of N = 11 accounts for two DS[k]'s of N = 22 – one which is identical in size to S[k] and the other which has the same 'scale' but 11 sides instead of 22. This will be explained below when we introduce the notion of scaling. The situation for composite odds is similar, but some of the S[k] and DS[k] may no longer be regular.

When N is even, D is the same as the generating N-gon and their local webs are conjugate, so the First Family is just the union of the S[k]'s and matching DS[k]'s arranged symmetrically around S[N/2 -2]. This means there are 2·(number of starpoints) -1 members = 2(N/2-1)-1= N-3. Below are the 19 members of the First Family for N = 22 arranged symmetrically about S[9]- which is a clone of the N= 11 (without S[1],S[2] and S[3] which are secondary resonances here). The symmetric DS[K] for the N-even case are sometimes called LS[k]'s instead of DS[k]'s.

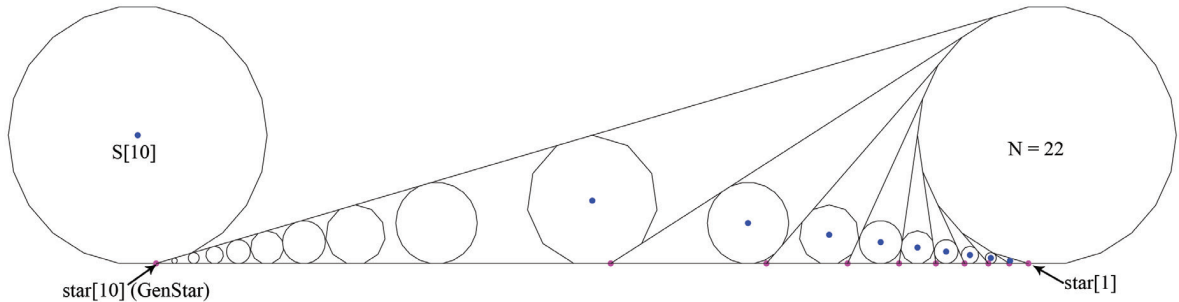

The twice-even and twice-odd cases have similar large-scale structure, but they differ in detail. In all twice-odd cases, N = 2k and N = k have congruent family structure as illustrated above for N = 22 and N = 11. When N is twice even, the First Family retains this same symmetry but there are no canonical N/2-gons so N = 4k and N = 2k may have very different dynamics. N = 16 is only distantly related to N = 8, and N = 12 bears little resemblance to N = 6 (which is congruent to N = 3 because it is twice odd).

The First Family for N = 12 is shown below. There are canonical tiles which 'see' only the embedded N = 3, 4 and 6, shown here so it is not surprising that mutated hexagons and octagons appear among the canonical tiles. The periods of S[2], S[3] and S[4] are 6,4 and 3 instead of 12. This type of orbit decomposition is generic for composite N-gons, but not all orbit decomposition yields mutated family members. Here only S[2] and S[3] are mutated. The web plot shows the self-similar nature of the dynamics. N= 5, N= 8 and N = 12 all have coordinates which can be defined quadratically and they have self- similar web structure. N = 12 has a self-similar sequence of S[1] tiles converging to the star[1] point (and the symmetric GenStar point).

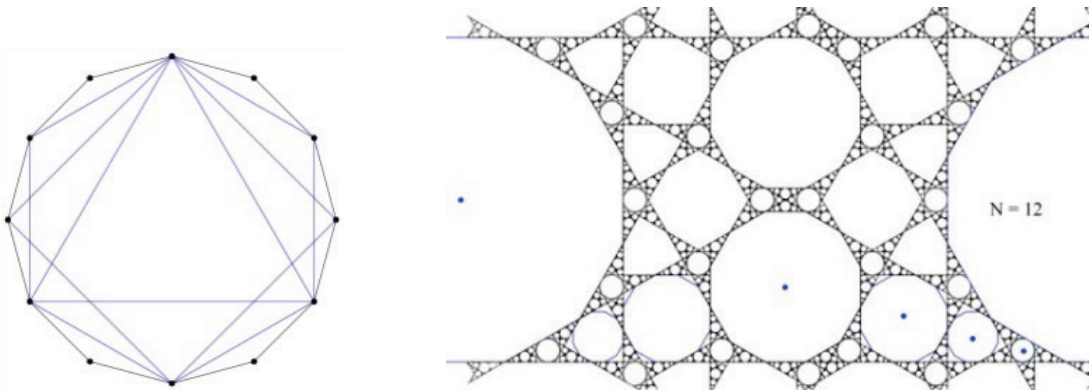

## Scaling

The one remaining issue in the definition of the First Families is the scale of the family members. Every star point defines a scale and conversely, but the scales are ratios and hence independent of size.

**Definition:** For a regular N-gon, scale[k]= s/(-2*star[k][[1]]) where star[k][[1]] is the horizontal coordinate of star[k] and s is the side length of N.

Since star[1] is a vertex of the generating polygon, every regular polygon has scale[1] = 1. The outermost star point is called GenStar[N] and the corresponding scale is called GenScale[N]. This is the scale for new 'generations' (if they exist). This will be explained in Section 3.

The one-to-one correspondence between scales and star points means that each step region has a matching scale, and the scales can be used to define the canonical S[k] tiles in each step region. Since these step regions define transitions in dynamics, the scales play the same role. For example with N = 7, S[1] and S[2] have scales which are determined by their dynamics and this is true for all subsequent generations. The critical issue is how these scales interact and this depends largely on whether they are commensurate. These scales have their own 'winding numbers' so they play the part of competing frequencies in a continuous Hamiltonian system. For regular polygons there should be $\varphi(N)/2 - 1$ incommensurate scales, so N = 7 and N = 9 are the first two regular polygons which have multiple incommensurate scales and their dynamics are a mixture of self-similarity and irregular motion.

In this section we will use the scales to define the parameters of the first family for the four possible cases: N = 11 (prime), N = 9 (composite and odd), N = 14 (twice-odd) and N = 16 (twice even).

**Case 1**: When N is prime, S[k] has radius rD·GenScale/scale[k] and DS[2k] has the same radius as S[k]. These are all 2N-gons, but the DS[k] for k odd are N-gons so they scale relative to M (and in reverse order), so rDS[2k+1] = rM·scale[[N/2]– k] = scale[[N/2]– k]. Therefore rDS[1] = GenScale and this tile is known as M[1] – the matriarch of the 'second generation'. DS[2] has radius rD·GenScale/scale[1] = rD·GenScale, so he is the matching D[1]. Each scale defines two matching 2N-gons and one N-gon, but only DS[1] and DS[2] have the M-D scale relationship, along with a local geometry which allows for the possibility of new generations.

Example: N = 11 has 5 star points so there are 5 scales. There are 14 tiles in the table (including D at S[5] and M at DS[9]). The graphic only has 13 tiles because we usually omit DS[10] – which is a symmetric D on the right.

| Scale | Tiles based on this scale |
|---|---|
| 1 | S[1] & DS[2] & DS[9] |
| 2 | S[2] & DS[4] & DS[7] |
| 3 | S[3] & DS[6] & DS[5] |
| 4 | S[4] = DS[8] & DS[3] |
| 5 (GenScale) | S[5] &DS[10] & DS[1] |

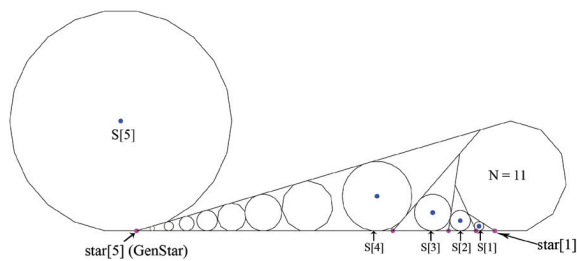

**Case 2**: When N is odd and composite, the formulas for the prime case still hold but some of the S[k] may be non-regular polygons based on the regular 'template'.

Example: The First Family for N = 9 is shown below. S[1] and S[2] have normal period 9 orbits, but the S[3] orbit only 'sees' one of the embedded N = 3 polygons, so it has period 3. This orbit decomposition yields a mutated S[3] which is composed of two nested regular hexagons of slightly different size. The matching D's for N = 3 are hexagons, so the embedded N = 3 is clearly the source of this mutation. DS[3] also feels the effect of this decomposition and becomes a non-regular hexagon composed of a regular 9-gon with extended edges.

These mutations are not false edges but rather incomplete edges caused by periodic orbits with shortened periods. The initial evolution is normal and the centers are unchanged from the regular case since they are determined by the level-0 web. Mathematica will draw perfect First Families for any regular N-gon, but for composites, a web is the best way to resolve possible mutations.

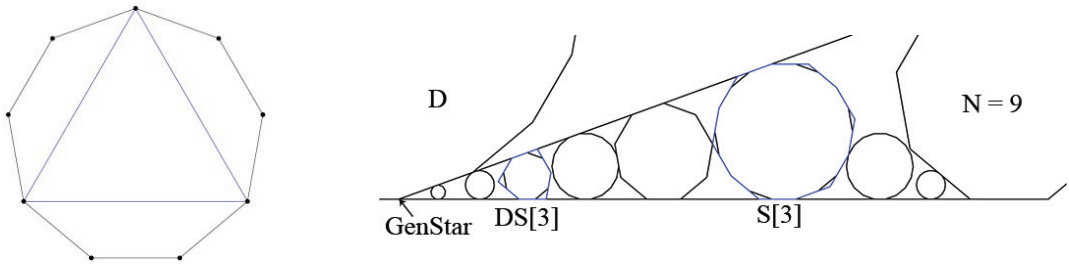

The two regular hexagons are shown in blue and magenta below. The regular S[3] template is in black. Note that it shares a vertex with S[2] on the right. This defines the radius of the magenta hexagon. The blue hexagon is defined by the vertex of DS[5] on the left. Together these two hexagons define a 'mod-3' version of the regular S[3] with every third edge extended. Polygons such as these occur in other contexts and we call them 'woven' polygons. There is a narrow window of convexity, based on the ratio of the radii. The ratio here is approximately .95418889.

As indicated above, N = 7 and N = 9 are the first two regular polygons with multiple incommensurate scales and the minimal polynomial for their vertices is cubic. For N = 9, the dynamics local to S[3] shown below are typical of the dynamics found elsewhere – namely a mixture of 'quadratic' self-similar dynamics (in red) and highly irregular dynamics (in blue) which is probably multi-fractal. N = 7 has a similar mixture.

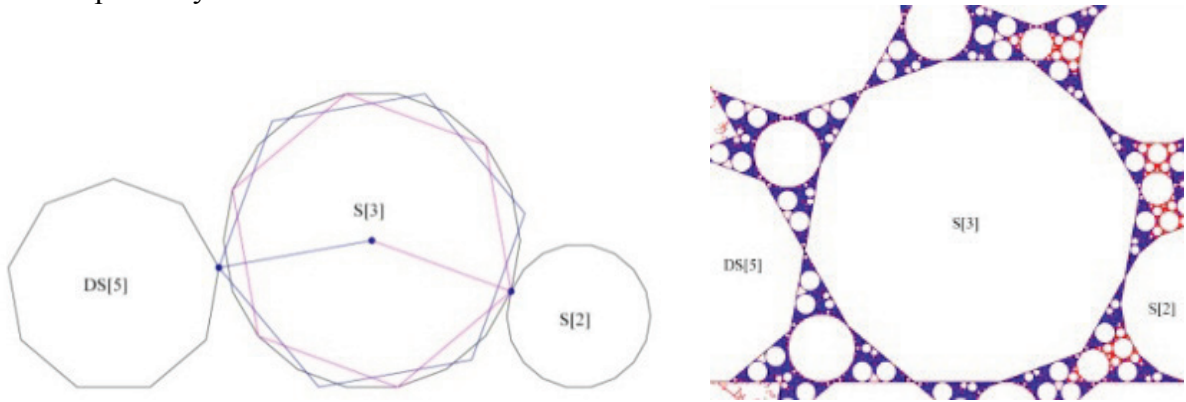

**Case 3**: N twice-odd. Whenever N = 2k, the generating N-gon is identical to the canonical D and by symmetry the First Families are symmetric about S[k-2], so the k -1 scales can be reduced to [k/2]. When k is odd, the central S[k-2] has dynamics conjugate to the k-gon 'M', so the First Family for N = 26 can be generated directly from the First Family of N = 13 (and conversely).

**Twice-odd Lemma**: (i) *For N = 2k for k odd, the number of scales is k-1 and scale[k-1] is GenScale[2k] and scale[k-2] is GenScale[k]. Therefore S[k-2] is a scaled copy of N = k and it shares the same edge length as N = 2k so it has a canonical M-D relationship with N = 2k in the sense that either tile will generate the other. (See the web evolution below)*

*(ii) As a consequence of part (i), N = 2k and N = k must have equivalent scales. Indeed the tiles of N = 2k are symmetric with respect to S[k-2], so the number of 'effective' scales is ⌊k/2⌋. To map the odd scales of N = 2k to the ⌊k/2⌋ scales. of N = k, GenScale[k]/scale[2n-1] = scale[⌊k/2⌋ –n].*

Below is the level-2 (inverse) web for N = 14 in blue together with the level-0 forward web in magenta, for reference. This evolution is similar to case of N = 7 discussed in section 1. Using the 'swap domain and range' algorithm, the ray $L_6$ generates the outer edges of D by a recursive slide and rotate. The rotation angles are the 'star' angles. For N even, these star angles are of the form $k\varphi$ where $\varphi$ is the exterior angle of N. Setting k = 1, gives the smallest angle $\theta_6 = \varphi = \pi/7$ as shown below. This will replicate N = 14. The evolution of S[5] mimics that of D, but the rotation angle is $\theta_5 = 2\varphi$ so S[5] will be a scaled copy of N = 7 with edge length the same as N = 14. This implies that S[5] and N = 14 have a canonical D-M relationship where either tile will generate the other.

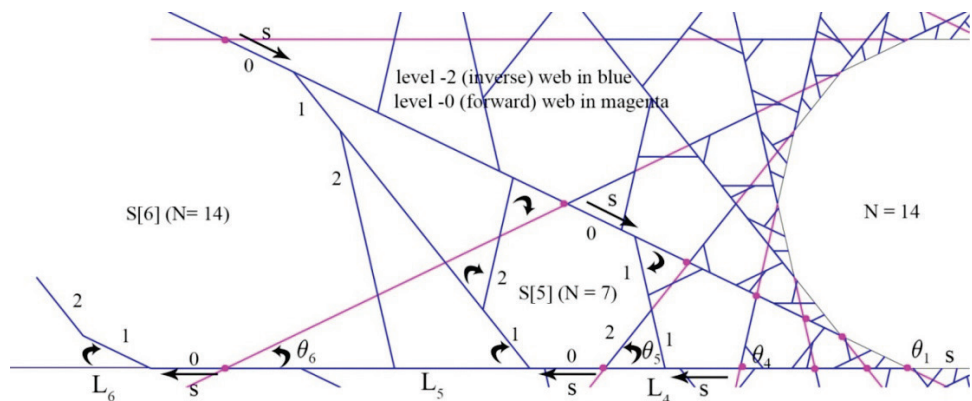

The twice-odd case to be discussed below is very similar to the twice even. In both cases the central S[k-2] tile evolves with a rotation angle of $2\varphi$ where $\varphi$ is the exterior angle of the N-gon - so there is the potential to construct an N/2-gon by a single cycle or an N-gon by two interwoven cycles. This latter case occurs iff N = 2k for k even because each of the level-0 edges (one from DRight and one from DLeft) will generate a closed step-2 cycle. Since these level-0 edges are one rotation apart, the two cycles will form the even and odd edges of S[k-2]. In this sense S[k-2] is similar to the generating N-gon- but formed in a step-2 fashion instead of a step-1 fashion – so for N = 16, the edges will be numbered 0,4,1,5,2,6,3,..This implies that S[k-2] and the generating N-gon will tend to have very different local dynamics. When N is twice-odd the two cycles collapse and N = 7 above is formed in a step-1 fashion just like N = 14 but with a different exterior angle. For more on the evolution of S[k-2], see the Df Theorem in Appendix F.

**Example**: N = 26 has 12 scales, but the 6 odd (or even) scales suffice. For example S[4] below has radius rD·GenScale[26]/scale[4] = rD·scale[9] ≈ 0.06419522834. (D's height is 1.)

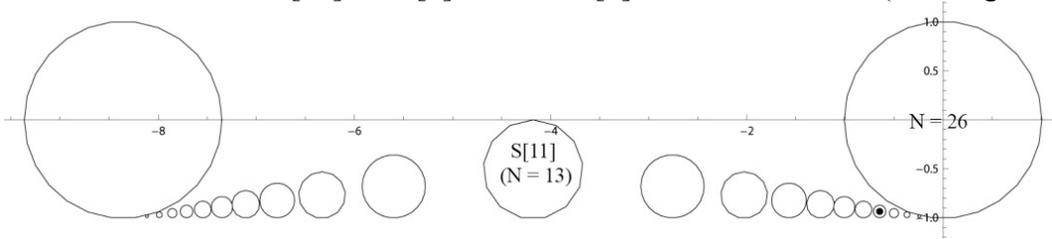

In the N = 13 world, S[4] is known as DS[4] and it has the same size as S[2] so rS[4] = rD(13)·GenScale[13]/scale[2] ≈ 0.12652506002. But GenScale[13]/scale[2] = scale[9] of N = 26, so rS[2] = rD(13)·scale[9]. This will match rS[4] above if we scale by the ratio of the two D's. This ratio is the same as the ratio of the M's which is rS[11] = RadiusFromSide[s,13] ≈ 0.5073716490048. Using this same ratio it is easy to import the rest of M's family – which are secondary family members here.

To construct the First Family for N twice-odd it is necessary to alternate N-gons and N/2-gons. The N-gons are all of the form rS[k] = rD·GenScale[26]/scale[k] as shown above for k = 4. To get the N/2-gons when k is odd, scale the corresponding tiles from N =13 using the equivalence of scales in the Twice Odd Lemma: rS[k] = rS[11]*scale[11]/scale[2k-1]

**Case 4**: N = 2k with k even. The web evolution for the twice-even case was discussed above. The symmetry of the twice-odd case is even more evident here because the canonical D-tile is now identical to N. This implies that the central S[k-2] is also an N-gon and in fact there are no canonical N/2 tiles – so DS[k] is identical to S[k] and sometimes they are called LS[k].

In terms of scaling, the basic precepts of the Twice-odd Lemma still apply. Once again the scales share the symmetry of the web, but now there are an odd number of scales. The scales map to each other by: GenScale[2k]/scale[n] -> scale[k-1-n], so odds map to odds and evens to evens. The central scale[k/2] maps to itself, which implies that GenScale[N] = scale[N/4]$^2$. (Therefore for N = 8, GenScale[8] is the only nontrivial scale.)

**Example**: N = 16 has 7 scales but by symmetry there are only 4 'effective' scales, which must be a mixture of even and odd scales. Scales 1,2,3 and 4 would suffice. For example rS[6] = rD·GenSscale/scale[6] = rD*scale[2]. In general rS[k] = rD·GenScale/scale[k].

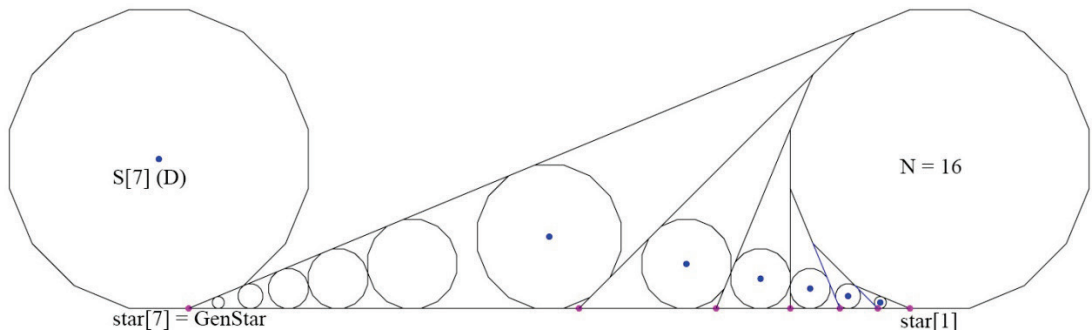

N = 16 is the first non-trivial member of the N = $2^k$ family, There appears to be little or no connection between the members of this family even though geometrically each one is a 'bisection' of the previous. N = 16 is a 'quartic' polynomial since the minimal polynomial for the vertices of any regular N-gon is degree φ(N)/2 where φ is the Euler totient function. Very little is known about the quartic case – which includes N = 15 and N = 20. In [LKV] the authors note that algebraic analysis in the quartic case appears to involve "great computational difficulties".

In twice-even cases. the canonical tiles are all N-gons – so it appears that the M-D distinction does not exist. However the notion of 'generations', which will be described in the next section for N-odd or twice-odd, still exists. For N odd, each generation (if it exists) is presided over by an M[k],D[k] pair playing the roles of matriarch and patriarch. These canonical M[k] and D[k] tiles form on the edges and vertices of the previous D[k-1], so they are always step-1 and step-2 respectively. N = 16 preserves this step-2 vs. step-1 dichotomy, so it would not be 'politicially incorrect' to associate S[1] with M[1] and S[2] with D[1]. These two tiles are shown on the left along with a magenta virtual First Family for S2[1].

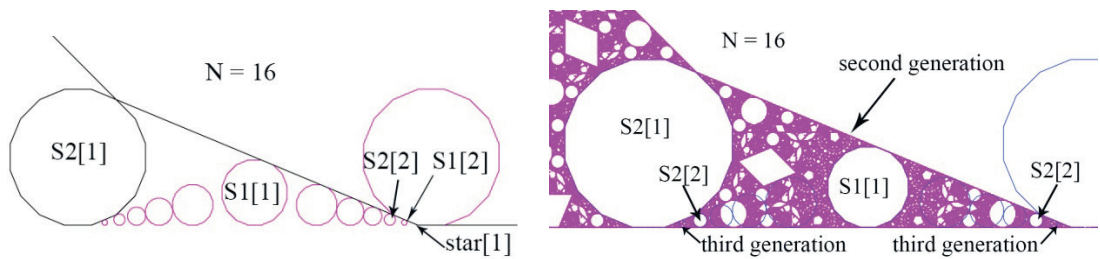

Note that S1[1] is a step-6 of S2[1] so it is natural to associate this tile with M[1]. However the web plot on the right makes it clear that S1[2] does not exist, so there is no 'M[2]]'. Since S2[2] does exist at star[1] and at the foot of S2[1], either tile can play the role of D[2]. We choose to study the foot of S2[1] because it provides valuable information about the local dynamics – which is often 'hidden' at star[1] and at GenStar. This region is enlarged below – note that the symmetry now is with respect to S2[2]. We have reproduced the virtual First Families of S2[2] to show the perfect match with four members of the third generation- including S1[3] which is the surrogate M[3]. The first 10 S2[k]'s in this sequence have periods 8, 32, 456, 2464, 20872, 110368, 974664, 5165216, 45423368 and 240668192 which gives ratios of about 4.66 and 8.8. As with N = 7, the ratios alternate high-low within the even and odd sequences. N=20 may also support generations – but they are complicated by the fact that S2[1] is a mutated decagon.

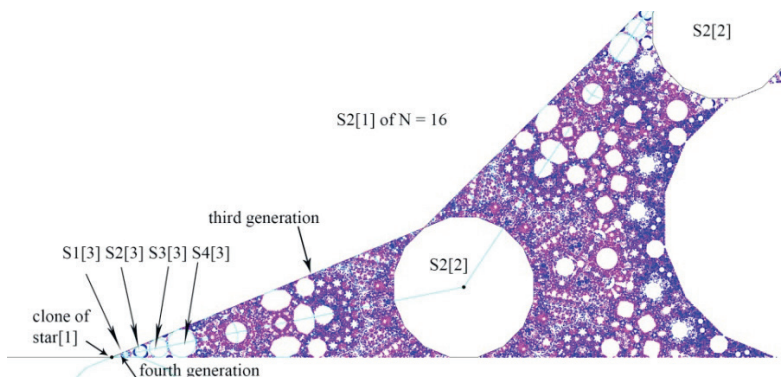

**Scaling and Virtual Tiles**

For N odd, it is easy to see that the only First Family tiles which have the canonical M-D relationship are M[0]- D[0] and M[1]-D[1] and for some polygons this M-D relationship continues at GenStar and appears to generate infinite families of M[k]-D[k] tiles. However all regular tiles have the potential to generate sequences of 'families' and it may be that all scaling parameters for regular polygons are generated in such a fashion.

Below we use N = 13 as an example of how the First Family members of any regular polygon can form their own M-D pairings as long as virtual tiles are allowed. These pairings must be nested inside the original M and D and share the same star[1] and GenStar vertices as shown below. The center lines will be normal and there is a natural extension to infinite series.

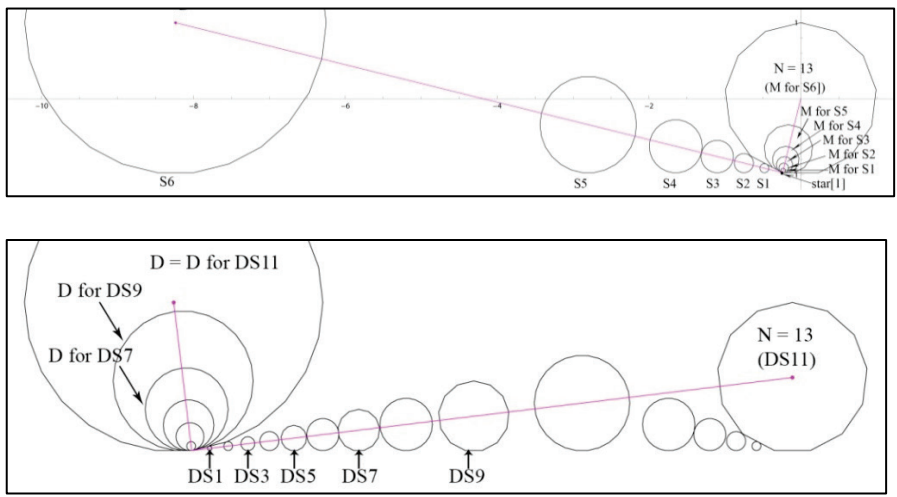

These M-D relationships carry with them the potential for complete First Families, and any alignment with actual tiles will have important dynamical consequences. The most significant alignment occurs with DS1 which is also known as M[1] because it has a perfect M-D relationship with DS2 – which is also called D[1]. As indicated earlier, this is the only case where the virtual D is a reflection of a real D - and it is an invitation for recursion. The case of DS3 is also very important because DS3 plays a major part in determining the dynamics of the $2^{nd}$ generation – which has M[1] & D[1] as matriarch and patriarch. Here DS3 has radius scale[5] so the step-5 tile will have radius rD*scale[5]*GenScale/scale[5] = rD*GenScale = rD[1]. Clearly all odd-steps of D have similar alignments. By definition, the step-1 of the virtual D's will yield the next generation, which is DS3[2] here. This will be a step-3 of the virtual D[1] but this does not guarantee that it will exist at GenStar. For N = 13 it does exist, so it also survives at D[1]- along with DS5 and DS7. But for N = 7 shown later, DS3[2] does not survive at GenStar.

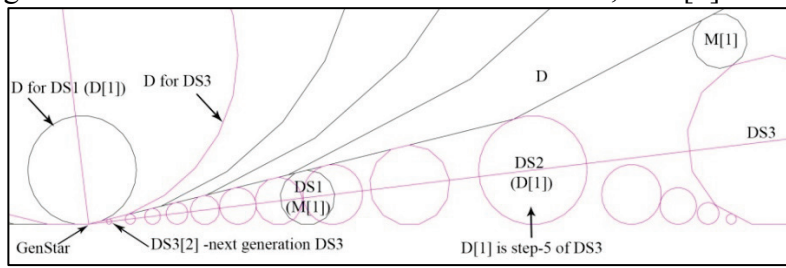

# Section 3 - Dynamics at Star Points

## Dynamics at Gen Star – the 4k+1 Conjecture

The 4k+ 1 conjecture given below claims that for 4k+ 1 prime polygons, there will always be 'chains' of generations converging to the GenStar point. Here we examine the local geometry of the GenStar point.

All regular polygons have a family structure, or extended family structure, that matches the geometry shown below, but we will assume for now that N is odd or twice-odd so the nucleus of the First Family is M[0] and D[0]. It does not matter which one is at the origin, but in the plot below D[0] on the right is assumed to be at the origin, so N is even and the local dynamics at GenStar are conjugate to the dynamics at star[1] of D[0]. But at GenStar, the candidates for D[1] and M[1] are shared by the two adjacent D[0]'s, so to recursively continue this chain of generations toward GenStar, it will be necessary to assume the existence of a virtual D[1] as shown here in magenta. These virtual D's were discussed above in the context of N = 13. They are geometric expressions of dynamical symmetry that may or may not exist. Here M[1] and D[1] are scaled by GenScale[N/2] relative to M[0] and D[0] because they are members of the M[0] First Family.

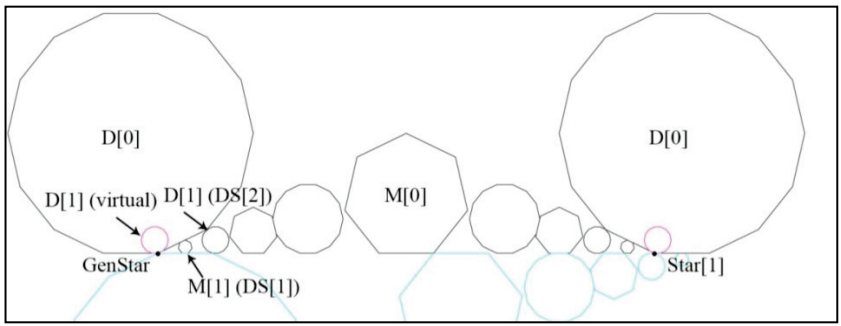

If the virtual D[1] exists at GenStar (that is, if there is symmetry between D[0] and D[1]) then there is hope for a continuation of this symmetry as shown in the enlargement below – where we have omitted M[1]. This D[2] will exist iff there is a matching D[2] in the star[1] region of D[1], so D[1] is assuming the role of D[0] and fostering a next-generation D[2]. However the self-similarity between D[k] and D[k+1] always alternates between right and left side, so the dynamics of even generations tend to differ from odd generations. There is even a slight 'bias' for N = 5 where the even and odd generations appear to be exactly self-similar.

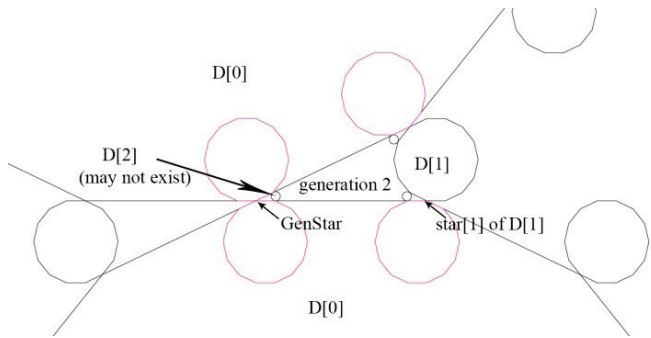

The issue is whether this generation process can be continued, so that D's and M's exist at all scales – even if the resulting generations are not self-similar. When there are converging sequences at GenStar it guarantees that there is symmetry, so it does not matter whether this convergence is studied at GenStar or star[1] of D[1]. Of course the convergence at star[1] of D[1] also depends on virtual D[k]'s as shown below, but the next-generation D[3] can be chosen on the right –side of D[2] – and once again the dynamics are conjugate to star[1] – so the limit point shown here will have dynamics conjugate to GenStar and star[1] of D[1].

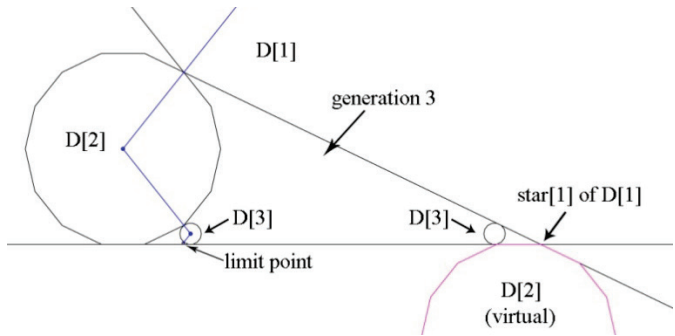

In every case studied, when N is odd or twice odd and there is a D[2] in these positions, there is evidence for an unbroken chain of generations. The 4k+1 Conjecture below states that 4k+1 prime polygons always have such chains. None have been observed for 4k+3 prime polygons – except for N = 7. However the 4k+3 prime polygons such as N = 11 and N = 19 do have M[2]'s at the correct positions, which imply that the virtual D's may have an effect even in these cases.

For N odd, the extended First Family of a regular N-gon can be used as a template for generations as shown below. To match the GenStar geometry, each generation needs to run from D to the adjacent D. For N even, both Dleft and DRight are already part of the First Family so the generation template is just the First Family itself.

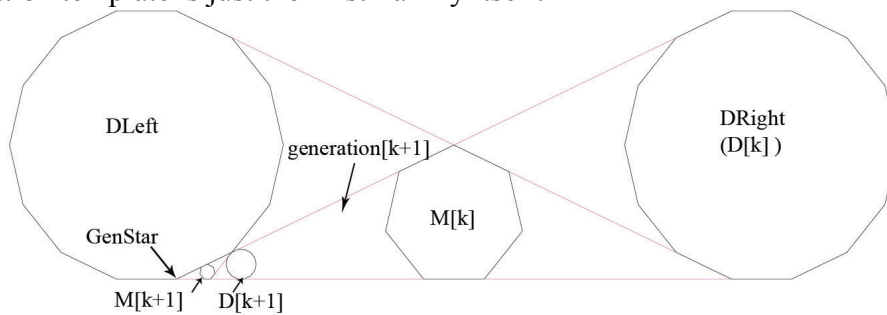

For N odd, the generating polygon is M[0] so the total horizontal span from DLeft to DRight is twice the distance from the origin to GenStar[N]. The scale that would put this extended First Family on an edge of D is GenScale[N] = s/(-2GenStar[[1]]) as defined in Section 2. For N twice odd, the scales are equivalent to the underlying odd polygon, so GenScale[N/2] can be used.

The N-even family has an extra level of symmetry which makes them valuable tools for analysis. For these families the local dynamics at the GenStar point is always congruent to the dynamics at star[1] of N.

# The 4k+1 Conjecture

The conjecture below applies only to the case where N is odd, but these results extend naturally to the twice-odd N-gons. The chains begin with the generator M[0] and the matching right-side D[0] which define the bounds of the 1st generation, but to continue this chain toward GenStar[N], it will be necessary to assume the existence of both left-side and right-side D's – some of which will be virtual.

The formulas for GenStar[N] and GenScale[N] are given below in the 4k+ 1 Conjecture. This conjecture says that these chains do in fact exist when N is prime of the form 4k+1 and moreover the ratio of the periods of consecutive D's or M's approaches N + 1.

**4k+1 Conjecture**: *Suppose M (M[0]) is a regular N-gon with N odd centered at the origin with a vertex at {0,1}. Define:*

(i) *GenScale[N] = (1-Cos[Pi/N])/Cos[Pi/N]*   (this is how generations scale under $\tau$ )
(ii) *GenStar[N] ={-Cot[Pi/N]·(1+Cos[Pi/N]), -Cos[Pi/N]}*   (the point of convergence )

*Suppose N is prime of the form N = 4k+1 for k a positive integer. Then there will be infinite sequences of regular N-gons M[j] and regular 2N-gons D[j] converging to GenStar. M[j] will have radius r[M[j]] = GenScale$^j$ and center at (1 − r[M[j]])·GenStar. The D[j]s will have radius rD[0]·GenScale$^j$ and center (1 − GenStar$^j$ ·(2+GenScale))·GenStar. The periods of these centers have ratios which approach N+1.*

**Examples**: In Section 4, we will derive the following difference equation for decagon periods for N = 5: $d_n = 5d_{n-1} + 6d_{n-2}$. With initial conditions 5 and 35, this gives decagon (center) periods of 5, 35, 205, 1235... with limiting ratio 6. This difference equation has a built-in high-low bias which occurs because the self-similar 'generation' process demands a left- right alternation as shown above. N= 8 has similar chains with a ratio of 9 – but it is unlikely that other twice-even N-gons obey the N + 1 rule. N= 7 is a 4k+3 prime but the even generations of M[k]'s seem to have ratio of 8.

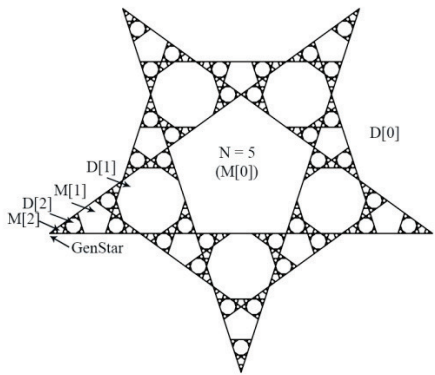

| N = 7 | Period | Ratios |
|---|---|---|
| M[1] | 28 | |
| M[2] | 98 | 3.5 |
| M[3] | 2212 | 22.57 |
| M[4] | 17486 | 7.905 |
| M[5] | 433468 | 24.789 |
| M[6] | 3482794 | 8.0347 |
| M[7] | 86639924 | 24.876 |
| M[8] | 696527902 | 8.0393 |

For N = 7 the combined ratio appears to be 200 and this would give a local Hausdorff dimension of -Log[200]/Log[GenScale[7]$^2$] ≈ 1.19978, but as shown below, other star points for N = 7 yield different growth rates. The full spectrum of dimensions for N = 7 has not been determined. See Appendices D and F.

**N = 13**: For all 4k + 1 prime N-gons, the distinction between even and odd generations creates a high-low alternation in the ratios as shown below for N = 13. Since the M's for each generation form on the edges of the D's from the previous generation, their alternation is reversed relative to the D's. The tables for N = 5 and N = 17 are very similar.

| Six generations of M's and D's for N = 13 | | | | | |
| --- | --- | --- | --- | --- | --- |
| Generation | Period | Ratio | Generation | Period | Ratio |
| M[1] | 10*13 | | D[1] | 9*13 | |
| M[2] | 182*13 | 18.20000 | D[2] | 119*13 | 13.22222 |
| M[3] | 2506*13 | 13.76923 | D[3] | 1673*13 | 14.05882 |
| M[4] | 35126*13 | 14.01676 | D[4] | 23415*13 | 13.99582 |
| M[5] | 491722*13 | 13.99880 | D[5] | 327817*13 | 14.00029 |
| M[6] | 6884150*13 | 14.00008 | D[6] | 4589431*13 | 13.99998 |

On the left below is the invariant outer star region for N = 13. In the last section we noted that there are step-relationships between D[1] and the odd-step tiles. This may explain why these odd-steps have a higher 'survival rate' in future generations than the even-steps. For N = 13, D[1] is a step-5 tile of DS3 but this is the only step-tile of DS3 that exists. However we noted that the virtual D for DS3 has a step-1 which survives as a DS3[2] and this true for DS5[2] and DS7[2]. These tiles can be seen below in the 2$^{nd}$ generation at D[1] and at GenStar. In subsequent generations DS3 seems to survive. This is also true for the 'twin' 4k+1 prime at N = 17. For N = 7, DS3[2] does not exist at GenStar but it returns on odd generations. For N = 11 and N = 19 there is no D[2] or DS3[2] and no sign of future generations.

For N = 13, it is not clear what are the effects of the step-5 relationship between DS3 and D[1] but for all odd N-gons, DS3 and M[1] are closely related because they share a vertex. In the web plot below it should be clear that M[1] and DS3 have similar dynamics. In fact all the tiles in this region are closely related. (The magenta-blue distinction occurs from two web scans centered on M[1] and D[1].)

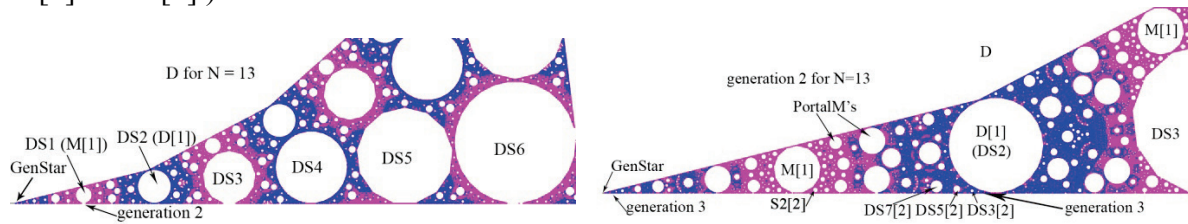

The second generation on the right continues this M-D dichotomy which extends all the way to the GenStar region. This is not surprising since there is a natural self-similarity between D[1] and D as described earlier. All of these statements should remain true for any odd N-gon so there is always the potential for further generations to exist at GenStar.

The second generation shown here is quite distict from the first and this is due in part to the interaction between M[1] and DS3. These are both M-type tiles but rDS3 = scale[5] and rM[1] = scale[6] (GenScale[13]). Since these scales are incommensurate, there is no common ground for sharing tiles, so the result is an unpredictable mixture of non-canonical tiles which we call PortalM's (although sometimes they are D's).

In the second generation shown above, M[1] has surviving step-1 and step-2 tiles which are both shared by DS3 - but they are the wrong size and position to be canonical tiles of DS3. As indicated above, DS3 has no surviving canonical tiles except for D[1] (which is also shared with M[1] since D[1] is step-6 of M[1]). The resulting PortalM's are named after similar tiles which occur with N = 7.

From the 4k+1 Lemma and the table above, there is reason to expect that future generations will settle down into predictable patterns - but, unlike N = 7, the DS3's for N = 13 exert influence on every generation. This means there is potential for an endless variety of small scale dynamics within the M-D family template. Since the new generations form at the foot of D[1], an important issue is the dynamical connection between DS3 and D[1]. As N increases, the gap between DS3 and D[1] increases, so the dynamical influence between these two should diminish.

For N = 13, it is clear that the second generation has altered the local dynamics around D[1] to produce the third generation shown on the left below. As expected, the resulting fourth generation shown on the right resembles the second generation more than the third – but it is still quite distinct. Within the 'right-left' dichotomy, every generation for N = 13 may be unique. N = 17 seems to have similar dynamics – with the DS3's surviving and small-scale dynamics distinct for each generation. This may be generic for 4k+1 primes.

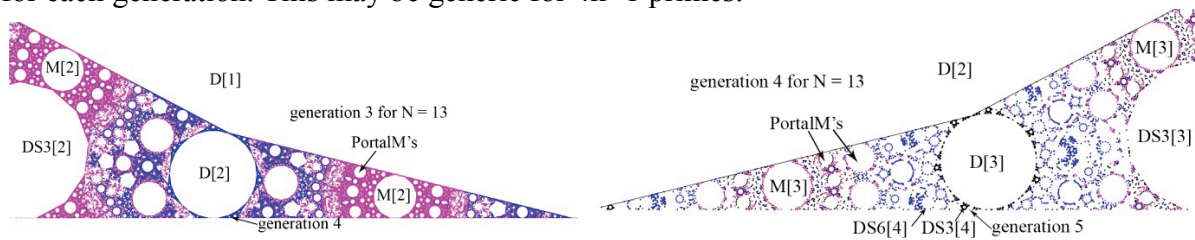

## Generation evolution for N = 7

It appears that for N odd, the interaction between M[k], D[k] and DS3[k] is an important driving force behind the dynamics of generations converging to GenStar. As indicated above, the case of N = 5 is fairly trivial because DS3 is M, so there is no compatibility issue between DS3 and M[1]. For the remaining odd N, DS3 and M[1] will have consecutive scales so they will have no canonical tiles in common. N = 7 is unique in that D[1] is a step-2 of DS3 - so the D[2] tiles generated by D[1] and M[1] are compatible with DS3 as shown below. However the step-2 tiles of M[1] are not compatible with DS3, so the edges of DS3 are dominated by non-canonical PM[2]'s which are shared by M[1] in place of the usual S2[2]'s. Therefore the next-generation DS3[2]'s only exist at a few places, and there are none at GenStar or at the foot of D[1].

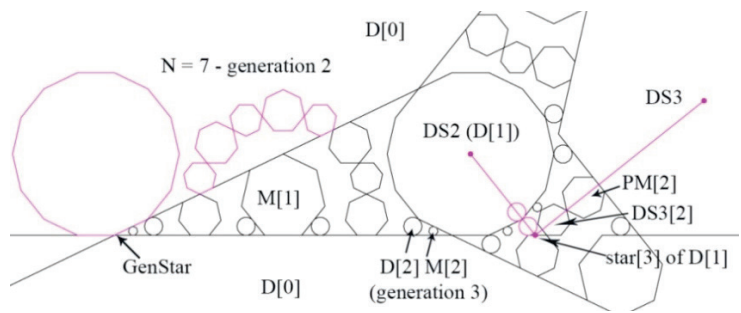

These missing DS3[2]'s actually play a positive role in the 3nd generation – because there is no longer a compatibility issue with M[2] – and the 3rd generation can evolve 'normally' to be self-similar to the 1st generation as shown below. Therefore the evolution process will be repeated with D[2] as the new D[0]. It would be reasonable to assert that the 'hybrid' nature of N = 7 is due to the natural 'cubic' symmetry combined with the step-2 relationship of DS3 and D[1], but N = 9 has both of these traits, and the resulting PM's do not prevent the formation of DS3's at GenStar. Therefore the 3rd generation is identical to the 2nd – and hence all subsequent generations are self-similar – with ratio of periods satisfying the N+1 rule. The larger separation between DS3 and D[1] may be a factor here and in N = 13.

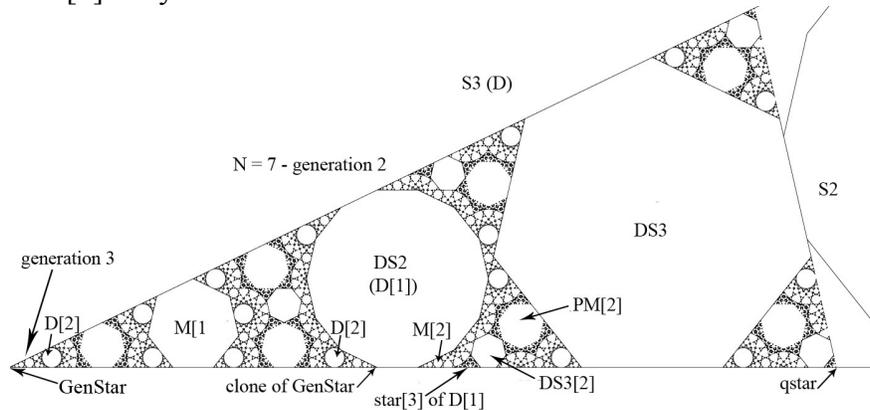

### Evolution at Star[3] of D[1] for N = 7

Since the star points of all regular N-gons are transverse intersections of forward and trailing edges, they are candidates for convergent sequence of tiles. It is possible that all invariant measures arise from sequences of tiles converging to primary and secondary star points. The canonical convergence at GenStar is just one of an infinite number of possible sequences. The vector plot above shows a sequence of D[k]'s converging to the star[3] point of D[1] – while from the other direction, a sequence of PM's and DS3's converge to the same point.

The D[k] sequence alternates real and virtual, so D[2] is virtual. It is shown in magenta in the vector plot above – along with its reflection. Both of these virtual tiles play a part in the dynamics. Below is an enlargement of this region. It is normal for rings of M's to form around D tiles and on the right-side of D[1] the rings are centered on a reflection of the virtual D[2]. The blue M[2]'s in this web will be real but truncated to match the colored image shown here. (The extended edges are not unusual for webs.) The result is a string of non-regular pentagons which are formed in a unique fashion. These strings can be seen in the web plots above and below.

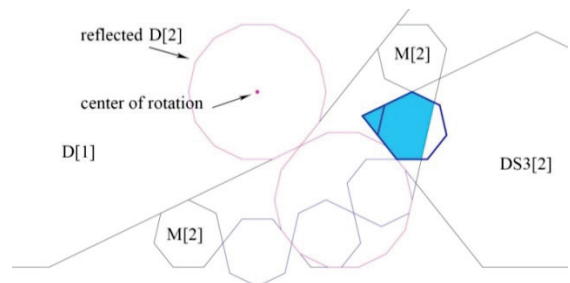

On the right side of star[3] (which is also star[2] of DS3), the DS3 chain has physical gaps – which are filled on even generations by PortalM's. The growth rate of periods is 113 and this matches the growth rate at qstar on the right side of DS3 and also the growth rate at star[2] of N = 7 – which we will investigate in Section 5. In all cases this growth rate skips one generation so the local Hausdorff dimension is $Ln[113]/Ln[1/GenScale[7]^2] \approx 1.071$ compared with 1.19978 at GenStar.

Note: There is no doubt that DS3 plays an important role in the dynamics of this region, but typically DS3 is not the 'dominant' tile in the outer star region. All regular polygons have invariant outer star regions with tiles that share similar dynamics. For N odd, the largest tile in this region is DS[[N/2]]. This means that 4k+1 primes will have a 2N-gon in this position and this may help to explain the difference in dynamics between these two classes. Having a 'D-type' tile in this position could influence the formation of next-generation M-D pairs at GenStar.

The 'towers' on the edges of the PM's below are bounded by the extended edges of the PM's, but the interior dynamics are determined by pairs of virtual D[2]'s such as those shown on the right. This geometry also exists in the first generation at GenStar. It is called a 'short family' because the central 'M' is really a DS3 and the rest of the family are DS1 and DS2- which are M[3] and D[3] tiles here. So this is not a traditional M-D relationship and the surrounding magenta webs are different from the traditional webs shown in blue. These magenta webs form inside the virtual D[2]'s. They are 4$^{th}$ generation webs – with embedded M[3]'s surrounded by rings of D[3]'s. This is the first place that such webs appear. This whole region highlights the conflict between D[1] and DS3 – which are step-2 and step-3 of D, so this is the 'outer-star' version of the S1-S2 conflict at star[2] of N = 7, which will be discussed in Section 5. (Note that the star[3] convergence shown below takes place on the edge of a tower sitting on top of a PM[2] which is off the screen at the bottom.)

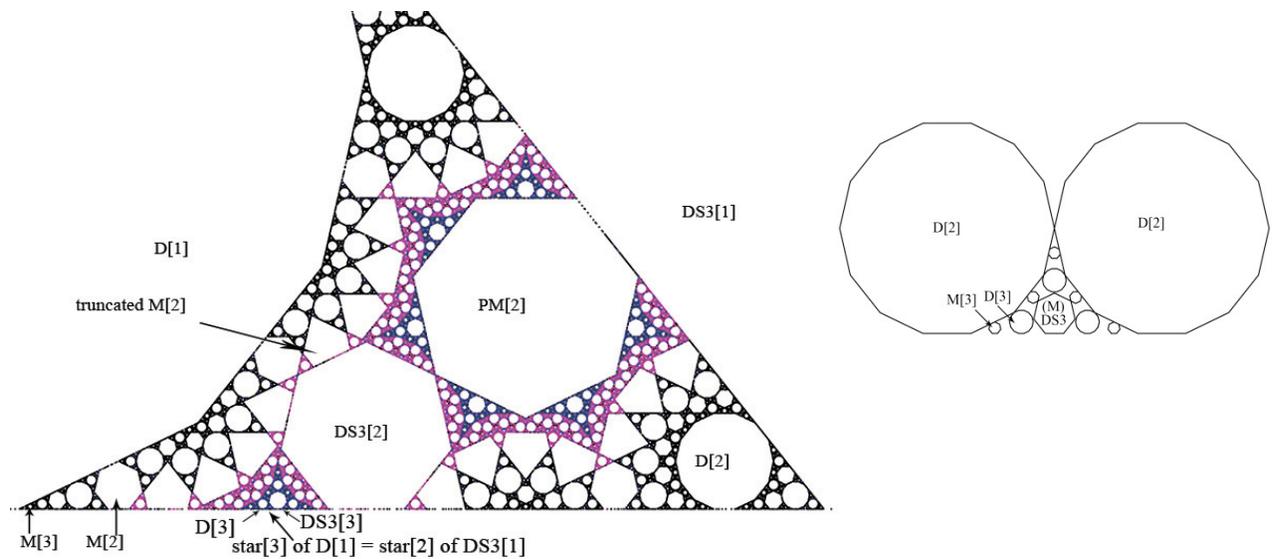

**Section 4 – The algebraic complexity of regular polygons**

The rank of the linear space generated by the vertices of a regular polygon is φ(N), but by symmetry, the minimum polynomial needed to specify the vertices is degree φ(N)/2 and we will define this to be the 'algebraic complexity' of a regular polygon. See Appendix B. The table below lists the degree of the first few regular polygons:

| Degree of minimal polynomial - φ(N)/2 | Polygons |
|---|---|
| 1 'linear' | N = 3,4,6 |
| 2 'quadratic' | N = 5,8,10,12 |
| 3 'cubic' | N = 7,14,9,18 |
| 4 'quartic' | N = 15,16,20 |
| 5 'quintic' | N = 11, 22 |

There appears to be a correlation between algebraic complexity and dynamical complexity- with members of the same algebraic 'class' sharing similar dynamics. The three linear cases N = 3, 4 and 6 are affinely equivalent to rational polygons – also known as lattice polygons. These are the only affinely regular lattice polygons and their webs have no accumulation points, so all orbits are periodic.

The quadratic polygons N = 5, 8 and 10 and 12 all have a fractal web structure with a single geometric scale factor and a corresponding 'temporal' scale which can be determined by renormalization methods. The two cubic polygons N = 7 and N = 9 have two and three competing scales, and a mixture of self-similar and irregular dynamics. The quartic cases of N = 15,16 and 20 retain some degree of self-similarity despite the 'mutations' that occur in critical tiles. The quintic case of N = 11 will be discussed below. It does not show the type of generation similarity found in the 4k+1 primes such as N = 13 – which is the lone sextic case.

N = 5 has a self-similar web so it is easy to find a set of difference equations which describe the evolution. These same 'renormalization' methods can be applied to any self-similar structure. In [T](1995) S.Tabachnikov used a 'template' like the one shown below to obtain a renormalization scheme for N = 5. This light-blue region 'tiles' the inner star and is invariant under $\tau^{10}$ so it can be used as a 'surrogate' for the entire web. The table describes how this region could be tiled with decagons and pentagons as they decrease in size by GenScale. The count is aided by fact that the dark blue region is also scaled by GenScale and it is self-similar to whole template.

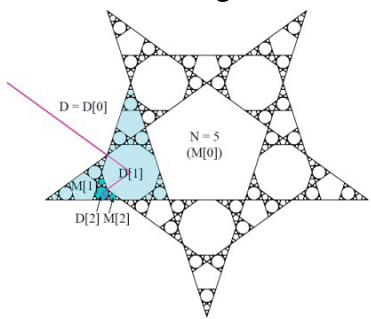

| Generation | decagons - $d_n$ | pentagons - $p_n$ |
|---|---|---|
| 1 | 1  (D[1]) | 2  (M[1]'s) |
| 2 | 7 = $3d_1 + 2p_1$ | 10 = $6d_1 + 2p_1$ |
| 3 | 41 = $3d_2 + 2p_2$ | 62 = $6d_2 + 2p_2$ |
| n | $d_n = 3d_{n-1} + 2p_{n-1}$ | $p_n = 6d_{n-1} + 2p_{n-1}$ |

It should be clear that the error in the tiling decreases at each new generation – and it is also true that the D's suffice. These difference equations can be used to find a coding sequence for the D[k] and M[k]. This will be done in Appendix C. Here we will find the periods - as discussed in the last section. These two equations can be combined together to give a second-order equation: $d_n = 5d_{n-1} + 6d_{n-2}$. On the full 'star' region these same relationships hold but the initial conditions are $d_1 = 5$ and $p_1 = 10$, so the solution is: $d_n = \frac{5}{7}[8 \cdot 6^{n-1} + (-1)^n]$. This gives decagon center periods of 5, 35, 205, 1235,.. and pentagon periods of 10,50,310,…(Our convention is to measure the periods of decagon tiles by the period of their centers, so the number of decagons matches the period – and this is also true of the pentagons.) These difference equations show that the ratio of the periods for the D's (and M's) approach 6 as in the 4k+1 conjecture. These D's can be used to 'cover' the star region and the error in this cover decreases with each generation, so the Hausdorff-Besicovitch fractal dimension is Ln[6]/Ln[1/GenScale[5]] ≈ 1.24114.

The magenta sequence of D[k] centers shown above defines a limit point s* with non-periodic orbit. The 'address' of s* is {2, 5, 2, 5,…} where the 10 'buds' of each D[k] are number using the convention that bud-1 is at the 3:00 position. This limit point is clearly on the extended forward edge of N = 5. An equivalent sequence is q* ={2,3,9,3,9…} which has co-ordinates {M[[5]][[1]], M[[4]][[2]]} – where {M[[1]][[1]], M[[1]][[2]]} = {0,1} by convention.

N = 8 has a self-similar web which can be analyzed by these same methods. Even though this is not a prime polygon, the ratio of the D periods approaches 9 so the fractal dimension is Ln[9]/Ln[1/GenScale[8]] ≈ 1.24648.

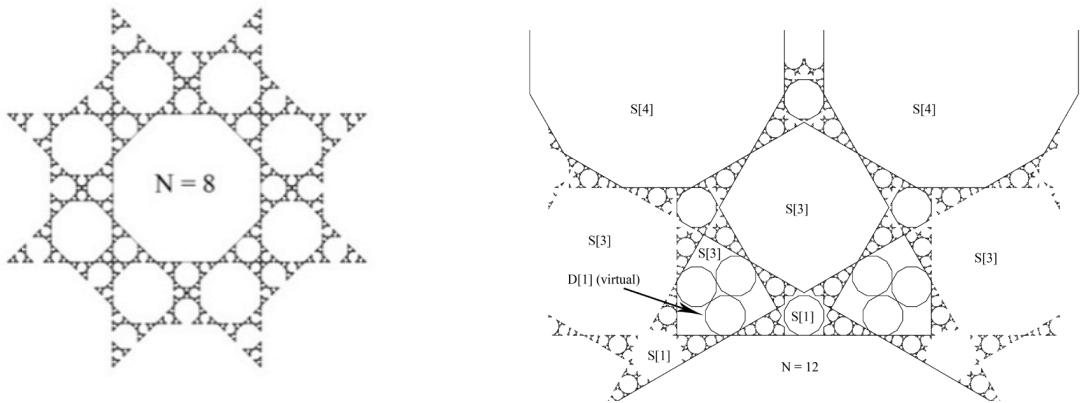

The case for N = 12 is not fully resolved but it clearly supports a form of global self-similarity. At GenStar and star[1] the dynamics are identical to the top edge of N = 12 shown above. Even though there are mutated tiles, the self-similarity is evident with S[1] playing the part of the new D[1] – scaled by GenScale[12]. The temporal scaling throughout is 27 so the fractal dimension should be Ln[27]/Ln[1/GenScale[12]] ≈ 1.2513.

It is probably no coincidence that these three fractal dimensions are increasing - since the fractal dimension is one measure of the complexity of the dynamics. For N = 7 and beyond, the dynamics are typically multi-fractal, but the maximal Hausdorff dimension can still be used as a measure of the complexity. In Appendix F this issue is discussed further and we conjecture that the maximal Hausdorff dimension will approach 2 as N→∞.

**N = 11**

N = 11 is the second 4k+ 3 prime polygon, but it is more 'typical' of that family than N = 7, which is clearly a 'hybrid', since it shares some characteristics of the 4k + 1 family. N = 11 (and the matching N = 22) are the only 'quintic' regular polygons. Their algebraic structure and matching Galois group are discussed in Appendix B.

The 4k+ 3 prime polygons beyond N = 7, all appear to have very complex dynamics with few signs of self-similarity. The generation structure breaks down very quickly – but all regular polygons have canonical First Families, which imply that there is a D[1] and M[1]. There is no evidence of a D[2] for N = 11 or N = 19, but both of these polygons have M[2]'s.

Below is a web plot of the second generation for N = 11 showing the dynamical connection between D[1] and GenStar that is apparently generic for N odd. This means there is still potential for future generations at GenStar. The symmetry relative to M[1] is preserved, so the geometry at GenStar is conjugate to star[1] of D[1].

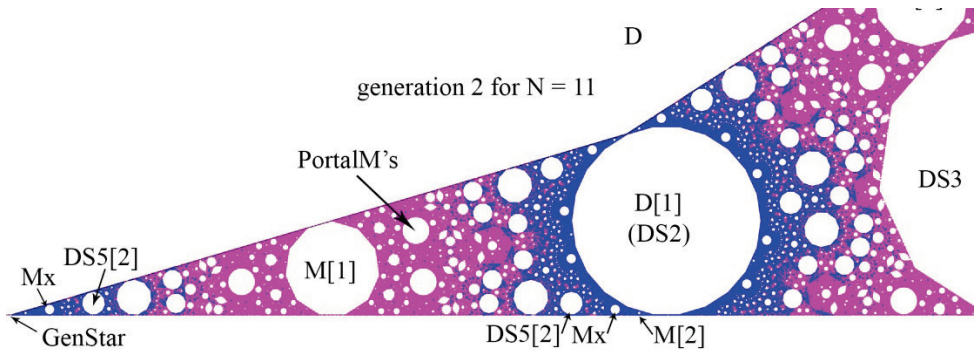

DS3[2] does not survive and in its place is a non-canonical M-tile which we call Mx. The missing D[2] would be to the right of Mx (DS3[2] and D[2] should replicate the step-4 relationship between D[1] and DS3[1] as shown above.) To see what goes wrong with the formation of D[2], we trace the local web below in the vicinity of star[1] of D[1] below. Mx is shown for reference. These are iterations 0,110,147 and 450. The new D[2] would be a step-2 of D[1] so it would share a vertex with D[1]. M[2] is maturing normally at step-1 but there is no canonical step-2 tile. It is clear that Mx is a failed DS3 because it has links to the edges of D[1] where a DS3[2] would share a vertex with a second M[2] - which does not exist here.

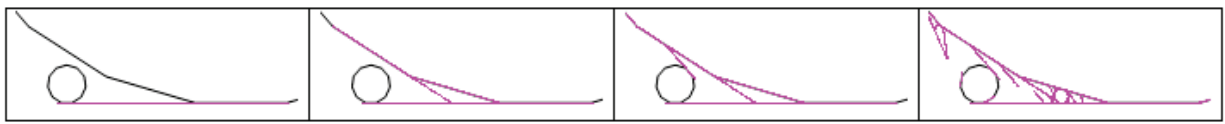

Below is the same star[1] region of D[1] for N = 13 showing the successful formation of M[2] and D[2]. These are iterations 450, 1100 and 2000. A perfect DS3[2] is also forming here on the left and it will have the correct relationship with the second M[2].

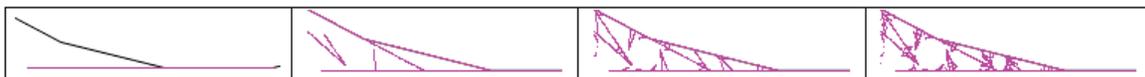

The region around Mx is enlarged below. (The lower left corner is at {-6.15937 ,-0.96645} and the region is .022 by .013). The magenta line of symmetry extends all the way from the center of M[1] to the local star[1] vertex of D[1] and it passes through the center of Mx (which can be found to high precision using the Digital Filter map which rotates interior points.) This is an 'orbit' plot using initial points which are either nonperiodic or the periods are so long that they are essentially nonperiodic. Mathematica tracks accuracy loss on all calculations and starting with 40 or more decimal places it can reliably track an orbit for 2 or 3 billion terms. Choosing initial points on a forward edge of N = 11 would guarantee that the orbits are not periodic, but 'almost all' such orbits terminate at a trailing edge. With no obvious self-similarity, it is a challenge to find true nonperiodic orbits to illuminate regions such as this.

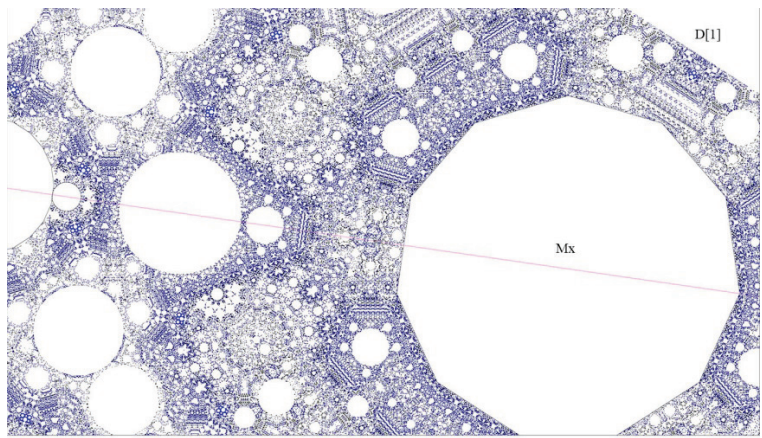

As indicated earlier, one measure of the complexity of the dynamics in any region, is to compute the percentage of web points which 'survive' on each iteration. These are the points which do not map to trailing edges. For 'well-behaved' regions this percentage would be low because it implies formation of 'stable' tiles. For this region, the percentage of survivors is very high – indicating a large population of nonperiodic points and a very complex residual set.

There are small invariant islands in this region – which are just visible on the plot above. One of these regions is enlarged below. These regions have fractal boundaries and their own unique dynamics. It should be clear from their shape that they originate near GenStar. Their scale is comparable to GenScale$^3$ ≈ .000075 but the tiles are not canonical in scale. Many of these islands have significant alignments. The two islands shown on the left below have local GenStar convergence points which lie on the edges of the small tile. This would imply that the edges of this tile may have local dynamics which replicate the dynamics at GenStar. Because N = 11 is a 4k+3 prime, very little is known about the GenStar dynamics.

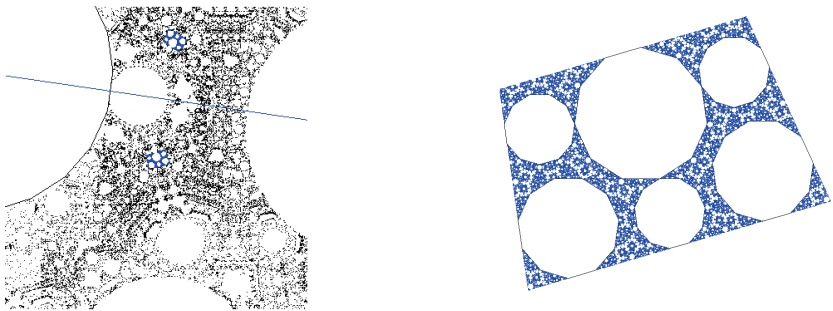

## Section 5 - Interaction of scales - dual roles

N = 7 is the first regular polygon with two incommensurate scales. This section will describe some of the ways in which these scales interact. This will give a hint about what to expect for larger N. The First Family for N = 7 is shown below.

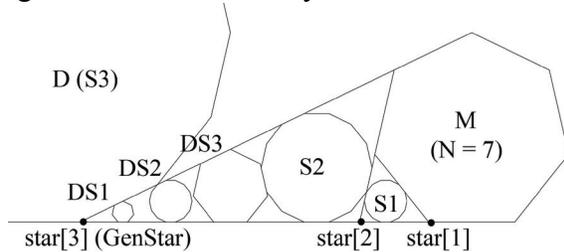

| Scale[1] | 1.000000 |
|---|---|
| Scale[2] | 0.3840429 |
| Scale[3] (GenScale) | 0.1099162 |

| Tile | radius | | Tile | radius |
|---|---|---|---|---|
| M | 1.0 (rM) | | DS1 (M[1]) | rM*GenScale |
| S1 | rD*GenScale | | DS2 (D[1]) | rD*GenScale |
| S2 | rD/scale[2] | | DS3 | rM*scale[2] |
| S3 (D) | 1.9498558 (rD) | | | |

Note that DS1 and DS3 are both scaled copies of M. For N odd there will always be a canonical M-type tile with radius equal to each scale – but of the First Family, only M-D and M[1]-D[1] will have the canonical GenScale relationship - which sometimes continues at GenStar. Since S1 is the same size as D[1], it can also foster canonical next-generation tiles. The S1 tile is representative of the scale[3] family in the same way that S2 represents scale[2]- so each step-region of M can be associated with a scale and it is safe to say that all of the dynamics of N = 7 are determined by the interaction of these scales. One such example will be presented below.

As generations evolve, any tile can foster 'families' of self-similar tiles so families can exist on all scales – but the basic families for N = 7 are the S1 family and the S2 family. In mixed generation plots it is sometimes impossible (or unnecessary) to tell the families apart because the only important issue in terms of dynamics is the relative scales. For example the scenario on the left below is from the first generation. This is sometimes called a 'short' family because DS3 is an M-type tile – but if DS3 was an M tile, the corresponding D's would have radius rD*scale[2] as shown in magenta, and these virtual D's are not compatible with M[1] – so the local dynamics would be unpredictable if DS3 was a step-3 and an M at the same time. At the top of S1 and at star[2], the roles are reversed and DS3 is an M[2] surrounded by real D[2]'s. Now the D's from the first plot are real or virtual S2[2]'s and the relative geometry is unchanged.

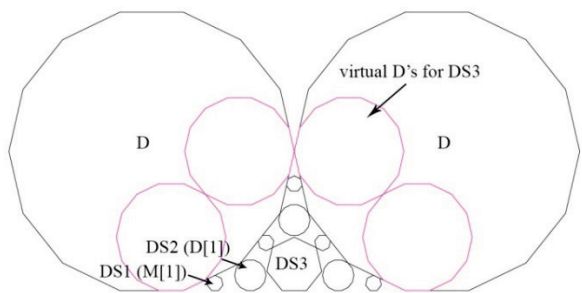
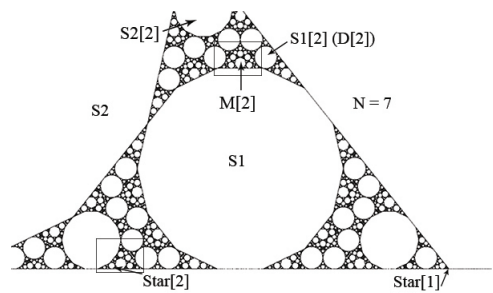

In the outer star region, DS3 is a step-3 of D and there is no conflicting role. The last section showed the important part that DS3 played in fostering the 2$^{nd}$ generation. The dynamics here could be resolved in the same way, if there were no S2[2]'s to create a conflicting role for M[2]. But the entire inner star region is dominated by the interaction of S1 and S2 dynamics. The conflict is most obvious at star[2] where M[2] has an S2[2] on the left and D[2]'s on the right. However the small scale dynamics at the top of S1 are conjugate to the dynamics at star[2] - so it would be safe to say that virtual S2[2]'s exist at the top of S[1] and virtual D[2]'s will play the dual role at star[2]. (Star[1] also has small-scale dynamics conjugate to star[2].)

Below are the two regions side by side. In both cases the dynamics are very 'directional' but there is a perfect match on the left side of M[2] - because at star[2], there is a matching virtual D[2] on the left (which is not shown here). So in both of these plots, there is a half-ring of four D[2]'s but on the right, two of these D[2]'s are virtual.

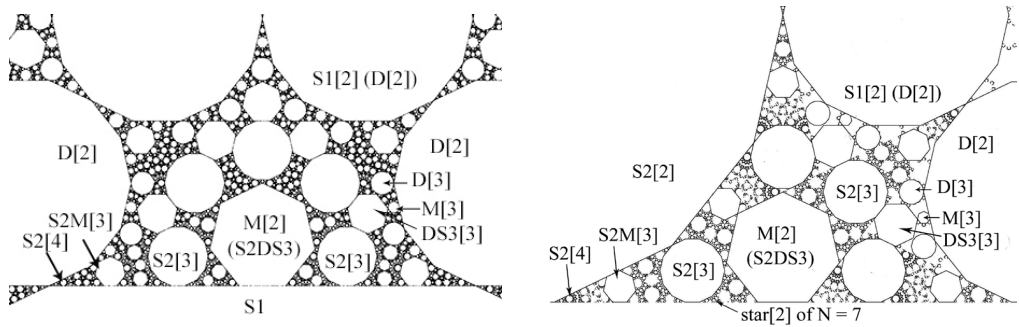

The central M[2] tile has radius GenScale$^2$ and this matches the step-3 of S2[2], which has radius (GenScale$^2$/scale[2])·scale[2]). These are the same numerically but not dynamically. The region to the left of M[2] is what would be expected of an S2-scaled DS3 – not an M[2]. It is a 'short' family - which includes a 4$^{th}$ generation with M-D pair S2M[3] and S2[3]. As with all 4$^{th}$ generations, it is dominated by PortalM's. This family structure continues past S2[4] which will be patriarch of a perfect S2-scaled 5$^{th}$ generation. However there is a break-down in dynamics at the rectangle because S2M[3] is not compatible with D[2]. This mismatch can be seen clearly on the first plot above where S2M[3] is not in any of the step regions of D[2].

On the star[2] plot, everything looks normal because S2M[3] is in the canonical step-1 position relative to S2[2], but the virtual D[2] still exerts disruptive influence – so both plots have identical irregular dynamics at the rectangles. However these irregular dynamics are confined to a small region and after S2[4], the dynamics return to normal. The rectangles should contain a PortalM, and the irregular dynamics can be regarded as the debris from this missing S2PM[4]- which is itself a failed S1 of S2DS3. That is why star[1] has similar issues.

The star[2] point has 'normal' dynamics with sequences of S2's alternating real and virtual in a fashion similar to star[3] of D[1] – so S2[4] at star[2] will be virtual. On the right side of star[2] the convergence also mimics that seen at star[3] of D[1] – with DS3's and matching PM's.

Whenever there are at least two incommensurate scales, there is the potential for an endless number of sub-scales as generations interact. Therefore it would be natural to expect a continuous spectrum of scales of the type observed by Lowenstein et al. for rotations by π/7.

**Section 6 – the global picture**

Case 1 – N odd

Below are iterations 0 – 4 of the inverse web in blue for N = 7, along with the level-1 forward web in magenta for reference. The evolution of the primary step-3 D tile was discussed in Section 2. This D and the adjacent D's form what we call Ring 0. Our primary interest here is the formation of the D's in Ring 1 and subsequent rings.

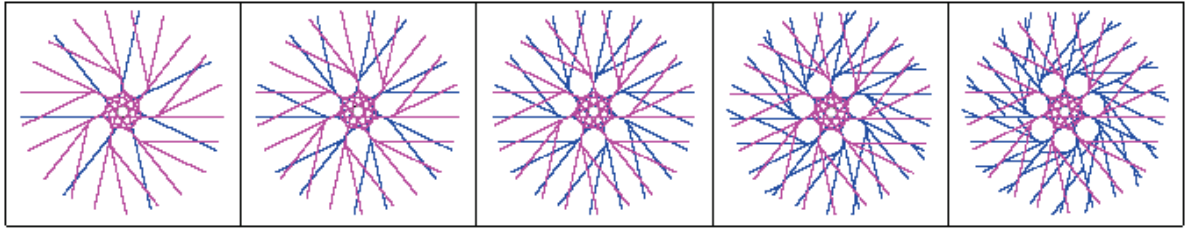

The enlargement below shows the level-6 (inverse) web where the primary D is almost a stable tile. Combining the inverse and forward webs 'speeds' up the evolution, but we retain the distinction because it provides insight into the web development. The plot below focuses on the region between the horizontal forward edge and adjacent forward edge. This region defines one of 7 congruent 'subdomains' of $\tau^{-1}$. In Section 2 we defined the 'star points' to be the intersection of the level-0 trailing edges with this horizontal forward edge. For N = 7, these star points define three 'step-regions'. The outer-most star point is the (level-0) GenStar point. The web develops recursively so every iteration of the forward web generates a new GenStar point and this defines a new ring of D's. When N is even, the spacing of all these points is identical, but here the first gap is essentially 'doubled' by the transition for N = 7 to N = 14, so the horizontal spacing between D's is $|2(GenStar[[1]] - s/2)| \approx 8.76257$.

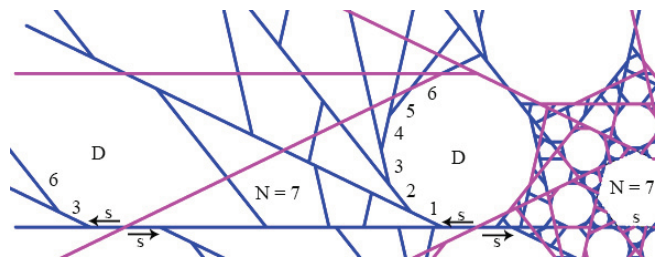

The level-1 GenStar point evolves in a fashion similar to D, but the edge progression is now step-3, so the step sequence of these D's is {3,3,4} compared with constant {3} for Ring 0. Each new ring adds another {3,4} to the sequence. For N odd the number of D's in ring k is always N(2k+1). This is also the period of the centers, so the periods are odd multiples of N. The D's in these rings are natural extensions of the canonical S[k] tiles – starting with the first ring of D's which are S[[N/2]].

Below are rings 0,1 and 2 for N = 7 with periods of 7, 21 and 35.

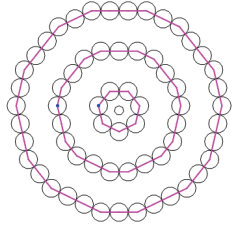

The magenta arcs connect the D centers. These magenta polygons are actually the orbits of the centers under the *return* map $\tau^2$, so under this map, each D center maps to its neighbor – clockwise here. Ring 0 is the only ring that can be replicated by simple rotation. The remaining rings form non-regular polygons. These non-regular polygons are of interest in their own right and we will study their dynamics in Appendix G.

Since the D's in Ring 0 have constant step sequence $\lfloor N/2 \rfloor$, their winding number is $\lfloor N/2 \rfloor/N$. This is an upper bound for the points inside this ring (the 'star' region). It is not clear whether there are points inside the star region with winding numbers arbitrarily close to D. All points outside of this star region have step sequences composed of just $\lfloor N/2 \rfloor$ and $\lfloor N/2 \rfloor + 1$. Clearly $\lfloor N/2 \rfloor + 1$ is the maximum possible step for a regular polygon with odd number of sides and the limiting horizon step-sequence for N odd is $\{\lfloor N/2 \rfloor, \lfloor N/2 \rfloor + 1\}$.

In any ring, the D's step sequences serve as bounds for the remaining points – just as they do in the star region. For example with N = 7, the upper bound for all points in the interior of the star region is 3/7. The D's in Ring 1 have period 21 orbits and the step sequences are constant $\{3,3,4\}$ so $\omega = 10/21$. This is an upper bound on $\omega$ for all points in the interior of this ring.

| Step sequences and periods for D's for N = 7 | | | |
|---|---|---|---|
| Ring | Step Sequence | Period of center | Winding Number - $\omega$ |
| 0 | {3} | 7 | 3/7 |
| 1 | {3,3,4} | 21 | 10/21 |
| 2 | {3,3,4,3,4} | 35 | 17/35 |
| k | {3} + k*{3,4} | 7(2k+1) | (3(k+1) + 4k)/7(2k+1) → 1/2 |

Each inter-ring region can evolve independently as shown below for N = 7

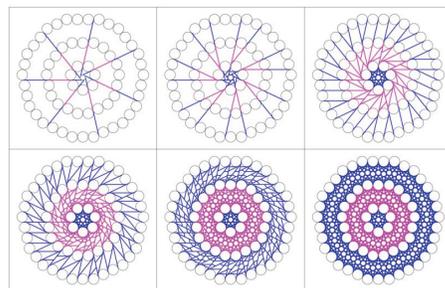

To study the global properties of the Tangent map for odd N-gons, it is helpful to use the return map $\tau^2$ and the related Pinwheel map pioneered by R. Schwartz [S2]. The Pinwheel map filters orbits to lie only on specific 'strips'. We will illustrate the Pinwheel map for N = 7 where there are 7 possible strips. They are shown below superimposed on the first three rings of D's. The strips define transitions in the dynamics of the return map, $\tau^2$. For N = 7, there are only 14 possible displacements given by the vectors V1,V2,..,V7 and their negatives. For example V1 = 2*(v1-v4) where v1 and v4 are the respective vertices of N = 7, so p2 = p1 + 2*V1 and p3= p2 + 3*V2. These are 'accelerated' orbits which can be used to analyze large-scale dynamics. (The point p1 is an off-center point in one of the D's, so it has period 70, but under the return map it is period 35 as shown here.) Under the Pinwheel map restricted to the primary strip, p1 is a fixed point, so the filter ratio is 35 to1 but no information is lost.

Recall that the winding number of each D in Ring 2 is 17/35 and for off-center points this translates into 34/70. These 34 steps are shown here. At a 'reasonable' distance from the origin (after Ring 1) any point in the interior of Ring 2 must visit each of the 14 regions on each transversal. In each region, the count is less than or equal to the count for Ring 2, so it can never exceed 34. This is an example of a 'winding number' trap - similar to the one defined by the first ring of D's – the inner 'star' region.

The green orbit shown below is typical. Only the first transversal is shown here. It visits each region at least once before returning to the primary strip. The total number of steps is 30. This orbit is period 182 so the return orbit is period 91 and the winding number is 616/(7*182).

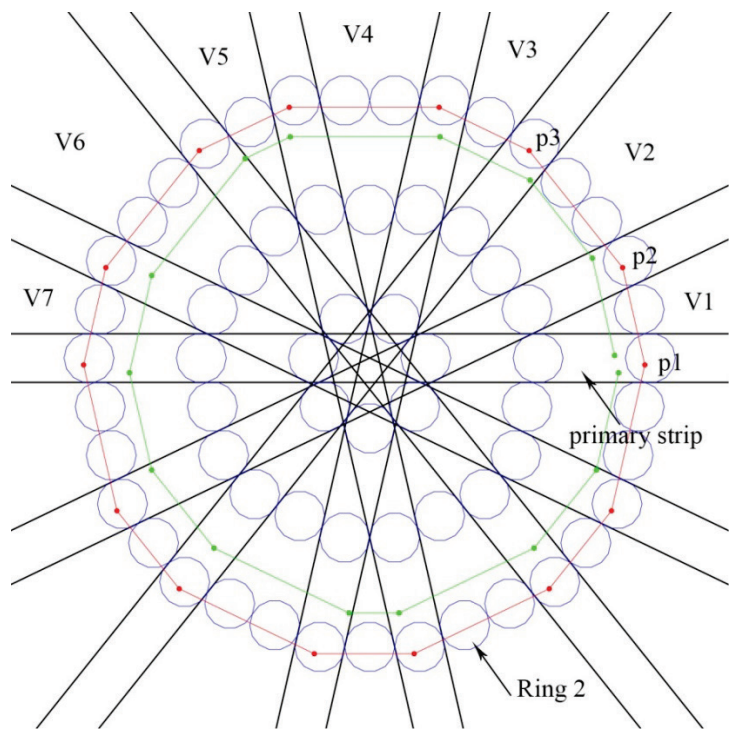

**Case 2: N twice odd**

The ring structure for N even is similar to the N odd case, except that the rings now form regular polygons. The Pinwheel map has no natural extension to the case of N even (or to the case of two parallel sides for N non-regular). To understand the new issues which arise for N-even we begin with the twice-odd case.

The first three rings for N = 14 are shown below. They have 14, 28 and 42 D's and ring k will have 14(k+1) D's. The ring centers are just multiples of the center of D so cRingk = cS[6](k+1)

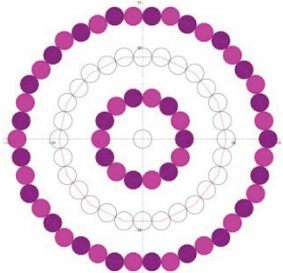

The 'even' rings 'decompose' into two distinct cycles as shown above. The Ring 0 decomposition mimics the natural factoring of N = 14 into two N = 7 heptagons. It is this factoring of N = 14 which causes the decomposition seen here. Each of the D's in the even rings, only 'see' one of the factor heptagons of N = 14 and this cuts their period in half, so ring 0 has two orbits of period 7, and ring 2 has two orbits of period 21. The 'odd' rings do not decompose, so in ring 1, each D has the full period 28 because the step sequence is {6,7} which is odd.

This complicates the issue of tracking orbits using the pinwheel 'strips'. For N = 14, instead of using $\tau^2$, it would be necessary to use $\tau^7$ for ring 0, $\tau^{21}$ for ring 1 and $\tau^{35}$ for ring 2.

Decomposition of this nature is common for regular N-gons when N is composite. For N = 14, it affects the orbits of other First Family members as well – all the S[k] for k even experience the same decomposition – the D's in even rings are natural extensions of the even S[k]. Appendix G discusses the general issue of 'decomposition'. Note that the step sequences in the table below alternate even and odd. When N is even, the primary D tile is S[N/2 -1] so for N twice-odd they will start out even but subsequent rings add N/2 – which is odd. For N twice-even, they start out odd and remain odd because N/2 is even.

| Step sequences and periods for D's for N = 14 | | | |
|---|---|---|---|
| Ring | Step Sequence | Periods of centers | Winding Number - ω |
| 0 | {6} | 7 & 7 | 6/14 |
| 1 | {6,7} | 28 | 13/28 |
| 2 | {6,7,7} | 21 & 21 | 20/42 |
| 3 | {6,7,7,7} | 56 | 27/56 |
| k even | {6} + k*{7} | 14(k+1)/2 & 14(k+1)/2 | (6 + 7k)/14(k+1) → 1/2 |
| k odd | {6} + k*{7} | 14(k+1) | (6 + 7k)/14(k+1) → 1/2 |

**Case 3: N twice even**

Below are the first three rings for N = 16. As in the twice-odd case, the centers of the D's form regular N-gons and Ring k has N(k+1) D's, so the rings below have 16, 32 and 48 D's and these are the periods. The physical ring spacing is the same as the twice-odd case so cRingk = cD*(k+1) and the horizontal spacing for N = 16 is -(GenStar[[1]] + s0/2) ≈ 5.22625.

Unlike the twice-odd case, the rings do not decompose because they are all odd-step – so the orbits visit every vertex. The first three rings are shown below with step sequences {7}, {7,8} and {7,8,8}.

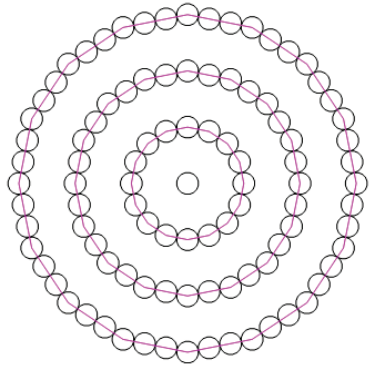

As with the twice-odd case, the centers no longer map to neighbors under the return map $\tau^2$, so to obtain the magenta polygons above, it was necessary to skip 7, 15 and 23 iterations respectively. Even though these rings do not decompose, we can mimic this decomposition by taking just the even or odd vertices. For ring 2, this yields two congruent non-regular 24-gons. (For N = 14, on the previous page, ring2 decomposes naturally into two congruent non-regular 21-gons – which are colored accordingly.) When placed at the origin, these 'Ring2' polygons have interesting dynamics. In the case of N = 16, the matching (3N/2)-gon has three regular octagons as factors, but one is slightly larger. To the astute observer, this should be evident from the ring structure above. The class of non-regular polygons which arise from rings of D's for N-even, are called Ring2 polygons. All the rings collapse by co-linearity to ring 1 so the Ring2 polygons are unique for a given N. See Appendix G.

| Step sequences and periods for D's for N = 16 | | | |
|---|---|---|---|
| Ring | Step Sequence | Period of center | Winding Number - ω |
| 0 | {7} | 16 | 1/16 |
| 1 | {7,8} | 32 | 15/32 |
| 2 | {7,8,8} | 48 | 23/48 |
| k | {7} + k*{8} | 16(k+1) | (7 + 8k)/16(k+1) → 1/2 |

**Rings composed of polygons with odd number of sides**

All the rings shown above are constructed from regular polygons with an even number of sides, but bounding rings can also be formed using M-type polygons with an odd number of sides.

When N is odd or twice odd, there will be unbroken rings of 'M's', which create secondary invariant regions between the rings of D's. Polygons in rings are typically odd-step relative to their neighbors so they are inverted. For the D's this makes no difference, but for polygons with odd number of edges, the centers may be displaced enough to create non-convex polygons.

Example 1: N = 7 on the left and N = 14 on the right

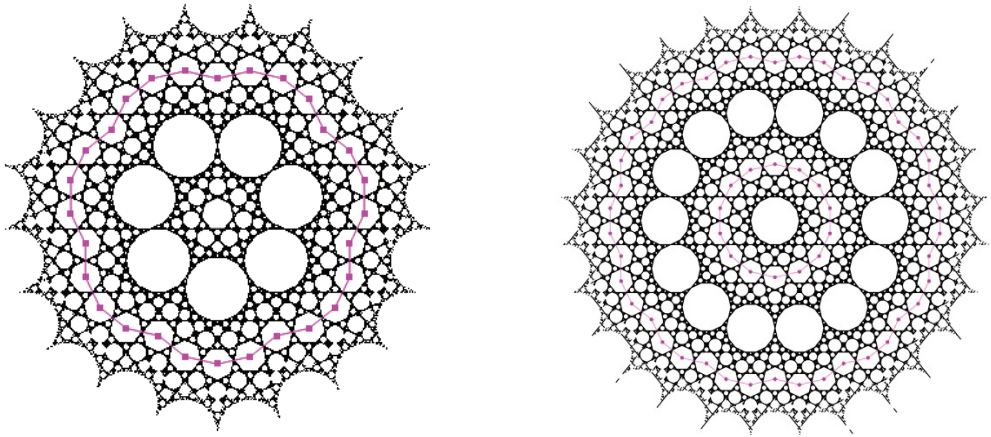

These rings of M's never decompose so the period of the N = 7 ring is 28. The step sequence is {3,3,3,4}, so the adjacent M's are step 13 from each other. In all cases the ring spacing is the same as the D spacing.

**Appendix A** - Open Questions

The first three questions below are from Richard Schwartz's Wikipedia article on Outer billiards. It is safe to say that what is not known about the dynamics of outer billiards far exceeds what is known.

(i) *Show that outer billiards relative to almost every convex polygon has unbounded orbits. (We present a possible counter-example in Appendix G.)*

(ii) *Show that outer billiards relative to a regular polygon has almost every orbit periodic. The cases of the equilateral triangle and the square are trivial, and S. Tabachnikov answered this for the regular pentagon.*

(iii) *More broadly, characterize the structure of the set of periodic orbits relative to the typical convex polygon.*

(iv) *If the algebraic complexity of a regular N-gon is defined to be $\varphi(N)/2$, then it appears that regular polygons with quadratic complexity or higher have structure on all scales. There are also non-regular polygons which have non-trivial webs. Is there a well-defined class of polygons with structure on all scales ? If a non-regular polygon has a 'factor' polygon (Appendix G) which is quadratic or higher, does it always have structure on all scales ?*

(v) *For a polygon with singularity set W, the 'residual set' $\overline{W}/W$ contains all the non-trivial dynamics, so it would be natural to call this the Julia set of the polygon. What is the topology of this set ? In particular does it have zero Lebesgue measure and what is the spectrum of dimensions ? Is the maximal Housdorff dimension strictly greater than 1 for regular polygons with at least quadratic complexity?*

(vi) *Every regular n-gon has a corresponding number field which is an algebraic extension of the rationals $\mathbb{Q}$. This number field is the cyclotomic field $K_n = \mathbb{Q}(z)$ where z is a primitive nth root of unity. To what degree does this field determine the dynamics of the polygon ?*

(vii) *Gutkin and Tabachnikov showed that outer billiards for any regular N-gon has polynomial complexity with maximal degree $\varphi(N) + 2$. Is the actual language complexity $\varphi(N)$ ? They show that this is true for rational N-gons.*

(viii) *All regular N-gons have a 'natural' scaling given by GenScale[N]. For the 'quadratic' polygons, N = 5, 8 and 12, the corresponding temporal growth rate is known and together they determine the fractal dimension. What part does GenScale[N] play in cubic cases and beyond where there are most likely a spectrum of scales ?*

(ix) *Because of the conjugacy between regular N-gons and 2N-gons for N odd, the dynamics of all regular polygons can be reduced to the dynamics of 2N-gons. There also appears to be a conjugacy between the dynamics of any regular 2N-gon and certain piecewise rational rotations acting on a torus. If the exact nature of the conjugacy could be determined, it could greatly simplify the study of the dynamics of regular polygons.*

*(x) What is the relationship between the π/5 two-triangle case and π/7 three-triangle case studied by Goetz, Lowenstein et al. and the corresponding regular polygons N = 5 and N = 7 ? Are three triangles sufficient for reproducing the non-trivial dynamics of any regular polygon ?*

*(xi) The web for a convex polygon generates 'tiles' which are themselves convex polygons. What is the relationship between the <u>in situ</u> dynamics of a tile and the <u>in vitro</u> dynamics – where the tile is the generating polygon ? For the maximal D tiles, the dynamics seem to be conjugate.*

*(xii) For regular polygons what is the relationship between geometric scaling and dynamics. Can this relationship be made explicit and does it persist for non-regular polygons ?*

*(xiii) Why is it true for 4k + 1 prime N-gons that generations of D's and M's converge to the GenStar point with geometric scaling GenScale[N] and temporal scaling N + 1 ?*

*(xiv) For regular polygons the stable tiles have measure which is bounded above by the D's but the lower bound is zero. The number of edges is also bounded by 2N for N odd and N for N even. Does the diversity of tiles and the complexity of the dynamics increase (within bounds) as the scale diminishes ?*

*(xv) For the regular hendecagon, N = 11, tiles of the $25^{th}$ generation would be smaller than the Plank scale of $3 \cdot 10^{-31}$m. What is the web structure at this scale ? Does the concept of 'generation' (determined by GenScale[11]) have any meaning for N = 11 ?*

*(xvi) In the Digital Filter map (Df), the 'twist' $\rho = 2\pi\theta$ lies in the range (0,1/4]. When $\rho = p/q$ with p and q relatively prime integers, the resulting Df dynamics are modeled by a regular q-gon with a step-p web when q is even or a regular 2q-gon with a step-2p web when q is odd. Therefore every regular 2k-gon, has [k/2] possible webs but some of these may collapse. When k is even, the maximal step-k/2 yields the null web of N = 4, so the effective 'maximal' step is k/2-1. In the Df Theorem in Appendix F, we identify just two of the possible webs when k is odd or when k is even and divisible by 4. In both cases the step-1 web is conjugate to the Tangent map web and the secondary web is the 'maximal' step. These secondary webs are also related to the Tangent map web, so a regular polygon like N = 24 has step-1 and step-5 webs which are 'predictable' but the 2,3,and 4 webs are not recognizable. All of these webs have the same algebraic complexity since they share the same minimal polynomial - but how are they related ? (Since rational rotations are 'dense' in the Df map, the regular polygons (with adjustable step-size) are also 'dense' – so understanding the dynamical progression of these step-sizes would be valuable.*

*(xvii) For the Tangent map with regular polygons, the rotations are always rational except for the limiting case as N→ ∞. For the Df map with $\theta = 2\pi\rho$, the 'twist' $\rho$ can be anywhere in the range of (0,1/4] and very little is known about the irrational rotations. An example using $\theta$ = ArcCos(3/4) was studied by Ashwin in [A] (2001) and he conjectured that it had a residual set with positive Lebesgue measure, but our analysis does not support this conjecture.(See Appendix F.) However it does seem likely that 'most' residual sets have a non-trivial spectrum of dimensions and this spectrum may vary continuously with $\rho$. Does the maximal Housdorff dimension approach 2 as $\rho \to 0$ and $N \to \infty$ ?*

# Appendix B – Cyclotomics and algebraic complexity

The algebraic structure of regular polygons has been a topic of great interest since Gauss laid the framework in *Disquisitiones Arithmeticae* of 1801. Since the vertices of a regular n-gon with radius 1 satisfy the *cyclotomic equation* $z^n = 1$, they are algebraic and not transcendental, so Gauss was able to apply his knowledge of modular arithmetic to cyclotomic theory.

In the case of N = 17, Gauss found a 'natural' grouping of the 16 non-trivial vertices to create a chain of nested equations of degree 16, 8, 4 and 2 - and this ensured constructability.

In modern Galois Theory, Gauss' chain of nested polynomials for N = 17 is replaced by a sequence of nested normal subgroups of the Galois group G – which is the set of automorphisms of the extension field $\mathbb{Q}(z)$ - where z is any primitive root of $z^{17} = 1$. Every automorphism of $\mathbb{Q}(z)$ must map these roots to each other, so they are of the form $\sigma(z) = z^k$ for k = 1,..,16. Therefore the Galois group G is a cyclic group of order 16 which is often written $C_{16}$. Each subsequent extension defines a subgroup of G, and if these nested subgroups are normal and the chain ends in {1}, the corresponding *cyclotomic polynomial* $\Phi(n)$ is solvable by radicals. This condition is guaranteed for all regular N-gons because G is $C_{\varphi(N)}$ and cyclic groups always have solvable chains.

**Definition**: The nth cyclotomic polynomial is $\Phi_n(x) = \prod_{k=1}^{n}(x - z_k)$ ($z_k$ primitive)

**Example**: For n = 6 shown below, $z_k = \cos(2\pi k/6) + i\sin(2\pi k/6)$. The only primitive roots are $z_1$ and $z_5$, so the cyclotomic polynomial is $\Phi_6(x) = (x - z_1)(x - z_5) = x^2 - x + 1$ since $z_1 + z_5 = 1$ and $z_1 z_5 = 1$. This cyclotomic polynomial is always degree $\varphi(n)$ and irreducible, but note that the minimal polynomial for $\cos(2\pi/6)$ is just 2x-1. By symmetry, the vertices of a regular polygon are always determined by a polynomial of minimal degree $\varphi(n)/2$

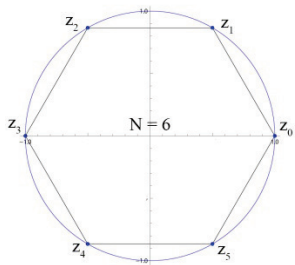

It is not trivial to prove that the cyclotomic polynomials are always irreducible. Gauss was one of the first mathematicians to recognized the importance of working with irreducible equations, although he did not use this terminology. In *D.A.* article 341 he shows that $\Phi_n(x)$ is irreducible when n is prime but he was frustrated by his inability to show that these cyclotomic polynomials are irreducible for all n, because without that proof, he could not show that his 'Fermat prime' condition for constructability was necessary. He stated that he had a proof but it was never published and it remains a matter of conjecture whether he had a proof or not. For a modern proof see section 8.4 of van der Warden's *Algebra*.

*The Nth cyclic polynomial $\Phi_N(x)$ is degree $\varphi(N)$. It is irreducible, but by symmetry, the degree of the minimal polynomial for $\cos(\pi/N)$ is $\varphi(N)/2$. Therefore the vertices of a regular polygon have 'algebraic complexity' $\varphi(N)/2$ and there appears to be a close relationship between this complexity and the dynamics of $\tau$.*

| Polygon (N) | 3 | 4 | 5 | 6 | 7 | 8 | 9 | 10 | 11 | 12 | 13 | 14 | 15 | 16 | 17 | 18 | 19 | 20 |
|---|---|---|---|---|---|---|---|---|---|---|---|---|---|---|---|---|---|---|
| $\varphi(N)/2$ | 1 | 1 | 2 | 1 | 3 | 2 | 3 | 2 | 5 | 2 | 6 | 3 | 4 | 4 | 8 | 3 | 9 | 4 |

Examples: For N odd, $\varphi(2N) = \varphi(N)$ and in fact $\Phi_{2N}(x) = \Phi_N(-x)$. This is consistent with the conjugacy in dynamics between N and 2N. For N even, $\varphi(2N) = 2\varphi(N)$, but members of this family are still related by the fact that 2N has an embedded copy of N, so there are typically vestiges of dynamics of N which survive in 2N. For example N = 16 has prominent square tiles – but the remaining dynamics appear to be unique and unrelated to N = 8 or N = 4.

N = 3, 4 and 6 are affinely equivalent to lattice polygons, so their web structure is bounded away from zero, and this implies that all orbits are periodic. N = 5, 8, 10 and 12 are 'quadratic' and they all have a relatively simple fractal structure. N = 7, N = 9 and their counterparts N = 14 and N = 18 are 'cubic' and their dynamics show a dichotomy which can be attributed to the fact that there are two or three competing non-trivial scales. This creates a residual set which is apparently multi-fractal.

N = 11 and the companion at N = 22 are the only regular polygons whose minimal polynomial is quintic. Mathematica says that the minimal polynomial of N = 11 is $1 + 3x - 3x^2 - 4x^3 + x^4 + x^5$. The Galois group for N = 11 is isomorphic to $C_{10}$, the cyclic group on 10 elements. As with all cyclic groups, $C_{10}$ is solvable, but it has a short chain because the only non-trivial subgroup is $C_5$ which is simple. (The general case for quintics is not solvable by radicals because the Galois group is isomorphic to $S_5$ - the symmetric group on 5 elements - and one subgroup is $A_5$ – the alternating group on 5 elements – which is simple but not Abelian. Therefore one chain is $S_5 \supset A_5 \supset \{1\}$ and all other chains are equivalent by the Jordan-Holder Theorem, so the chains are not normal and cannot be used to construct quotient groups.)

As expected, the minimal polynomial for N = 11 given above, has Galois group $C_5$, so it is the second polynomial in the chain. Gauss would have found it by matching up the complex conjugate pairs $s_1 = \{z + z_{10}\}$, $s_2 = \{z_2 + z_9\}$, $s_3 = \{z_3 + z_8\}$, $s_4 = \{z_4 + z_7\}$, $s_5 = \{z_5 + z_6\}$, where $z_k = \cos(2\pi k/11) + i\sin(2\pi k/11)$. It is easy to verify that these $s_k$ are the (real) solutions to the minimal polynomial given above. The sum of the $s_k$ (trace) must be -1, but the alternating sum $s_k(-1)^k$ is equal to GenScale[N] when N is a 4k + 3 prime and -GenScale[N]-2 for 4k+ 1 primes. This gives a relationship between the scales and the minimal polynomial for prime N.

Any algebraic analysis of the dynamics can be carried out in real or complex form but for regular N-gons there are advantages to working with the minimal polynomial for $\cos(2\pi/N)$ rather than the full cyclotomic polynomial and this applies to any mapping which is based on rotations by $\theta = 2\pi k/n$. When $\theta$ has this form the minimal polynomial will be monic of degree $\varphi(n)/2$ so $\cos(\theta)$ is an algebraic integer. Therefore it is possible to do exact arithmetic within the ring $\mathbb{Z}[\cos(\theta)]$ - which is finitely generated. See Appendix F.

**Appendix C** – Symbolic dynamics and language complexity

*The algebraic complexity of the Tangent map is closely related to the 'language complexity' – that is the complexity of the corner sequences or step sequences. For regular N-gons, this complexity is polynomial of degree no greater than φ(N) + 2.*

Conservative dynamics and the related symplectic mappings do not exhibit exponential divergence of nearby orbits, so τ is not technically 'chaotic'. Any mapping with polynomial complexity cannot be ergotic. Yet the phase space and webs are similar to the Poincare cross-sections of chaotic systems which have positive Lyapunov exponents. Authors have used the term 'weak chaos' or 'pseudo-chaos' to describe the dynamics of these mappings.

One measure of the complexity of any mapping is the complexity of the mapping as a language. This is the science of symbolic dynamics which was pioneered by Jacques Hadamard, Emile Artin, George Birkhoff, Claude Shannon, Steven Smale and Yakov Sinai. Any discontinuous mapping with a divided phase space is a candidate for symbolic dynamics.

A number of authors have applied symbolic dynamics to inner billiards with some success, but the case of outer billiards has proven more difficult. In [T](1995) S. Tabachnikov applied renormalization methods and symbolic dynamics to the regular pentagon N = 5 to obtain a complete description of the dynamics. In [B] (2001), Buzzi showed that any piecewise isometry will have zero entropy – which implies that the complexity is sub-exponential. In [GT] (2006), Gutkin and Tabachnikov showed that outer billiards for any regular N-gon has polynomial complexity with maximal degree φ(N) + 2. In particular, they show that rational N-gons have quadratic complexity, and this means that in general, the upper bound of the degree might be φ(N).

**Example**: To obtain a coding sequence for the nonperiodic orbits of N = 5, any limiting sequence of D's or M's would suffice. We will use the same D[k] sequence discussed in Section 4, so we have reproduced the plot and table. In Section 4, these equations were used to obtain a recursive count of the decagons and pentagons. However these same equations apply equally to periods and to step sequences.

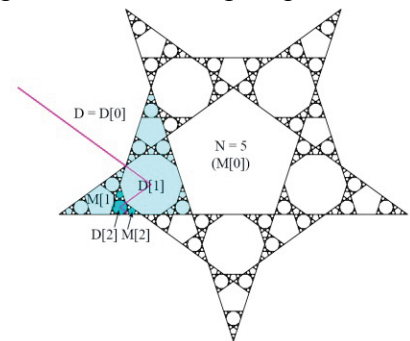

| Generation | decagons - $d_n$ | pentagons - $p_n$ |
|---|---|---|
| 1 | 1 | 2 |
| 2 | 7 = $3d_1 + 2p_1$ | 10 = $6d_1 + 2p_1$ |
| 3 | 41 = $3d_2 + 2p_2$ | 62 = $6d_2 + 2p_2$ |
| n | $d_n = 3d_{n-1} + 2p_{n-1}$ | $p_n = 6d_{n-1} + 2p_{n-1}$ |

The step sequence of D[1] is simply {1} and the step sequence of M[1] is {2,1}. According to the difference equations, D[2] will have step a sequence of the form 3*{1} ⊕ 2*{2,1} or a cyclic permutation. The actual sequence is {1,2,1,2,1,1,1}. Likewise M[2] has step sequence {1,1,1,1,2,1,2,1,1,1}. Therefore {2,1} → {1,1,1,1,2,1,2,1,1,1} which yields the generation substitution rule:   σ: 2→{1,1,1} and 1→{1,2,1,2,1,1,1}

For example, this says that the step sequence for D[3] is the third in the chain shown below:
(1} → {1,2,1,2,1,1,1} → {1,2,1,2,1,1,1} ⊕{1,1,1} ⊕{1,2,1,2,1,1,1} ⊕{1,1,1} ⊕{1,2,1,2,1,1,1} ⊕{1,2,1,2,1,1,1}⊕{1,2,1,2,1,1,1}. D[3] has period 205 and the step sequence has period 205/5 = 41 as shown here. Each iteration involves a scale reduction of GenScale[5].

Note that the region colored above is invariant under $\tau^{10}$, and each iteration inverts the triangles, so these 10 iterations can be broken down into rotations about the center of D[1] (the 2's) and rotations around D[0] ( the 1's). This was the method used in [T].

In the language of symbolic dynamics, the alphabet here is A= {1,2} and the limiting sequence u = $\lim_{n\to\infty} \sigma^n$ (1) is a called a *substitution sequence*. This sequence is a 'fixed point' of σ, since σ (u) = u .

Tabachnikov points out that σ is similar to the 'Morse' substitution which is defined on the alphabet {0,1}: $\upsilon$(0) = 01, $\upsilon$(1) = 10, so  0→01 →0110 → 01101001…

This limiting sequence is 'recurrent' because every word occurs infinitely many times. It is also nonperiodic  (note that 101010 never occurs)

**Definition**: A substitution is 'primitive' if there is an integer k such that $\sigma^k$(i) contains every symbol for each i ∈A.

**Theorem**: Every primitive substitution dynamical system is minimal and uniquely ergotic. (The associated mapping of these systems is the shift map.)

Both σ and $\upsilon$ are clearly primitive, so they are minimal (they have no proper subsystem) and uniquely ergotic (there is a unique probability measure). This in turn implies that both are uniformly recurrent because minimal is equivalent to uniformly recurrent.

A simple proof that u is nonperiodic will follow from the two facts given here - which the reader can verify: (i) the limiting winding number of  u is  ¼ and (ii) no periodic orbit for N = 5 can have a winding number of ¼.

Because u codes the decagons which are dense in $W^C$, it follows that it must lie arbitrarily close to any other nonperiodic orbit. Therefore u generates all of the residual set and there is only one invariant measure. By contrast for N = 7 there appear to be an infinite number of invariant measures.

N. Bedaride and J. Cassaigne [BC] (2011)  reproduce Tabachnikov's results in a more modern context of language theory. They obtain a complete language analysis for N = 3,4,5,6, and 10 and show that 5 and 10 have equivalent sequences.

Extending these results to the cubic cases such as N = 7 or N = 9 is far from trivial but some progress has been made recently using simplified maps – which are still very complex. These maps will be discussed below, but first here is an example of a nonperiodic orbit for N = 14.

**Example**. A nonperiodic orbit for N = 14 (or N = 7)

Because N = 7 and N = 14 have high degrees of self-similarity it is easy to find limit points by following sequences of converging cells in the same fashion as N = 5. One example is shown below. N = 14 has three basic invariant regions as shown on the left below. The inner region shown on the right is sub-divided into 14 congruent triangles such as the short family 'tower' outlined in blue.

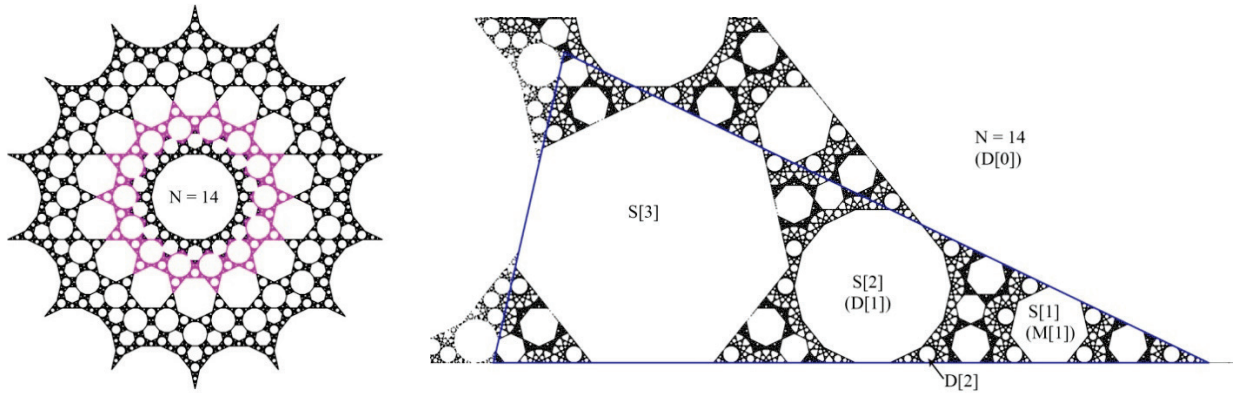

There is a sequence of D[k] cells converging to a limit point s* on the forward edge extension of N = 14 as shown in the enlargement below. The 'address' of s* is {6,3,6,3,…} where the 'buds' of N = 14 are numbered with the same convention used for N = 10 earlier - starting with the 3:00 position.

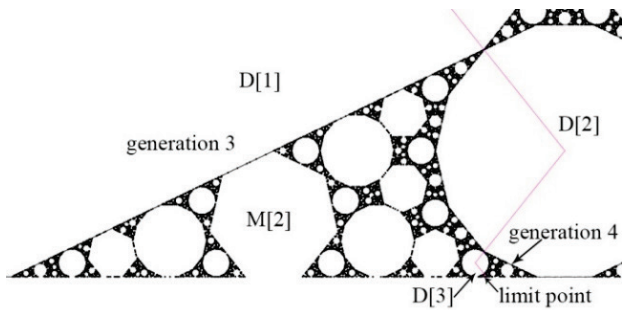

The D's (and matching M's) are scaled by GenScale[7] – and the local geometry is conjugate to the GenStar region of N = 7 via a reflection about the central S[5] of N = 14. If s* is iterated under $\tau$ the resulting non-periodic orbit will no doubt define an invariant measure but it is unlikely that the step sequence of this orbit could be described by a simple substitution sequence as for N = 5. The first 20 terms in the sequence are {3,3,3,2,1,1,1,1,2,3,3,3,3,2,1,1,2,3,3,3,2} and the 'tally' for the first 500,000 terms is {{3,179121}, {2,167200},{1,153679}} with winding number $\omega \approx .14634$.

**Appendix D – Periodicity**

Below is a summary of periodicity of orbits for regular polygons. This is a short list because very little is known about periods for regular polygons. Any tile of a regular polygon will have a bounded orbit and it is likely that the periodic orbits are dense, but this has not been proven except for a few simple cases. Lattice polygon with rational vertices will have all orbits periodic and this is true for any affine image of a lattice polygon – but this yields just N = 3, 4 and 6.

The periods of the rings of D's was discussed Section 6. For prime N-gons, the periods of the D's are odd multiples of N, so the D's in Ring k will have period N(2k+1). (By convention, these are the periods of the center and the first ring is Ring 0.) When N is twice-odd, the even rings will decompose into two orbits with period N(k+1)/2 each. The odd rings do not decompose so they will have period N(k+1). The twice-even case is the same as the twice-odd case but there is no decomposition, so the periods are N(k+1).

Appendix G discusses some of the issues arising from decomposition of orbits. Below we give formulas for periods of the canonical family members, but these formulas are only valid for prime N-gons. For composite N, the periods may be smaller. In general an S[k] tile will decompose whenever it is possible - so whenever GCD[k,N] > 1. This means it is not trivial to relate the periods of 'conjugate' pairs such as N = 7 and N = 14. They share webs but this does not imply that they share the same dynamics.

| N-gon | Periods |
|---|---|
| N:Prime | The canonical 2Ngons S[1],S[2],…, S[[N/2]] all have period N (mapping centers) The canonical D step-k has period N(N−(k+2)) . For example M[1] is D step-1 so it has period N(N−3) and D[1] is step-2 with period N(N−4). |
| N: Non prime odd | The canonical periods are the same as above, but the prime periods may be smaller due to decomposition of orbits. For example with N = 9, DS3 is step-3 of D so it should have period 9(9-5) = 36 but the prime period is 12 and there are three groups of 12 to make up the 36. |
| N: 4k+1 prime | For (4k+1) primes there appear to be M's and D's at all generation scales and the ratio of periods D[k]/D[k-1] →N+1 This limiting ratio applies equally to the chain of M's. (See example of N = 13 below) |
| N = 5 | If d[n] is the period of decagons and p[n] is the period of pentagons, then  d[n_] := (5/7)(8*6^(n-1) + (-1)^n); p[n_] := 5*( ((3*d[n]/5)+(-1)^n)/2+(-1)^(n-1)); d[Range[10]]= {5, 35, 205, 1235, 7405, 44435, 266605, 1599635, 9597805, 57586835} p[Range[10]]= {10, 50, 310, 1850, 11110, 66650, 399910, 2399450, 14396710, 86380250} It is easy to translate these values to the rings of D's so this is the only non-trivial case where all the periods are known. |

| | | |
|---|---|---|
| N = 7 | Even though 7 is a 4k+3 prime it does support infinite chains of D's and M's converging to the GenStar point, but the even and odd families have different dynamics. If 7 was a 4k+1 prime, the ratio of these periods would approach $N + 1$, but here there are two ratios which seem to be 8 and 25 (at least for the M's). The limiting overall growth rate of 200 must be the same for the D's and M's, but the individual ratios for the D's appear to be 10 and 20. Note that the M's for generation k are dependent on the D's from the previous generation | |

|  | Period | Ratios |  | Period | Ratios |
|---|---|---|---|---|---|
| M[1] | 28 |  | D[1] | 21 |  |
| M[2] | 98 | 3.5 | D[2] | 336 | 16 |
| M[3] | 2212 | 22.57 | D[3] | 4151 | 12.354 |
| M[4] | 17486 | 7.905 | D[4] | 27720 | 6.6779 |
| M[5] | 433468 | 24.789 | D[5] | 829941 | 29.940 |
| M[6] | 3482794 | 8.0347 | D[6] | 5524568 | 6.656 |
| M[7] | 86639924 | 24.876 | D[7] | 54678813 | 9.897 |
| M[8] | 696527902 | 8.0393 | D[8] | 1104888456 | 20.206 |

There are an endless number of local 'star' points which can support some type of convergence, but the only sequences studied in detail are the three primary star points and three of the secondary star points. Here is a summary of their apparent temporal ratios:

| Star point | GenStar (star[3]) | star[2] | star[1] | star[1] of M[1] | star[2] of DS3 (left and right) |
|---|---|---|---|---|---|
| Ratio of periods | 200 | 113 | 1254 | 1254 | 113 |

DS3 has two star[2] points which appear to have different dynamics. The star[2] point on the right of DS3 is on the dividing line between the inner and outer star region and the star[2] point on the left side is shared with star[3] of D[1]. However the ratios of periods (skipping one generation) is 113 for both. In both cases generations alternate real and virtual.

At the star[1] point of M[0] and M[1], the ratio of 1254 is a 'best guess' because this chain involves skipping 4 generations. This makes it very difficult to track periods – because they grow rapidly. Note that 113 and 1254 are both 1 mod 7 so there might be a version of the 4k+1 conjecture which applies to these secondary star points.

| | |
|---|---|
| N = 8,9,12 | These are not prime but they support infinite strings of D's and M's. For N = 8 and N = 9, the ratio of periods $D[k]/D[k-1] \rightarrow N+1$ just as in the 4k+ 1 prime case. For N = 12 this ratio appears to be 27. For N = 8 and N = 12 these D's are dense so they can be used to find the fractal dimension in the same fashion as N = 5. |

| N = 13 | This is the second 4k+1 prime and it supports families converging to GenStar but unlike N = 5, these families are only self-similar in a 'mod-2' fashion. |
|---|---|
| | <table><tr><th colspan="6">Six generations of M's and D's for N = 13</th></tr><tr><th></th><th>Period</th><th>Ratio</th><th></th><th>Period</th><th>Ratio</th></tr><tr><td>M[1]</td><td>10*13</td><td></td><td>D[1]</td><td>9*13</td><td></td></tr><tr><td>M[2]</td><td>182*13</td><td>18.20000</td><td>D[2]</td><td>119*13</td><td>13.22222</td></tr><tr><td>M[3]</td><td>2506*13</td><td>13.76923</td><td>D[3]</td><td>1673*13</td><td>14.05882</td></tr><tr><td>M[4]</td><td>35126*13</td><td>14.01676</td><td>D[4]</td><td>23415*13</td><td>13.99582</td></tr><tr><td>M[5]</td><td>491722*13</td><td>13.99880</td><td>D[5]</td><td>327817*13</td><td>14.00029</td></tr><tr><td>M[6]</td><td>6884150*13</td><td>14.00008</td><td>D[6]</td><td>4589431*13</td><td>13.99998</td></tr></table> |
| | The only hint of 'odd-even' generation dichotomy in the table above is the high-low alternation of ratios. Note that that the self-similarity never seems to include the first generation. As with N = 7, the M's ratios are out of sync with D's because the M's for generation k form on the sides of the D's of the previous generation. The table for N = 17 is very similar. |
| N even | The periods of the First Family members will depend on the factors of N as in the odd composite case above. However there is an extra level of symmetry here which can be used to relate the S[k] and the LS[k] periods. See Appendix H. |

**Period doubling**

When a tile is periodic, the Propensity Lemma guarantees that all points in the tile will have even periods except possibly the center. Since our convention is to use the center point of tiles to define the period, tiles with odd periods are also known as 'period doubling' tiles. There is no obvious geometric criteria that can be used to determine if a tile has period doubling, but symmetry excludes tiles with an odd number of sides, so no M-type tile can have odd period. For N = 9, the canonical S[3] tile is non-regular but the center is period 3 so it has period doubling. When N is odd the primary D's all have period doubling, but secondary D-type tiles may fail to have period doubling – for example in the chain of D's converging to GenStar, the even D[k]'s fail to have period doubling (see table above). The rings of D's for N = 14 also alternate doubling and non doubling and this alternation exists for the S[k] as well.

## Appendix E – Winding Numbers

To analyze orbits, we use a tool borrowed from classical analysis. The winding number (rotation number) of an orbit is a measure of the average rotation per iteration - on a scale from 0 to 1.

**Definition**: For any polygon N, let Ind[p,m] be the indices of the vertices visited by the first m iterations of the point p. The depth m-1 step sequence of p is the set of first differences of Ind[p,m]

**Definition**: For any N-gon, suppose a point p has step sequence $S = s_1, s_2,..$ The winding number of S is defined to be

$$\omega(S) = \frac{1}{N}[\lim_{m \to \infty} \frac{1}{m} \sum_{j=1}^{m} s_j]$$

If S is a periodic step sequence with period k then $\omega(S) = \frac{1}{N}[\frac{1}{k}\sum_{j=1}^{k} s_j]$ so in this case the winding number is 1/N times the mean number of steps in one period of the step sequence.

For a regular N-gon with N odd, the canonical tile S[k] has step sequence {k} and winding number k/N. For periodic tiles with period doubling, it is not necessary to map centers to get step sequences – all points in a periodic tile must have the same step sequence because the corner sequence is a property of tiles. Step sequences are all immune to 'decomposition' of orbits.

When N is prime, the prime period of the S[k] is always N, but otherwise it could be a factor of N. This does not affect the step sequence (or the winding number), so they are impervious to period doubling and to decomposition. On the negative side, two tiles may have the same step sequence and different dynamics.

**Examples**:

(i) For N = 7, D has step sequence {3} so his winding number is 3/7 and this is an upper bound for the star region. The members of the First Family for N odd have constant step sequences so their winding numbers form a sequence 1/N, 2/N.., [N/2]/N which approaches ½ as N increases. This is the 'horizon' winding number for a circle and it is an upper bound for any polygon. The minimum step sequence for a regular polygon is the canonical 1-step orbit with $\omega$ = 1/N.

(ii) The winding number of the dense non periodic orbit for N = 5, approaches .25. It often happens that self-similar dynamics leads to highly rational $\omega$. It is easy to prove that no periodic orbit for N = 5 can have a winding number of .25.

(iii) Winding numbers can be used to define the boundaries of invariant regions. Below are the four invariant regions for N = 11. As expected, winding numbers tend to increase with distance from the origin. The winding number of S[4] is 4/11 but DS[7] has step sequence {4,5} so the winding number is 9/22 – which is higher than any point on the S[4] side of the boundary.

There are points on the S[4] side which have 5's in their step sequences, but these 5's are always followed by at least two 4's, so the largest possible winding number on the S[4] side corresponds to the step sequence {5,4,4} and this is achieved by the green point p below. The point q is a close second at {5,4,4,4} and this is why it is a relative maximum in the 3D winding number plot below. The adjacent tile on the DS[7] side has sequence {5,5,4} so it is on high ground.

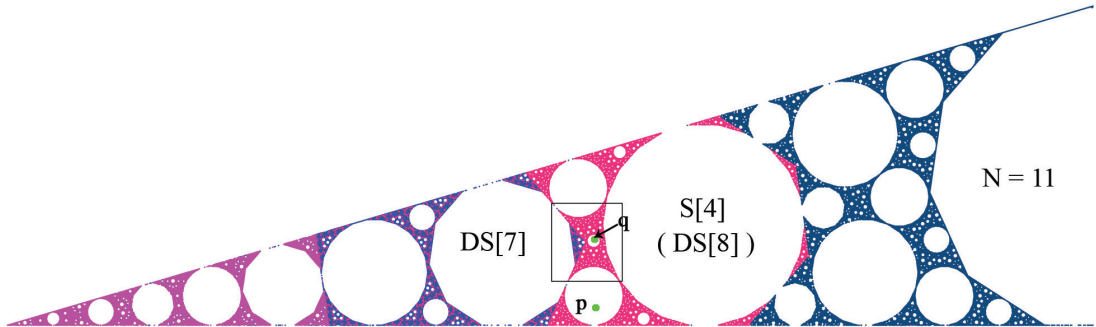

The 3D plot is a 'density' plot using the winding number as height. The view here is from S[4] looking toward DS[7]. S[4] and DS[7] are barely visible in the foreground and background because they are surrounded by points with larger winding numbers. These are sometimes called 'sticky orbits' or 'pseudo-hyperbolic' orbits because they accumulate around local star points which play the role of hyperbolic points. Not all first family members are relative minimum. D is an absolute maximum for the inner star region at ω = 5/11. It is not clear whether there are winding numbers arbitrarily close to this value.

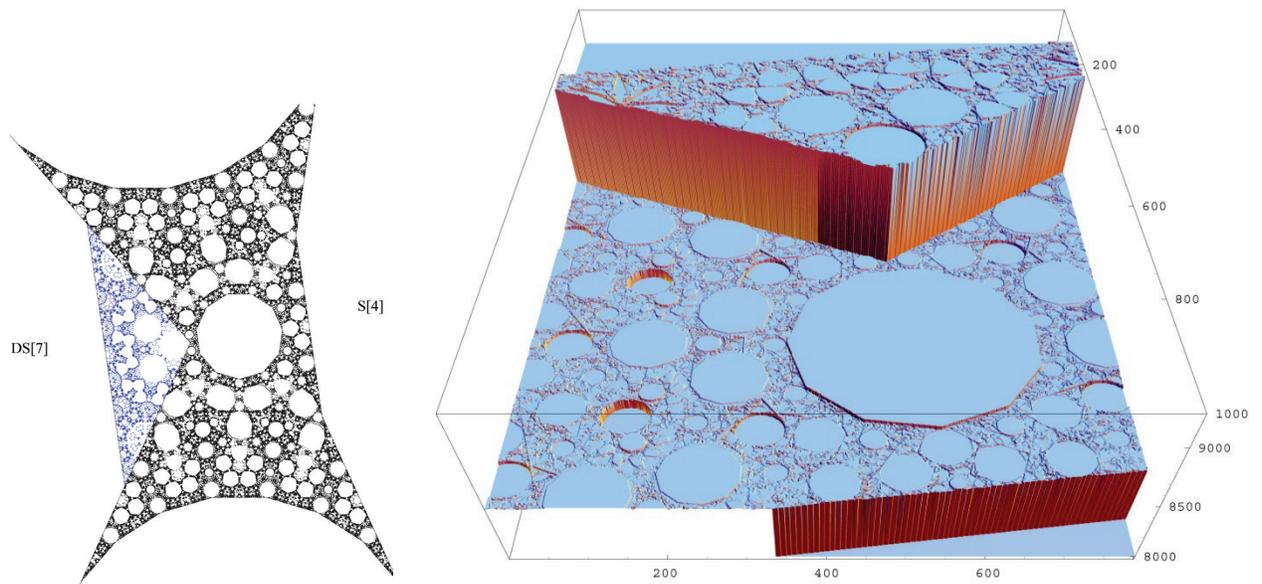

## Appendix F – Rational rotations and the Digital Filter map

For regular polygons, the Digital Filter map described below has a singular set W which is locally conjugate to the Tangent map. Since any conjugacy preserves temporal scaling and geometric scaling, it may be feasible to use maps of this kind to study topological issues such as fractal dimension and scaling for the regular polygons.

The Digital Filter map is an example of an *affine piecewise rotation.* These are two-dimensioal analogs of interval exchange maps- where rotations typically play the part of exchanges. The linear form for an affine piecewise rotation by $\theta$ can be written as:

$$T\begin{bmatrix} x \\ y \end{bmatrix} = \begin{bmatrix} 0 & -1 \\ 1 & 2\cos\theta \end{bmatrix} \begin{bmatrix} x \\ y \end{bmatrix} \quad \text{or in complex form } T[z] = \rho z \text{ where } |\rho| = 1$$

Typically the phase space X is partitioned unto a finite number of mutually disjoint 'atoms' $A_i$ and $T_i$ acts on atom $A_i$ as a rotation followed by a translation. The combined map $T: X \rightarrow X$ is defined as $T(z) = T_i(z)$ iff $z \in A_i$. For the Tangent map, each vertex of the polygon defines a $\tau_i$. When the map T is bijective, the image of a partition is a partition. When the rotation angle $\theta$ is a rational multiple of $\pi$, T is called a rational rotation or a polygonal rotation.

Recently X. Bressaud and G. Poggiaspalla [BP] used computer analysis to categorize the possible bijective polygonal piecewise isometries on two and three triangles. For two triangles there is only one case which yields non-trivial dynamics. This is what they call the 'tower' case.

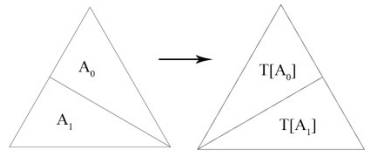

There are only two ways to implement a reflection such as this using orientation preserving rotations and A. Goetz had found them earlier. The first involves rotations by $\pi/5$ and the second $\pi/7$. We will present both of them below – starting with the $\pi/7$ case.

Example (a) This is a two-triangle map based on rotations by multiples of $\alpha = \pi/7$, with N = 14 shown for reference. $A_0$ has angles $\alpha$, $\alpha$ and $5\alpha$ and it is rotated about its center by $6\alpha$ and translated to yield $T[A_0]$. $A_1$ has angles $2\alpha$, $2\alpha$, and $3\alpha$ and it is rotated by $-2\alpha$ and translated to yield $T[A_1]$. The web shown here is obtained by iterating the x axis with T to depth 5000.

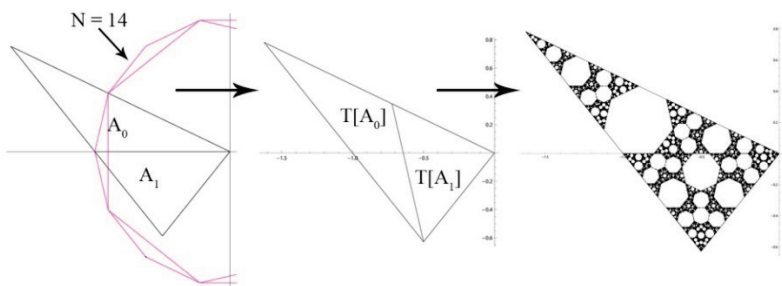

This map has not been studied in depth, but it is clear that it shares some of the characteristics of the cubic regular cases N = 7 and N = 9 – namely a limited form of self-similarity and a sequence of maximal tiles converging to the tip. The three-triangle Example (c) below is a more natural setting for analysis and this was the 'cubic' case studied by Goetz, Lowenstein et al.

Example (b) This is a two-triangle map based on rotations by $\pi/5$. This 'quadratic' example was one of the first maps studied by Goetz [G1]. It can be implemented in complex form as follows: set $a = \cos(\pi/5) + i\sin(\pi/5)$ and define triangles $A_0$ and $A_1$ with vertices $\{0, a^2 + a^4 + a^6, -1\}$ and $\{0, -1, a^6\}$ as shown below (with N = 10 for comparison). The piecewise rotations and corresponding translations are $T_0(z) := a^4 z + a^2 + a^4 + a^6$ and $T_1(z) = a^6 z + a^6$ (so the rotations are $\pm (\pi - \pi/5)$). The combined map is $T(z) = T_0(z)$ if $\text{Im}[z] > 0$, otherwise $T_1(z)$. The first iteration is shown below along with the residual set obtained by iterating a non-periodic point. (The x axis could be 'scanned' as in Example (a) above – but the results are less satisfactory.)

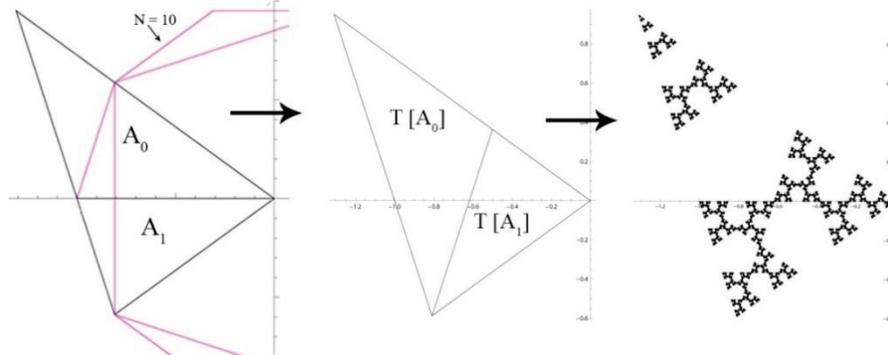

The geometric scale factor is $k = 1/\gamma$ where $\gamma$ is the golden ratio. All the limiting tiles are pentagons and their growth factor is 2, so the Hausdorff dimension is $-\text{Log}[2]/\text{Log}[k] \approx 1.4404$. To relate this to N = 5, GenScale[5] = $k^3$, and the temporal scale factor is 6, so the Hausdorff dimension for N = 5 can be written as $-(\text{Log}[2] + \text{Log}[3])/3\text{Log}[k] \approx 1.241$

Adding a third triangle yields 810 cases involving rational rotations– including 258 'cubic' cases based on rotations by $\pi/7$. One of these three-atom $\pi/7$ solutions is a simple extension of example (b). It has been studied by Goetz, Poggiaspalla, Kahng, Lowenstein, Vivaldi and Kapustov. This example is presented below. (Later in this section we will discuss a three-triangle family of maps from [BG] which seems to reproduce 'most' of the dynamics for regular 2N-gons.)

Example (c) Set $a = \cos(\pi/7) + i\sin(\pi/7)$ and define triangles $A_0$ and $A_1$ and $A_2$ with vertices $\{0, a^5-1, -1\}$, $\{0, -1, -a^3\}$ & $\{0, -a^3, -a^3+a^2\}$. The corresponding transformations are $T_0(z) = za^6 + a^5 - 1$, $T_1(z) = -az - a^4 + a^5 - a^6 - 1$ and $T_2(z) = za^6 - a^3$. Then $T[z] = T_i(z)$ iff $z \in A_i$. (Note that the rotations are $\pm (\pi - \pi/7)$, so swapping $\pi/7$ and $\pi/5$ will make example (c) identical to (b) when restricted to $A_0$ and $A_1$.) This $\pi/7$ map has a 'recursive tiling' property that enabled Lowenstein et al. to determine the scaling sequences and hence the spectrum of dimensions with maximal Hausdorff dimension $\approx 1.652$.

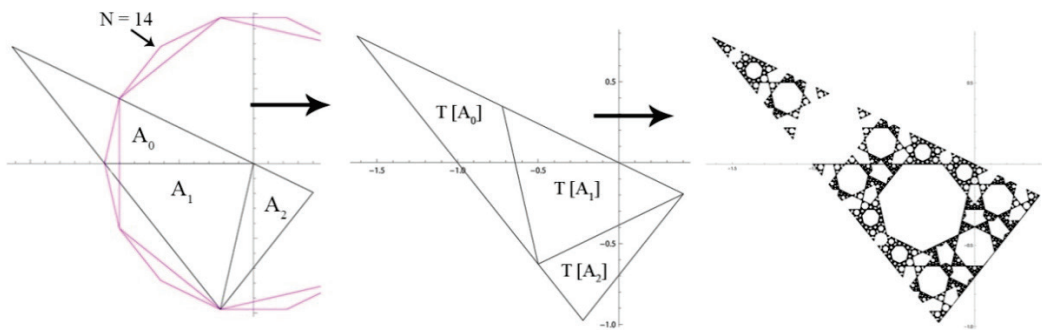

Below is a comparison of N = 5 and N = 7 alongside the respective maps from (b) and (c) above. The two residual sets for $\pi/5$ and N = 5 were generated under identical conditions with Mathematica using 100,000 points in a single nonperiodic orbit. The maps show the contrast in density between the mixed tiling of N = 5 and the purely pentagonal packing of $\pi/5$. The $\pi/7$ and N = 7 images were also generated in similar circumstances by Mathematica, but by necessity, these are 'web' plots showing $\overline{W}$. It seems likely that there will be a similar contrast in density between the limiting heptagonal tiling of $\pi/7$ and the mixed tiling of N = 7.

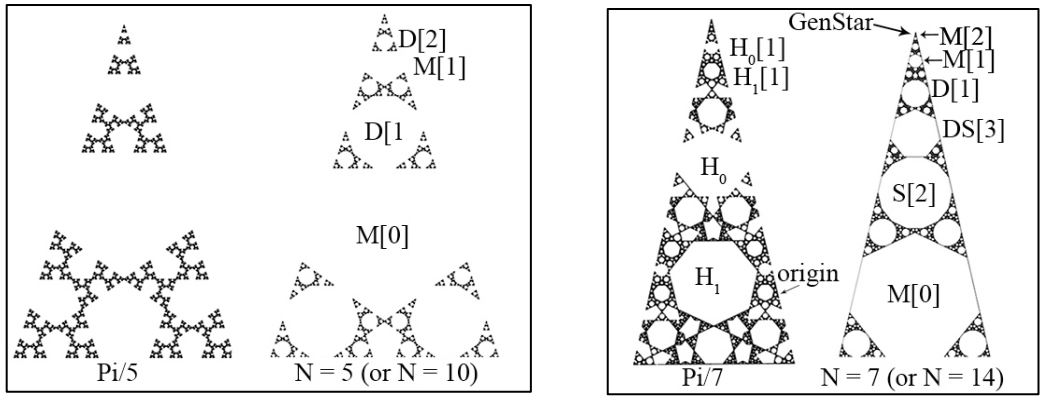

In the N = 7 'tower' shown here, the convergent point is GenStar[7] and the self-similarity skips one generation so it begins with M[2]. The geometric and temporal scaling of this M[k] sequence are apparently GenScale[7]$^2$ and 200, but there appear to be an infinite number of distinct scaling sequences, and Lowenstein et al. have shown that this applies to the $\pi/7$ case as well.

The $\pi/7$ tower is self-similar beginning with $H_1[1]$. The heptagons $H_0$ and $H_1$ get their names because their centers are fixed points of $A_0$ and $A_1$. Goetz and Poggiasapalla [GP] used symbolic complex analysis to aid in the analysis of these sequences of $H_0$'s and $H_1$'s. They determined that the sequences share a geometric scale factor of $\upsilon = 4\sin^2(\pi/14)$ (which is $\lambda$Genscale[7] where $\lambda = 2\cos(\pi/7)$ is the trace of the rotation matrix). The temporal scale factors of these sequences are $(2^{2n+1}+1)/3$ and $(4^{n+1}-1)/3$, so the temporal scaling approaches 4 in either case - one from the top and the other from the bottom.

This is similar to the scenario for typical 4k+1 primes such as N = 13 where the high and low sequences converge to the same N + 1 ratio. These sequences correspond to even and odd generations which are eventually self-similar. N = 7 also has two distinct self-similar generations with corresponding sequences of M[k]'s (and D[k]'s) converging to GenStar.

The $\pi/7$ and N = 7 cases illustrate the fact that 'quadratic-type' self-similarity can coexist with an infinite number of distinct scaling sequences. The coexistence of 'regular' and 'irregular' dynamics is a hallmark of Hamiltonian dynamics. To see some of the diversity of dynamics for the $\pi/7$ case, click on the tower above. (The 'tower' from [L] pages 32-34 is embedded in the dynamics. It can be seen in the enlargements, and of course these secondary towers exist at all scales. This occurs because the parameters of the Lowenstein tower were based on the $\pi/7$ case described above. These new parameters were designed to simplify the task of generating a 'catalog' of distinct domains. In Lowenstein's words, this example is "the smallest recursively tiled catalog found for a piecewise rotation with cubic irrational parameter.")

## The Digital Filter Map

The circuit shown below is a second order digital filter with two feedback loops. It consists of three registers with a time delay of one unit between them. The intermediate outputs y(t) and y(t+1) are multiplied by *b* and *a* and fed back in where they are added to the contents of the accumulator.

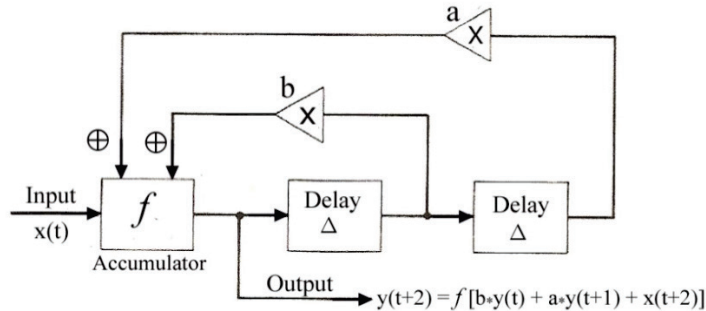

The equation for the output is $y(t+2) = f[\,by(t) + ay(t+1) + x(t+2)]$ which is a second-order difference equation. To study self-sustaining behavior, the input x(t) is set to 0. Following studies by Chua & Lin [CL] and A. Davies [D], we are interested in the self-sustaining oscillations which can occur even when there is no input, so we will assume x(t) = 0 for all t.

The corresponding second-order equation can be reduced to two first order equations by setting $x_1 = y(t)$ and $x_2 = y(t+1)$. Then at each time tick, $x_1 \to x_2$ (by definition) and $x_2 \to f[bx_1 + ax_2]$. If the accumulator has no overflow, $f$ is the identity function and the state equations are

$$\begin{bmatrix} x_1(k+1) \\ x_2(k+1) \end{bmatrix} = \begin{bmatrix} 0 & 1 \\ b & a \end{bmatrix} \begin{bmatrix} x_1(k) \\ x_2(k) \end{bmatrix}$$

This is the linear form of the Df map. The b term represents damping so when b = -1 this models a 'lossless resonator' or digital filter. For b = -1 and $a \in (-2,2)$ the eigenvalues are complex with unit absolute value so the motion is a generalized rotation with $\theta = \mathrm{ArcCos}(a/2)$. By symmetry we can restrict *a* to [0,2) so $\theta \in (0, \pi/2]$. This 'elliptical' rotation can be conjugated to a pure rotation as shown below. The conjugate rotation angle is $-\theta$, so it is clockwise.

$$\begin{bmatrix} 0 & 1 \\ -1 & 2\cos\theta \end{bmatrix} \sim \begin{bmatrix} \cos\theta & \sin\theta \\ -\sin\theta & \cos\theta \end{bmatrix}$$

The rectified linear plots will yield circles or regular polygons depending on whether the 'twist' $\rho = \theta/2\pi$ is irrational or rational. For a regular N-gon, $\rho = 1/N$ and the parameter $a = 2\cos(2\pi/N)$. willll always be an algebraic integer with minimal polynomial of degree $\varphi(N)/2$, so there is the potential for purely algebraic analysis as long as the dynamics are not rectified by the conjugacy above.

Under ideal conditions the function *f* would be the identity function, but since the registers have finite word-length, there is the issue of possible overflow. Assuming that negative numbers are stored in two's complement form, the overflow function *f* can be modeled by a sawtooth function of the form $f(z) = \text{Mod}[z+1,2]-1$ as shown here. In reality the sawtooth ramps should have as many teeth as the length of the registers.

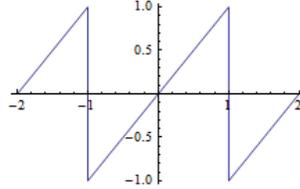

**Definition**: The Digital Filter map Df: $[-1,1)^2 \to [-1,1)^2$ is defined as

Df[{x,y}]:={y, *f* (-x + ay)} where $f(z) = \text{Mod}[z+1,2]-1$

For a given *a*, Df is a piecewise isometry with three 'atoms' which can be labeled 1 (overflow), 0 (in bounds), or -1 (underflow). The equations for these atoms are:

$$S[\{x,y\}] = \begin{cases} 1 & \text{if } -x+ay \geq 1 \\ 0 & \text{if } -1 \leq -x+ay < 1 \\ -1 & \text{if } -x+ay < -1 \end{cases}$$

**Example**: With $a = 2\cos(2\pi/14)$, the three atoms A,B and C are shown below with A being the overflow region and C the underflow. The Df map applies a (clockwise) elliptical rotation of $\theta = 2\pi/14$ to each region and the sawtooth nonlinearity *f* provides the corresponding translation – which is vertical by -2, 0 and + 2 respectively for A,B and C as shown on the right below.

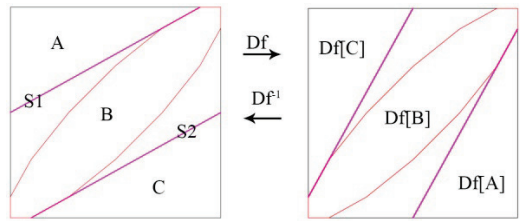

The seperatices S1 and S2 define the maximum extent of the linear center rotation, so they become extended edges of the limiting 14-gon. The bounds of the 2-torus are also extended edges of the limiting 14-gon so S1 and S2 map to these bounds. Any subset of these can serve as the initial level-0 singularity set –depending on what region is of interest. At this time there is no theory that can be used to determine a 'minimal' level-0 set for Df or for τ. It could be a point set.

The slope of S1 (or S2) determines the number of 'steps' involved in a rational rotation and of course each 'step' is a rotation by θ. Our convention (described in detail below) is to define the number of steps using the (clockwise) difference between S1 (or S2) and the 'next' bounding edge – so the example above is step-1. This matches the 'twist' ρ which is 1/N. Below is the same map in Euclidean 'rectified' space where the Df elliptical rotation becomes a true rotation and the translations have magnitude 2/sinθ to match the edges of the bounding rhombus. In the rectified plots shown here, we have added a rotation by π/2 to make them match traditional Tangent map space. This means that the displacements of A and C are now horizontal of magnitude +2/sinθ and -2/sinθ respectively.

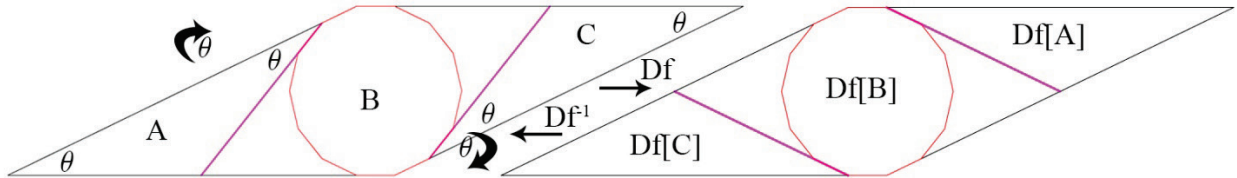

Therefore on each iteration of the web, the seperatices are rotated by θ. In this example, θ was the exterior angle of N = 14, so the resulting web is step-1 which is consistent with the web generation algorithms used for the Tangent map.

**Example**: Below are the first few iterations of the Df web in magenta and the rectified web in blue for $a = 2\cos(2\pi/14)$. This web is obtained by mapping the two magenta seperatices under Df. The last iteration is the level 100 web. The Mathematica code for level-100 is given below.

**S1=Table[{x,(x+1)/a},{x,-1,a-1,.001}]; S2=Table[{x,(x-1)/a},{x,1-a,1,.001}]; S0=Join[S1,S2];
DfWeb = Flatten[Table[NestList[Df,S0[[k]],100],{k,1,Length[S0]}],1];
TrWeb = DfToTr[DfWeb]; Graphics[{AbsolutePointSize[1.0], Magenta, Point[DfWeb],
Blue, Point[TrWeb]}].**

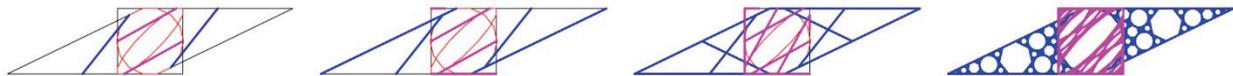

This web was generated using ρ = 1/14, but 2/14 and 3/14 are also valid parameters and they generate Df webs with steps of two and three. Both of these are significant because the step-2 web corresponds to N = 7 at θ = 2π/7 and the maximal step-3 web generates a subset of the First Family for N = 14.

It is possible to specify an 'odd' rotation like ρ = 1/7 using θ = 2π/7 but this will not generate the web for N = 7 because the symmetry of the Df map demands that the opposite sides of the limiting polygon be parallel. This is illustrated on the left below. The linear map does yield heptagons but they are not compatible with the seperatices so the limiting linear polygon is a regular 14-gon as shown in the second plot. The resulting web is generated at step-2 so it does not match the N = 14 web. Below are iterations 0,6 and 100 followed by the rectified web.

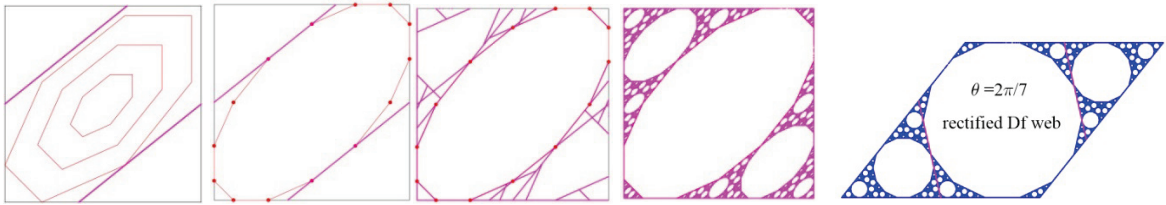

The step-2 relationship between S1 (or S2) and the 'next' edge of the rhombus can be seen in the blue rectified web on the right - where we have reproduced the image of S1 and S2 in magenta.

Note that there is also a step-3 relationship using the 'retrograde' counter-clockwise rotation, so the step-2 web based on $\tau$ (and Df) is equivalent to a step-3 web using $\tau^{-1}$ (and $Df^{-1}$). In the limit these webs will be conjugate – but our default will be clockwise.

In the case of N = 7 and N = 5, the resulting step-2 webs are recognizable subsets of the traditional Tangent map webs, but for the most part, these attempts at odd rotations are only distantly related to the desired web. As indicated earlier, N = 14 also supports a step-3 web. This is discussed in the Df Theorem below.

In general the Df map defines a one-parameter family of maps and that parameter could be the winding number (rotation number) $\rho$ or the angle $\theta = 2\pi\rho$. We will define a rotation to be 'rational' iff $\rho$ is rational. Below are some examples of rational rotations.

| $\rho$ (winding #) | 0 | 1/14 | 1/12 | 1/10 | 1/8 | 1/5 | 3/14 | 1/4 | 1/2 |
|---|---|---|---|---|---|---|---|---|---|
| Rotation $\theta$ | 0 | $2\pi/14$ | $2\pi/12$ | $2\pi/10$ | $2\pi/8$ | $2\pi/5$ | $2\pi(3/14)$ | $2\pi/4$ | $2\pi/2$ |
| $2\cos\theta$ | 2 | 1.80194 | $\sqrt{3}$ | $(\sqrt{5}+1)/2$ | $\sqrt{2}$ | $(\sqrt{5}-1)/2$ | 0.44504 | 0 | -2 |
| polygon | circle | N = 14 | N = 12 | N = 10 (& N = 5) | N = 8 | N = 10 (short) | N = 14 (short) | N = 4 | line |

By default we restrict $\theta \in (0, \pi/2]$ so $\rho \in (0, 1/4]$. This covers the full range of regular polygons since the odd cases are included in the twice-odd cases. The highlighted region is what we call the 'quadratic' range because it includes the four 'quadratic' regular polygons- N = 5, 8, 10 & 12. Below is a series of plots showing the dynamics in this range as $a = 2\cos\theta$ increases from $\sqrt{2}$ (N = 8) to $\sqrt{3}$ (N = 12) with N = 10 in between at $(\sqrt{5}+1)/2$. These are all 'rational' rotations because the winding numbers are decrements of 1/120- starting with N = 8 at $\theta = 2\pi(15/120)$.

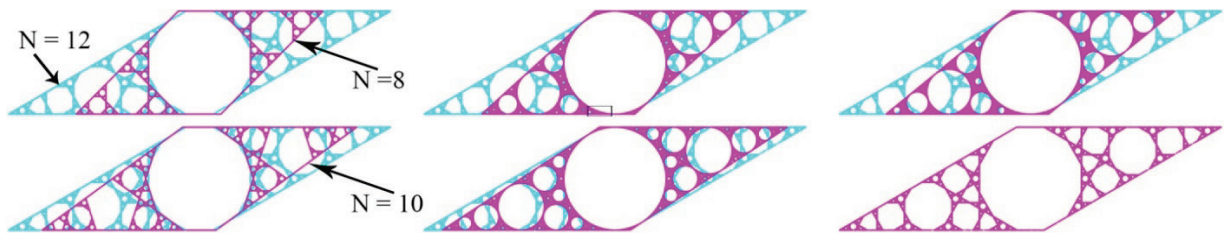

Therefore the 2nd magenta plot above has $\theta = 2\pi(14/120)$. As a twist map it is period 60 so it corresponds to a Df map with N = 60 and a 'step-7' web as shown in the detail below. This detail matches the small rectangle in the plot above.

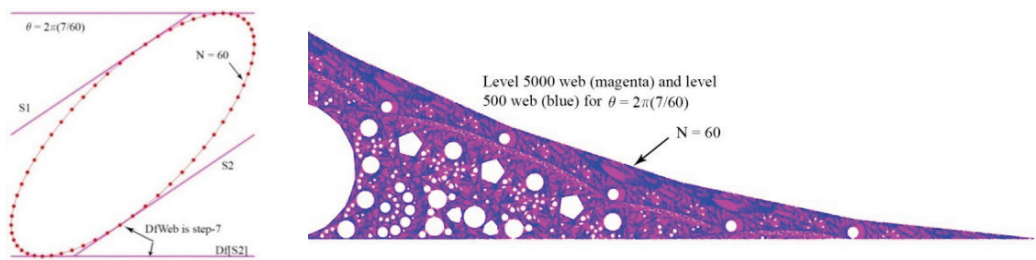

N = 60 has step-k webs ranging from 1 to 15, but the step-15 case reduces to N = 4 with null web, so the practical range extends to 14. It is not clear which of these webs are related to the traditional step-1 web - but the theorem below says that the 'maximal' 14-step web is special in that it reproduces part of the 'in-sito' dynamics of the central S[28] tile of N = 60.

**Lemma:** For the Df map with parameter $\theta = 2\pi\rho$, when the 'twist' $\rho = p/q$ with p and q relatively prime integers and $p/q \leq \frac{1}{4}$, then the resulting Df web is conjugate to that of a regular q-gon with a step-p web when q is even or a regular 2q-gon with a step-2p web when q is odd.

Since $\rho \leq \frac{1}{4}$, it follows that the maximal step size for any regular N gon is $\lfloor N/4 \rfloor$ - but as noted above, when N is twice-even this reduces to the null case of N = 4, so we will define the 'maximal' step size to be N/4-1 when N is divisible by 4. Most of these 'rational' webs are unrecognizable but there are two special cases which are described in the Theorem below.

**Theorem (Df map for Regular Polygons):** For a regular polygon N = 2k, the 'effective' range of rotation values for the Df map with $\theta = 2\pi\rho$ are $\rho = 1/N$ to $\lfloor k/2 \rfloor/N$ for k odd or $(k/2-1)/N$ for k even. Therefore the Df webs range from step-1 to step-$\lfloor k/2 \rfloor$ or step-$(k/2-1)$. These minimal and 'maximal' step values have webs which are described below:

(i) When $\rho = 1/N$ the Df web will be step-1 and this web is locally conjugate to the Tangent map web for N.

(ii) When k is odd the web in (i) will have a central S[k-2] tile which will be a regular N/2-gon. The Twice-Odd Lemma implies that the local web of this polygon will be conjugate to the Tangent map web for N/2, so the web for N includes a scaled copy of the web for N/2. When $\rho = \lfloor k/2 \rfloor/N$, the Df web will have maximal step-$\lfloor k/2 \rfloor$ and the rectified web will yield a shortened version of the step-1 rhombus from (i). This secondary rhombus will be locally identical to the step-1 Tangent map rhombus but contain just the first $\lfloor k/2 \rfloor$ tiles of the First Family of N. This is known as the Inner Family or 'short family'.

(iii) When k is even, the web in (i) will again have a central S[k-2] tile but now this tile will be a regular N-gon instead of a N/2-gon. Since N is a now a multiple of 4, the 'maximal' k/2-1 step will be odd iff N is a multiple of 8. In this case $\rho = (k/2 -1)/N$ will have a Df web which is conjugate to the local web for the central S[k-2] tile.

To see why the short webs from (ii) and (iii) are faithful to the tangent map web, note that in both cases the gap between consecutive Df boundary edges is $\lfloor k/2 \rfloor+1$. In case (ii) this implies that the maximal step-$\lfloor k/2 \rfloor$ is retrograde step-1 as in the traditional Tangent map web. In case (iii) the 'maximal' k/2-1 step is retrograde step-2 and this matches the evolution of S[k-2] in the Tangent map web when N is twice even. This evolution was described in section 2. S[k-2] is generated by two interwoven step-2 cycles which define the even and odd edges of this N-gon. When N is twice-odd, these cycles collapse to a single N/2-gon and yield the Twice-Odd Lemma.

**Examples:** The first few terms in the 'twice-odd' series of ρ are 1/6, 2/10, 3/14, 4/18 so the numerators are the positive integers and the denominators increase by 4. On the resulting rectified web we have included the image of the seperatices S1 and S2 in magenta. Note that these webs are all (retrograde) step-1 on the shortened rhombii.

| N = 10 : ρ = 1/10, a = 2Cos[2π/10] | ρ = 2/10, a = 2Cos[2π/5] |
|---|---|
| 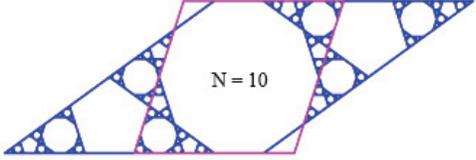 | 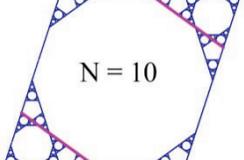 |
| N = 14 : ρ = 1/14, a = 2Cos[2π/14] | ρ = 3/14, a = 2Cos[3π/7] |
| 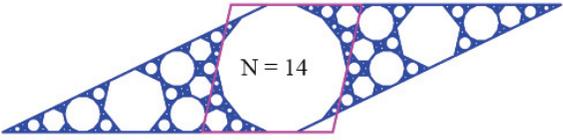 | 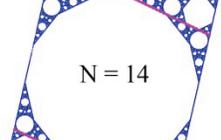 |
| N = 18 : ρ = 1/18, a = 2Cos[2π/18] | ρ = 4/18, a = 2Cos[4π/9] |
| 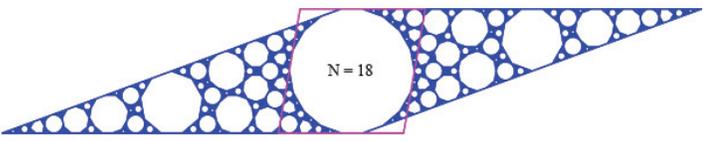 | 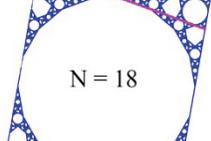 |

The first few terms of the 'twice-even' series of ρ are 1/8, 3/16, 5/24, 7/32, 9/40 so the numerators are the odd integers and the denominators increase by 8. As above, we have included the images of S1 and S2 in magenta, and here they are (retrograde) step-2.

| N = 16 : ρ = 1/16, a = 2Cos[2π/16] | ρ = 3/16, a = 2Cos[3π/8] |
|---|---|
| 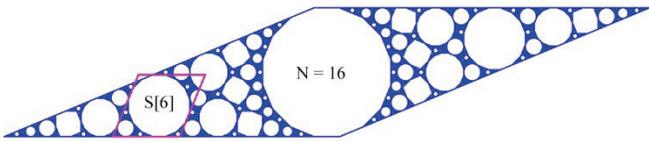 | 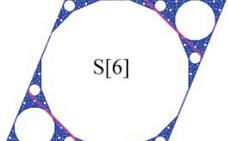 |
| N = 24 : ρ = 1/24, a = 2Cos[2π/24] | ρ = 5/24, a = 2Cos[5π/12] |
| 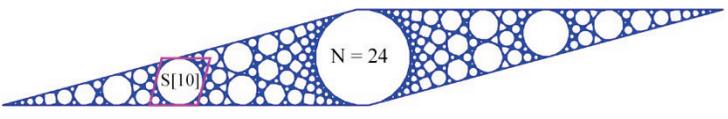 | 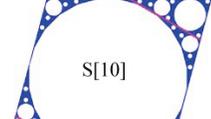 |
| N = 32 : ρ = 1/32, a = 2Cos[2π/32] | ρ = 7/32, a = 2Cos[7π/16] |
| 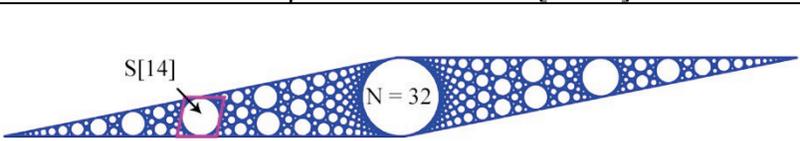 | 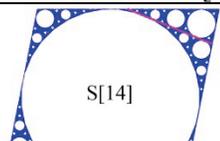 |

Of course 'most' parameter values correspond to irrational rotations and a number of studies have centered on these cases. In [A] (2001), P. Ashwin discusses the case of $2\cos(\theta) = 3/2$ which lies in the 'quadratic range' between N = 10 and N = 8. He conjectures that the residual set has positive Lebesgue measure and that this measure changes continuously with $\theta$ so that a full measure of parameters yield positive Lebesgue measure. In the grid from earlier, $\theta$ = ArcCos[3/4] is closest to $\theta = 2\pi(14/120)$ which is the case studied in some detail here. For comparison, below is a scan of this irrational rotation showing a region comparable to $\theta = 2\pi(14/120)$ shown earlier.

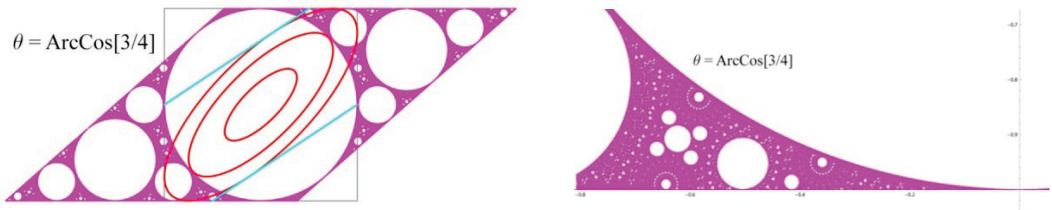

Even though the web is dense, there appears to be structure on all scales and it is unlikely that the residual set has positive Lebesgue measure. However we agree with the second part of Ashwin's conjecture- that the scaling appears to vary continuously with $\theta$. It would also be natural to expect that 'most' parameter values yield a multi-fractal spectrum of dimensions and this spectrum should vary continuously with $\theta$ and $\rho$, so the maximal Hausdorff dimension of 2 would be expected to occur in the limit of $\rho \to 0$. Indeed the case of small $\rho$ is very complex as it mimics the case of large N.

**Can the Df map be reduced to a piecewise isometry on three triangles ?**

Most piecewise isometries are formulated using triangular atoms. For the Df map, the central B region is a hexagon, but it may be possible to use the symmetry between A and C to reduce the 'essential dynamics' to a three-triangle 'dart' as suggested by Broussard and Poggiaspalla [BP]. One possible 'dart' for N = 14 is outlined on the left below. Of the 15 different configurations for three triangles, the dynamics of this dart corresponds to mapping a 'type 5' to a 'type 12' as shown on the right below.

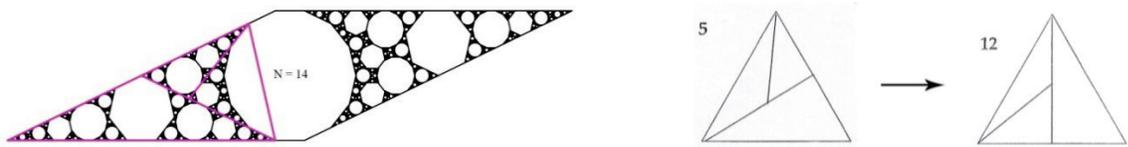

Their computer analysis found 96 possible 'solutions' for this mapping, but only six corresponded to isolated piecewise isometries, and of these, the one that is most promising for regular polygons is their solution 35 which for parameter $\pi/n$ reproduces 'most' of the dynamics of the tangent map for a regular 2n-gon – but the missing dynamics may be an issue for larger N values.

## Appendix G – Decomposition of orbits and quasi-regular polygons

**Definition**: Suppose a regular n-gon has dihedral group $D_n$ where $C_n$ is the subgroup of n rotations. If a set S of congruent tiles has a periodic orbit with period k, then every element of $C_n(S)$ will also have a periodic orbit with period k. When S is a proper subset of $C_n(S)$ then S is called a *factor* (or *decomposition*) of a combined 'orbit' formed from $C_n(S)$.

**Example**: Below are three groups of period 3 tiles for N = 9. Together they form a combined 'orbit' of period 9. Note that all the tiles share the same step sequence, which is {3}. Each group of tiles only 'sees' one of the three regular triangles embedded in N = 9.

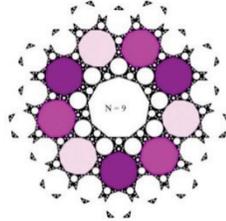

For the S[k] tiles of regular polygons, orbits factor in this fashion if and only if GCD[k,N] >1 so every orbit that can factor, will factor. This occurs because each group only 'sees' an embedded regular k-gon. For N = 9, the points in S[3] only see the embedded blue triangle .

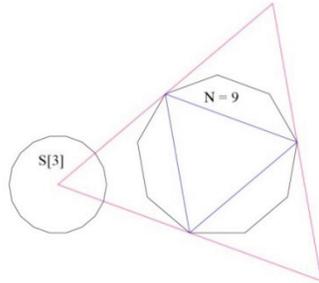

**Example**: For periodic tiles outside the 'star' region, the same principle applies where 'k' is replaced by the total number of steps in a periodic sequence. For example with N = 14, the D's in Ring 2 have step sequence {6,7,7}. Since GCD[20,14] = 2, the points in these D tiles will only 'see' one of the embedded N = 7 heptagons, so Ring 2 will decompose into two distinct orbits of period 21 each as shown below. All the 'even' rings will experience the same decomposition.

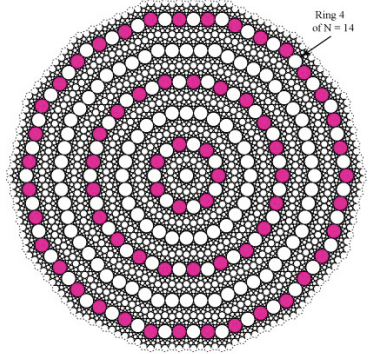

In ring 2 above, the magenta (or white) centers form what we call a 'quasi-regular' polygon. It should be clear that subsequent even rings collapse down to ring 2 (and ring 0), so this quasi-regular 21-gon is unique to N = 14. This defines a class of non-regular polygons which we call Ring2 polygons. This class will be discussed below.

Every polygon has a corresponding *factor graph* which shows the embedded regular polygons which share all of their vertices with the 'parent' polygon. Below are the factor graphs for N = 9, N = 12, N = 14 and N = 20. In some cases there are multiple congruent factors.

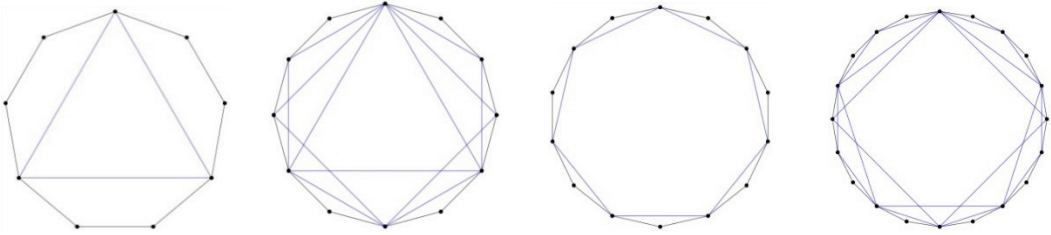

When N is regular, every embedded regular k-gon has a corresponding k-step orbit which 'sees' only the vertices of that polygon (or a rotated version) and conversely every constant k-step orbit sees just the vertices of one constituent factor. When N is prime there are no embedded regular polygons and all the k-step orbits are 'entrained' with D so the centers visit every vertex of N just once and hence the S[k] all have period N.

**Quasi-regular polygons**

Non-regular polygons can also have regular factors which play an important part in the dynamics. We will define a *quasi-regular polygon* to be a convex polygon which contains at least one regular factor. (Other authors have used the term 'quasi-regular' for different purposes.) So the regular polygons are also quasi-regular and so are the 'woven' polygons which often arise in mutations of canonical orbits. The examples below show three of the ways in which quasi-regular polygons arise: (i) from decomposed orbits of D's with N even (Ring2 polygons), (ii) from orbits of D's with N odd (Riffle polygons) and (iii) from 'mutations' of S[k] for composite N, caused by the decomposition of their orbits ('woven' polygons).

Examples:
(i) When N is even all regular polygons have rings of D's whose centers form regular polygons. As indicated in Section 6, ring k will have N·(k+1) D's and for N twice-odd the even rings always 'decompose' into two factors. For N twice-even there is no decomposition, but we will impose one in the same fashion as the twice-odd case. Starting with ring 2, these factors will be non-regular. We call them *Ring2 polygons*. (As indicated earlier, subsequent even rings collapse down to ring 2 (and ring 0)). These Ring2 polygons will have 3N/2 vertices and three regular factors with N/2 vertices each (so when N= 4, these 'regular' factors reduce to line segments.). We will call these factors A, B and C and scale the largest factor A to have radius 1. B and C will have equal radii . Below are the first four members of the Ring2 family of quasi-regular polygons with A in green, B in magenta and C in blue.

| N4Ring2 | N6Ring2 | N8Ring2 | N10Ring2 |
|---------|---------|---------|----------|
|         |         |         |          |

Note that these Ring2 polygons can be generated by taking any regular N-gon with N even and dividing the sides into 4 equal segments and then selecting the even or odd vertices. This mimics the 3N/2 process of ring 2 and of course generates both even and odd number of vertices.

Below is ring 2 for N = 6, showing the unscaled quasi-regular N6Ring2 in magenta. This is also called N9NonRegular when the context is understood. The rings are clones of N9NonRegular so they are M-type, but the key issue for invariance is their local geometry and not the M or D issue. As can be seen on the left, the M-type triangles can form rings as well as the hexagon D's. These are called type-2 rings in [VS] – although they did not cover the N-even case in that paper.

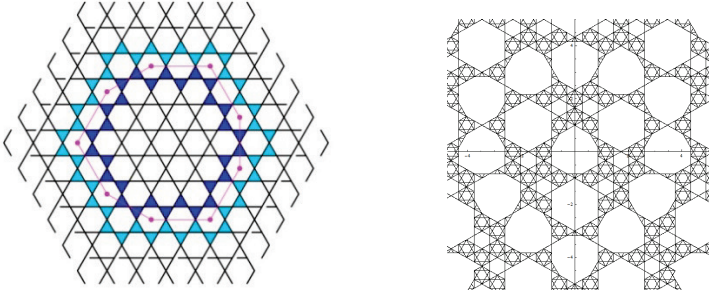

The case of N = 10 is shown below with the web for N15NonRegular again showing rings of congruent polygons playing the role of D's. See the enlargement below.

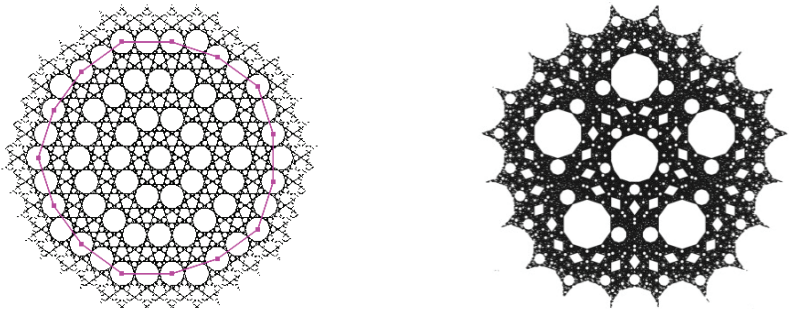

**Conjecture**: The members of the Ring2 family derived from N twice-odd have bounded dynamics.

These are the Ring2 members with an odd number of sides. The first member is N9NonRegular described above. One common element in all of the Ring2 webs is the recurrence of small-scale dynamics conjugate to that of the factors and the 'parent', so they are typically multi-fractal. Below is the level-10 (inverse) web for N15NonRegular. N21Non-Regular has similar dynamics. (Length $s_3$ is defined by the edges of the yellow regular pentagon – which technically does not exist.)

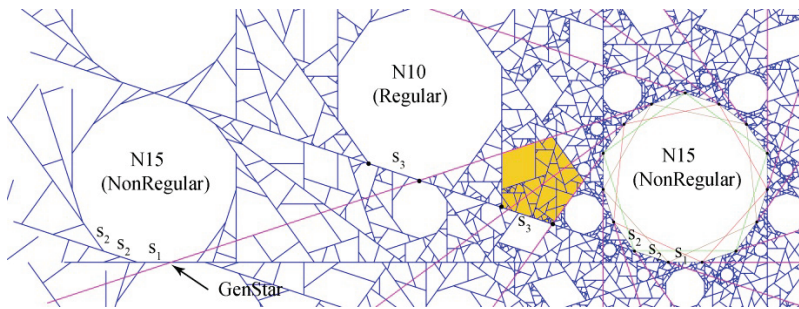

(ii) **Riffle Polygons and Semi-Regular Polygons**

The rings of 'D' tiles for N odd generate non-regular polygons directly. Below are rings 0,1 and 2 for N = 7. Ring 0 is the only regular case. Ring 1 is period 21 but colinear points reduce it to the non-regular 14-gon shown in the middle below. It has two identical N = 7's as factors, so these mimic the two congruent factors of N = 14. The dynamics appear to be bounded and this may be generic for this class of quasi-regular polygons.

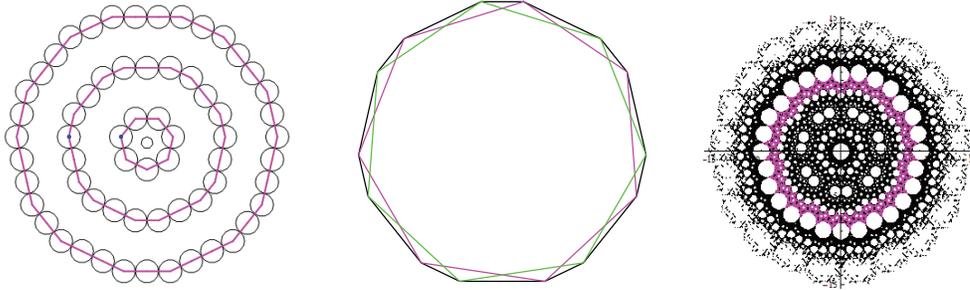

This example is special since the two factor heptagons A and B are simple rotations of each other. We call these *Riffle polygons* because the resulting quasi-regular 14-gon can be obtained from N = 14 and the two canonical N = 7 factors, A (green) and B (magenta), by first rotating B to obtain C and then using the Mathematica Riffle[A,C] command to interweave them. (The actual rotation here is an irrational factor of $\pi/7$ since it should be clear that the ratio of sides above is exactly 2. This ratio is fixed for rings obtained from N-gons for N odd.)

| N=14 is Riffle[A,B] | C = Rotate B by ≈.3383π/7 | Riffle[A,C] |
|---|---|---|
|  |  |  |

Using the rotation angle between heptagons A and B as a parameter, this yields a 'continuum' of quasi-regular Riffle polygons – all of which share similar dynamics. (However the class of all quasi-regular polygons remains a class of measure zero among all polygons.)

**Definition**: Suppose that P is a regular N-gon with N even, where vertex k is written P[[k]]. Separate the 'odd' and 'even' vertices as A = Table[P[[k]],{k, 1, N, 2}] and B = Table[P[[k]],{k, 2, N, 2}]. Then Riffle[A,B] will alternate vertices to yield P. Define C =RotationTransform [theta] [B] where the angle of rotation theta is of the form $\rho 2\pi/N$ with p ∈ [0,1]. Then the family Riffle[A,C] runs from P to A and this is called the **Riffle family** of P. When N is odd, we will define the Riffle family of N to be the same as that of 2N.

**Lemma**: Every member of the Riffle family of a regular polygon has bounded dynamics.

The proof of this Lemma is based on the fact that the large-scale web evolution is similar to the regular case so that rings of maximal D tiles will exist at all distances from the origin. For twice-even Riffle polygons, these rings mimic closely the regular case, but the twice-odd regular polygons have alternating even and odd parity rings (see Section 6), and this alternation carries over to the Riffle family and interferes with the formation of the even rings.

Below is an example using N = 14 and $\rho$ = .25. The extended forward and trailing edges of this N14 Riffle form the level-0 inverse and forward webs shown the right below in blue and magenta. This level-0 web tells the whole story because for any value of $\rho \in [0,1)$ the angles between adjacent vertices are unchanged from the regular case, so all the members of this Riffle family have the same large-scale web structure as N = 14 (except for $\rho$ = 1 which is N = 7).

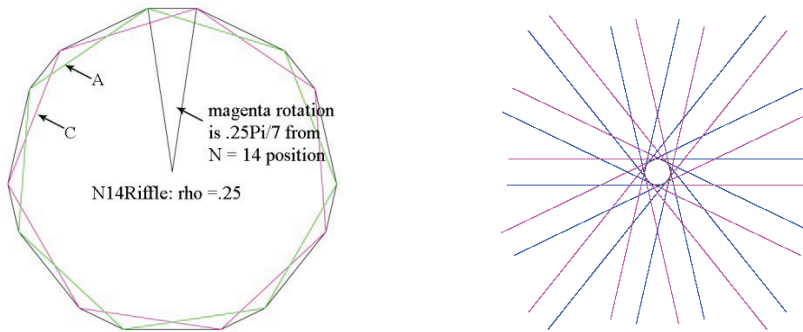

In the level-0 web, the region between two consecutive blue forward edges forms one of the 14 congruent subdomains of $\tau^{-1}$. We will look at one such subdomain below. Each of these subdomains is divided into 6 regions by the magenta level-0 trailing edges. The corresponding intersection points are called the 'star' points and the regions are the 'step' regions discussed in the Star Point Theorem.

All the points in a given region must share the same 'next' vertex so they map together for at least one iteration of $\tau$. The Star Point Theorem claims that for a regular polygon, this selection process can be repeated so that there is always a canonical tile in each region where all points have a constant step-k orbit. Some of this structure carries over to Riffles as well, but our interest here is in the unbounded 6$^{th}$ region defined by the (level-0) GenStar point.

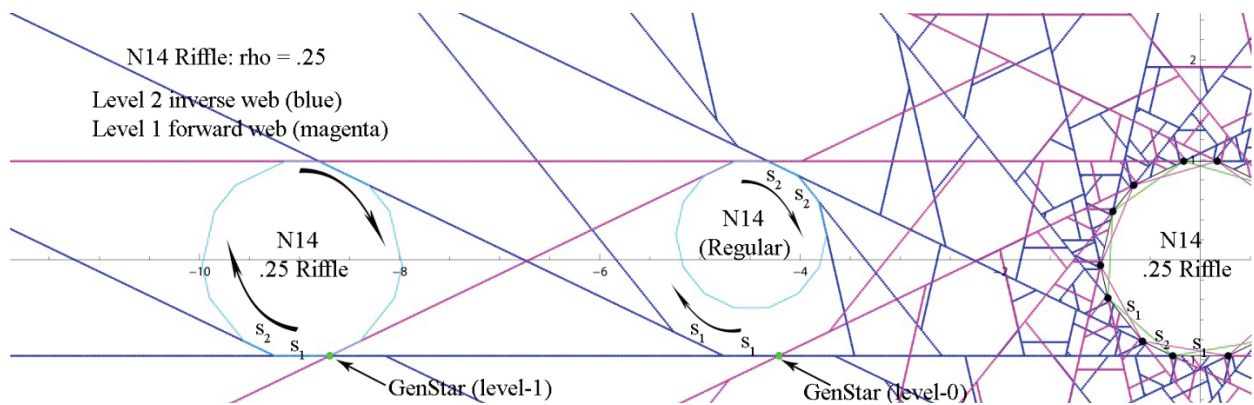

The evolution of the GenStar point for N = 14 is discussed in Sections 2 and 6 where we applied the 'swap domain and range' concept to show that all star points follow a similar 'slide' and

'rotate' algorithm that would form N= 14 starting from a single vertex. For the tangent map, τ, the 'slide' is a shear caused by the side of N = 14 and the rotation depends on the step region so it is of the form kφ where φ is the exterior angle of N = 14. For the GenStar region, k = 1, so the slide and rotate algorithm applied recursively will generate an identical copy of N = 14.

For the N14Riffle given here, the angles are unchanged from the N = 14 case, but now there are two competing shears. Recall that the shears act outwards on the lower edges and inwards on the top edges so the points in the step-6 region evolve with two different side lengths – which means that the shorter side will prevail as shown above. This evolution is matched by the subdomain region above because now the 'bottom' shear is $s_2$. These regular tiles form a ring with gaps, so the region inside is not invariant. We call this Ring 0.

The magenta web shown above has a second GenStar point because this is a level-1 web and each iteration acts on the previous to generates a new GenStar point with the same relative displacement. For the regular case these GenStar points generate unbroken rings of maximal D tiles which guarantee boundedness.

Dynamically these GenStar regions are one step apart, so the points in the second GenStar region will have step sequence {6,7}. This is an 'odd-step' region as opposed to a constant {6} for the previous region. A step-6 orbit will only experience $s_1$ or $s_2$ shears and in fact the points in the regular tile only 'see' the green 'A' heptagon. By contrast the points with {6,7} will visit all the vertices and reproduce an identical copy of the N14Riffle. Therefore the rings come in pairs, alternating even and odd step. In the regular case this odd-even alternation still occurs but it has little effect because all the sides of N = 14 are equal. However it does imply that even rings must 'decompose' because a given D only 'sees' one of the two embedded heptagons.

A portion of the level-20 combined web is shown below in blue. By convention, the periods given here are based on the center points, so the even rings will have two periods. The N14 regular tiles generate perfect local families so as ρ gets small, these 'seeds' generate multiple rings of tiny N = 14 D's.

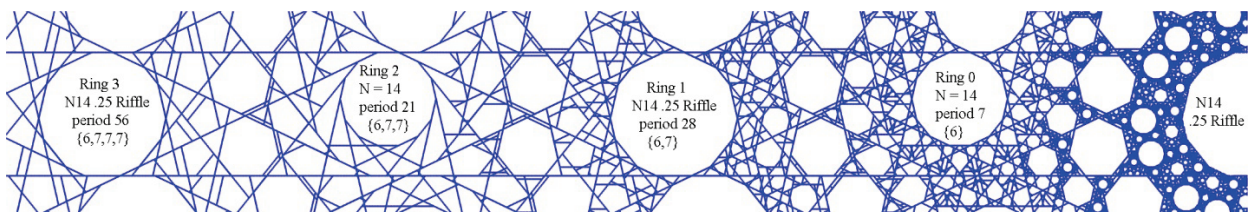

**Semi-Regular polygons**

These 'Riffle' polygons are also known as semi-regular polygons – although there is no formal definition of 'semi-regular' when applied to polygons. The definition given below is based on 'historical usage' and the recent work of Patrick Hooper [H1] and Richard Schwartz [S3] on polytope exchange transformations - PETS.

**Definition**: For N-even, an N-gon is **semi-regular** iff it has dihedral symmetry group $D_{N/2}$.

Since the regular polygons have dihedral symmetry group $D_N$, they can be regarded as special cases of semi-regular polygons. It appears that the class of semi-regular polygons is the same as the class of Riffle polygons.

**Example**: The graphic below is an overview of the semi-regular 'family' for N = 8 using Riffle rotations. The rotation in each case is $\rho 2\pi/8$, with limiting case N = 4 at $\rho = 1$ when the rotating embedded Red square, coincides with the fixed Green. The blue inverse webs are level-100 and the magenta forward webs are level-1. These coarse webs do little justice to the intricate structure. In [S3], Schwartz uses the equivalence between semi-regular octagons and octagonal PETS to obtain a complete characterization of the dynamics of this family. The most 'interesting' cases are often highly rational rotations. Clearly $\rho = 1/2$ gives a partial alignment with N = 16 and $\rho = 1/3$ gives a similar alignment with N = 12. At this point no one knows the relationship between the scaling parameters of N = 8, N = 12 and N = 16 so these semi-regular 'mixed' cases may provide some insight into how these parameters are related. The alignment with N = 12 apparently yields the highest value of the fractal dimension – which (surprisingly) exceeds that of N = 8. As indicated earlier, the algebraic relationship between N = 8 with $2\cos\theta = \sqrt{2}$ and N = 12 with $2\cos\theta = \sqrt{3}$ encompasses the full quadratic range. The fourth–order N = 24 alignment at $\rho = 2/3$ is also algebraically interesting.

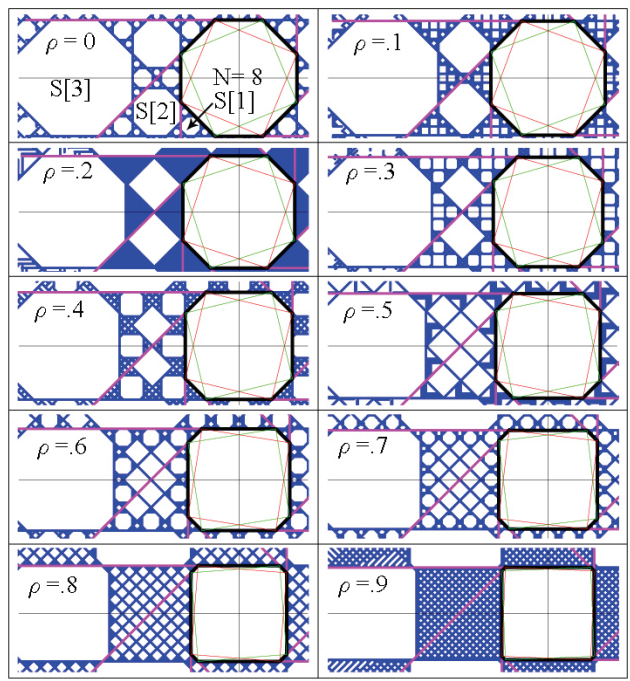

(iii) There are many quasi-regular polygons where the (regular) factors are identical but they have different radii as in the Ring2 polygons. We call these 'woven' polygons. These woven polygons also arise as 'mutations' of canonical S[k] or DS[k] tiles. Whenever the orbit of an S[k] or DS[k] decomposes into smaller periods, the dynamics are constrained to lie on the corresponding embedded factor polygon so the web may form an incomplete version of the regular polygon with fewer edges than the corresponding S[k] or DS[k]. The centers are unchanged from the regular case because these are a function of the level-0 star points.

For example N = 20 has a mutated S[2] consisting of two interwoven regular pentagons. The first occurrence of a mutated S[k] is the S[3] tile of N = 9 as described in Section 2. Because of the embedded N = 3, S[3] becomes a quasi-regular dodecagon consisting of two interwoven regular hexagons at slightly different radii as shown on the left below. On the right is the web which results when this mutated S[3] is the generator. In general it is not clear what is the relationship between the *in situ* dynamics and *in vitro* dynamics of any polygon. (See (xi) in Appendix A.)

The large-scale web for S[3] shows rings of D's which are identical to the generator, but the rings have gaps so the 'star' region is not invariant. The gaps contain seeds for secondary rings which become dominant in a cyclic fashion – so each of the embedded regular hexagons exerts influence which varies periodically. This is typical of the 'woven' polygons studied. In most cases there is no sign of large-scale invariance so they may support unbounded orbits. The Ring2 polygons derived from N twice-even seem to have dynamics similar to the traditional 'woven' polygons, but the Ring2 polygons derived from N twice-odd, may be an exception.

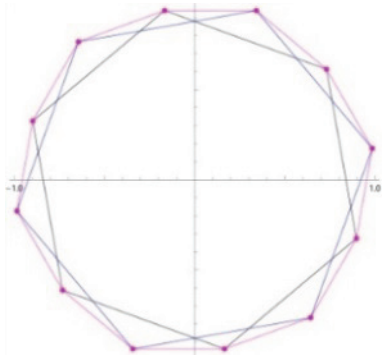
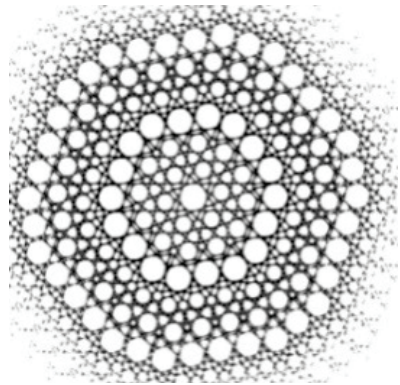

# Appendix H - using the Digital Filter map to reconcile decomposed orbits

When N is even it is always possible to subdivide the inner star into N/2 congruent rhombi as shown below for N = 16. Note that the overlap is identical to the region not covered so each rhombus contains exactly 1/8 of the tiles.

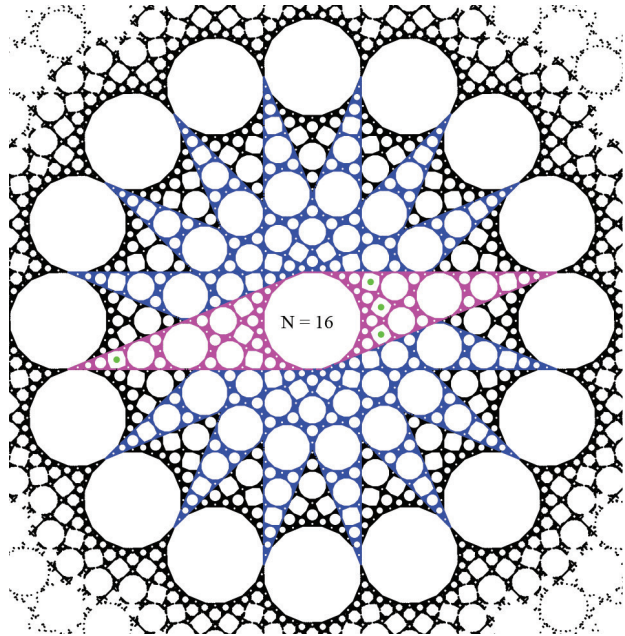

The Df map performs this partition naturally and it also mimics the reflective symmetry of each rhombi, so if the green point shown above has period 4 in the Df map, then there will be another 4 congruent regions in the magenta rhombus. They are easy to spot here. This makes a total of 8 congruent regions in the magenta tile and there are 8 such tiles in the star region for a total count of 64 tiles.

The green point chosen here is the center of LS[4] for N = 16 and it has period 12 back in tangent space. There are 4 groups of 12 in the outer star region, for a total of 48 tiles. The Df map unites these 48 with the 16 congruent S[4] tiles from the inner star. In this way the Df map can be used to reconcile orbits. LS[4] and S[4] are symmetric tiles for N = 16 as shown below, yet they have very different periods (12 and 16). The Df map sees them as part of one united orbit because the 16 D's in the first ring act as one. This correspondence between Df orbits and the number of congruent tiles in the star region is called the 2kN Lemma.

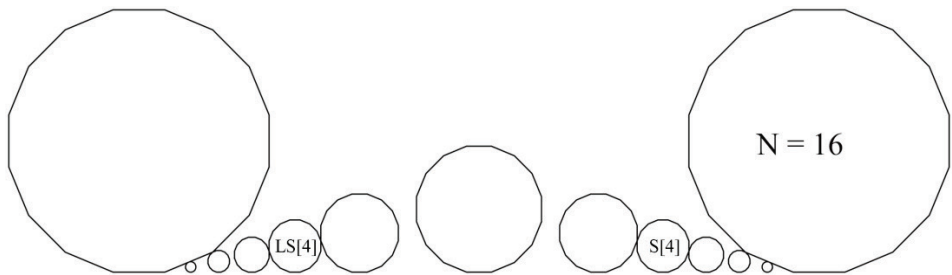

**The 2kN Lemma**: For the Df map with θ = 2Pi/2N, every periodic tile in Df with period k >2 accounts for 2k congruent tiles in Df space and 2kN congruent tiles in the star region of τ. In this star region, the centers of these tiles may have different periods but the sum of the periods must be 2kN. (If a tile in Df space has period 2, then the center point p must map to –p so there are only 2 congruent tiles in Df space and therefore only 2N congruent tiles in the star region.)

Example: Note that for regular 2N-gons, the Df periods of the S[k] are just the number of congruent tiles to the right and left of the central tile as shown below for N = 11 using θ = 2π/22. These periods are consecutive integers from 2 to 10. In general the Df period of S[k] is N-k.

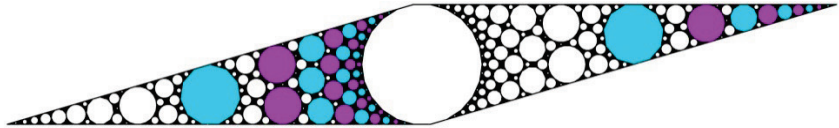

This puts the central S[9] in a special position as the only tile which maps to 'itself', so the total number of congruent S[9] tiles in τ- space is 2N, but all the other Df periods get doubled by symmetry. For example the number of tiles in τ-space congruent to S[5] is 2*6*11. Note that these tiles are typically partitioned into multiple periodic orbits, so it is not trivial to derive periods from these counts. However there are simple formulas for the periods of the members of the First Family for prime N. The next example illustrates the issues which arise for composite N.

Example: The 2kN Lemma can be used to help organize the family periods for N = 16. Use Df with θ = 2π/16 (so N = 8 in the Lemma). Earlier we found that the Df period of S[4] was 4 and hence the number of congruent tiles in τ-space was 2*4*8. In this case it was possible to account for these 64 tile as 12*4 + 4*4 where LS[4] had period 12 and S[4] had period 4. The orbit of S[4] would normally be 16 but this decomposes into 4 groups of 4. This decomposition complicates the periodicity issue for 2N-gons – and that includes twice-odd and twice even.

Below are the periods of the N= 16 first family where S[6] = LS[6] occupies the central position. Note that it is period 8 because of decomposition. However the Df count for the number of congruent S[6] tiles is correct at 2*8.

| Tile | S[1] | S[2] | S[3] | S[4] | S[5] | S[6] | LS[1] | LS[2] | LS[3] | LS[4] | LS[5] | LS[6] |
|---|---|---|---|---|---|---|---|---|---|---|---|---|
| Period | 16 | 8 | 16 | 4 | 16 | 8 | 96 | 40 | 64 | 12 | 32 | 8 |

The 2kN Lemma can make some sense out of a table like this by uniting congruent tiles such as S[k] and LS[k] unto a combined count using the known Df periods of the S[k]. This is shown below.

| Tiles | S[1]& LS[1] | S[2]& LS[2] | S[3]& LS[3] | S[4]& LS[4] | S[5]& LS[5] | S[6]& LS[6] |
|---|---|---|---|---|---|---|
| Count | 2*7*8 | 2*6*8 | 2*5*8 | 2*4*8 | 2*3*8 | 2*8 |

**Links:**

(i) The web site of the author at DynamicsOfPolygons.org is devoted to the Tangent map and related maps from the perspective of a non-professional. There are three versions of this file available there for download. They differ only in the embedded graphics.

RegularPolygonsBasic.pdf (about 12Mb – no embedded graphics)
RegularPolygonsWord.doc (about 30Mb – no embedded graphics)
RegularPolygonsArxiv.pdf  (the published document – about 39Mb at arxiv:1311.6763)

Just click on the link to PDFs. These should be downloaded (right click) because most on-line viewers cannot handle large documents.

More specific links are given below.

(ii) Outer billiards, digital filters and kicked Hamiltonians: arXiv:1206.5223

(iii) A chronology of the Tangent map

(iv) Assorted images.

(v) Mathematica notebooks.  The notebooks cover the four possible cases: Nodd, NTwiceOdd, NTwiceEven and Nonregular. (The Digital Filter map is included in the N-even notebooks, but it is independent of the notebooks.)

(vi) For more on pinwheel maps and related projections go to PDFs and select Projections.pdf